\documentclass[final,onefignum,onetabnum]{siamart220329}

\usepackage{braket,amsfonts,amsmath,amssymb}

\usepackage{array}

\usepackage[caption=false]{subfig}
\usepackage{pgfplots}

\newsiamthm{claim}{Claim}
\newsiamthm{conjecture}{Conjecture}
\newsiamremark{remark}{Remark}
\newsiamremark{hypothesis}{Hypothesis}
\crefname{hypothesis}{Hypothesis}{Hypotheses}

\usepackage{algorithmic}

\usepackage{graphicx,epstopdf}

\Crefname{ALC@unique}{Line}{Lines}

\usepackage{amsopn}

\usepackage{xspace}
\usepackage{bold-extra}
\usepackage[most]{tcolorbox}

\colorlet{texcscolor}{blue!50!black}
\colorlet{texemcolor}{red!70!black}
\colorlet{texpreamble}{red!70!black}
\colorlet{codebackground}{black!25!white!25}

\lstdefinestyle{siamlatex}{%
  style=tcblatex,
  texcsstyle=*\color{texcscolor},
  texcsstyle=[2]\color{texemcolor},
  keywordstyle=[2]\color{texemcolor},
  moretexcs={cref,Cref,maketitle,mathcal,text,headers,email,url},
}

\tcbset{%
  colframe=black!75!white!75,
  coltitle=white,
  colback=codebackground, 
  colbacklower=white, 
  fonttitle=\bfseries,
  arc=0pt,outer arc=0pt,
  top=1pt,bottom=1pt,left=1mm,right=1mm,middle=1mm,boxsep=1mm,
  leftrule=0.3mm,rightrule=0.3mm,toprule=0.3mm,bottomrule=0.3mm,
  listing options={style=siamlatex}
}

\newtcblisting[use counter=example]{example}[2][]{%
  title={Example~\thetcbcounter: #2},#1}

\newtcbinputlisting[use counter=example]{\examplefile}[3][]{%
  title={Example~\thetcbcounter: #2},listing file={#3},#1}

\DeclareTotalTCBox{\code}{ v O{} }
{ 
  fontupper=\ttfamily\color{black},
  nobeforeafter,
  tcbox raise base,
  colback=codebackground,colframe=white,
  top=0pt,bottom=0pt,left=0mm,right=0mm,
  leftrule=0pt,rightrule=0pt,toprule=0mm,bottomrule=0mm,
  boxsep=0.5mm,
  #2}{#1}

\patchcmd\newpage{\vfil}{}{}{}
\begin{tcbverbatimwrite}{tmp_\jobname_header.tex}
\title{Coupled cluster theory: Towards an algebraic geometry Formulation\thanks{Submitted to the editors \today.
\funding{This work was partially supported by the Air Force
Office of Scientific Research under award number FA9550-18-1-0095 and by the Simons Targeted Grants in Mathematics and Physical Sciences on Moir\'e Materials Magic (F.M.F.)}}}

\author{Fabian M. Faulstich\thanks{Department of Mathematics, Rensselaer Polytechnic Institute, Troy, NY (\email{faulsf@rpi.edu})}
\and Mathias Oster\thanks{IGPM, RWTH Aachen, Germany (\email{oster@igpm.rwth-aachen.de})}}

\headers{CC Theory: Towards an Algebraic Geometry Formulation}{Fabian M. Faulstich \& Mathias Oster}
\end{tcbverbatimwrite}
\title{Coupled cluster theory: Towards an algebraic geometry Formulation\thanks{Submitted to the editors \today.
\funding{This work was partially supported by the Air Force
Office of Scientific Research under award number FA9550-18-1-0095 and by the Simons Targeted Grants in Mathematics and Physical Sciences on Moir\'e Materials Magic (F.M.F.)}}}

\author{Fabian M. Faulstich\thanks{Department of Mathematics, Rensselaer Polytechnic Institute, Troy, NY (\email{faulsf@rpi.edu})}
\and Mathias Oster\thanks{IGPM, RWTH Aachen, Germany (\email{oster@igpm.rwth-aachen.de})}}

\headers{CC Theory: Towards an Algebraic Geometry Formulation}{Fabian M. Faulstich \& Mathias Oster}

\ifpdf
\hypersetup{ pdftitle={Title} }
\fi

\newcommand{\dd}{\rm d}

\usepackage{import}
\definecolor{airforceblue}{rgb}{0.36, 0.54, 0.66}

\usepackage{graphicx,epstopdf} 
\usepackage[caption=false]{subfig} 

\begin{document}
\maketitle

\begin{tcbverbatimwrite}{tmp_\jobname_abstract.tex}
\begin{abstract}
Coupled cluster theory produced arguably the most widely used high-accuracy computational quantum chemistry methods.
Despite the approach's overall great computational success, its mathematical understanding is so far limited to results within the realm of functional analysis.
The coupled cluster amplitudes, which are the targeted objects in coupled cluster theory, correspond to solutions to the coupled cluster equations, which is a system of polynomial equations of at most degree four.
The high dimensionality of the electronic Schr\"odinger equation and the non-linearity of the coupled cluster ansatz have so far stalled a formal analysis of this polynomial system.
In this article, we present algebraic investigations that shed light on the coupled cluster equations and the root structure of this ansatz. 
This is of importance for the {\it a posteriori} evaluation of coupled cluster calculations. 
To that end, we investigate the root structure by means of Newton polytopes. 
We derive a general v-description, which is subsequently turned into an h-description for explicit examples.
This perspective reveals an apparent connection between {\it Pauli's exclusion principle} and the geometrical structure of the Newton polytopes. 
We also propose an alternative characterization of the coupled cluster equations projected onto singles and doubles as cubic polynomials on an algebraic variety with certain sparsity patterns.
Moreover, we provide numerical simulations of two computationally tractable systems, namely, the two electrons in four spin-orbitals system and the three electrons in six spin-orbitals system. These simulations provide novel insight into the root structure of the coupled cluster solutions when the coupled cluster ansatz is truncated. 
\end{abstract}

\begin{keywords}
homotopy continuation, Newton polytope, Bernstein--Khovanskii--Kushnirenko theorem, coupled-cluster theory, quantum many-body problem, Schr\"odinger equation 
\end{keywords}

\begin{MSCcodes}
12D10, 14Q20 , 90C53, 81-08, 81-10  
\end{MSCcodes}
\end{tcbverbatimwrite}
\begin{abstract}
Coupled cluster theory produced arguably the most widely used high-accuracy computational quantum chemistry methods.
Despite the approach's overall great computational success, its mathematical understanding is so far limited to results within the realm of functional analysis.
The coupled cluster amplitudes, which are the targeted objects in coupled cluster theory, correspond to solutions to the coupled cluster equations, which is a system of polynomial equations of at most degree four.
The high dimensionality of the electronic Schr\"odinger equation and the non-linearity of the coupled cluster ansatz have so far stalled a formal analysis of this polynomial system.
In this article, we present algebraic investigations that shed light on the coupled cluster equations and the root structure of this ansatz.
This is of importance for the {\it a posteriori} evaluation of coupled cluster calculations.
To that end, we investigate the root structure by means of Newton polytopes.
We derive a general v-description, which is subsequently turned into an h-description for explicit examples.
This perspective reveals an apparent connection between {\it Pauli's exclusion principle} and the geometrical structure of the Newton polytopes.
We also propose an alternative characterization of the coupled cluster equations projected onto singles and doubles as cubic polynomials on an algebraic variety with certain sparsity patterns.
Moreover, we provide numerical simulations of two computationally tractable systems, namely, the two electrons in four spin-orbitals system and the three electrons in six spin-orbitals system. These simulations provide novel insight into the root structure of the coupled cluster solutions when the coupled cluster ansatz is truncated.
\end{abstract}

\begin{keywords}
homotopy continuation, Newton polytope, Bernstein--Khovanskii--Kushnirenko theorem, coupled-cluster theory, quantum many-body problem, Schr\"odinger equation
\end{keywords}

\begin{MSCcodes}
12D10, 14Q20 , 90C53, 81-08, 81-10
\end{MSCcodes}


\section{Introduction}
Electronic structure theory models a quantum mechanical system of electrons moving in an exterior potential, mostly given by the electrical field of (clamped) nuclei. 
This model covers most of the quantum physical effects of chemistry and bulk crystals. 
{\it Ab initio} electronic structure calculations provide numerical simulations based on the first principles of quantum mechanics. 
The persistent difficulty in the study of a quantum many-particle system is the high dimensionality. The object of interest, i.e., the system's (ground state) wavefunction, is a function of potentially many variables---one variable for each particle. 
Consequently, the dimensionality of the underlying function space grows exponentially in the number of particles. 
This exponential explosion of cost with respect to the number of particles is called the {\it curse of dimensionality}.

Coupled cluster (CC) theory is a well-established and widely used approach in computational chemistry that is employed to circumvent the curse of dimensionality.
It provides highly accurate approximations to the electronic Schr\"odinger equation. 
Although the CC approach makes molecular systems of tens and some hundreds of correlated electrons numerically tractable, it still faces severe scaling. For example, the CC method truncated at the singles and doubles level, which is the subject of this article, has sixth-order scaling with respect to the system's size.

The working equations called the {\it CC equations}, correspond to a system of polynomial equations. Under certain assumptions, the roots of this polynomial system describe the exact solutions of the considered system.
Therefore, it is important to understand the structural properties of said roots. However, due to the high dimensionality and non-linearity of the CC equations, the algebraic and physical properties of the corresponding solution set are still far from being thoroughly understood.

The CC equations arise from the exponential parametrization of the wavefunction, which was independently derived by Hubbard~\cite{hubbard1957description} and Hugenholtz~\cite{hugenholtz1957perturbation}. 
Pioneering works that applied the exponential ansatz to the fermionic many-body problem include Coester~\cite{coester1958bound}, Coester and Kümmel~\cite{coester1960short}, Čížek~\cite{vcivzek1966correlation,vcivzek1969use}, and Čížek and Paldus~\cite{vcivzek1971correlation}, where the projective nature of the CC equations was proposed and where the polynomial CC equations were derived for the first time.

\subsection{Previous works}

Previous mathematical investigations related to the CC theory were mostly performed within a functional analytic framework, addressing the approximation properties and local convergence results~\cite{schneider2009analysis,rohwedder2013continuous,laestadius2018analysis,laestadius2019coupled,faulstich2019analysis,hassan2022analysis}. 
Within the chemistry community, the root structure of the CC equations has been studied at a fundamental level with the goal of including homotopy methods in the CC methodology.
However, it appears that widespread applications have been hampered by the high dimensionality and non-linearity of the CC equations, as well as the steep scaling of the employed algebraic methods. 
The first study on this topic dates back to 1978 when Živkovič and Monkhorst investigated the singularities and multiple solutions of the single-reference CC equations~\cite{vzivkovic1978analytic}.
This was followed by mathematical and numerical studies of multiple solutions of the single-reference and state-universal multi-reference CC equations and their singularities and analytic properties in the early 1990s by Paldus and coworkers~\cite{piecuch1990coupled,paldus1993application}. 
In 1998, Kowalski and Jankowski revived the homotopy methods in connection with the single-reference CC theory and used them to solve the CC equations with doubles for a minimum-basis-set four-electron problem~\cite{kowalski1998towards}.
This was followed by a fruitful collaboration of Kowalski and Piecuch, who extended the application of the homotopy methods
to the equations defining the CC approaches with singles and doubles (CCSD), singles, doubles,
and triples (CCSDT), and singles, doubles, triples, and quadruples (CCSDTQ)~\cite{piecuch2000search}, again using
a four-electron system described by a minimum basis set as a target. 
They also introduced the formalism of $\beta$-nested equations and proved the {\it Fundamental Theorem of the $\beta$-NE Formalism}, which enabled them to explain the behavior of the curves connecting multiple solutions of the various CC polynomial systems, i.e., from CCSD to CCSDT, CCSDT to CCSDTQ, etc. 
In~\cite{kowalski2000complete2}, Piecuch and Kowalski used homotopy methods to determine all solutions of nonlinear state-universal multireference CCSD equations based on the Jeziorski-Monkhorst ansatz, proving two theorems that provided an explanation for the observed intruder solution problem. 
In a sequel work~\cite{kowalski2000complete}, they used homotopy methods to obtain all solutions of the generalized Bloch equation, which is nonlinear even in a CI parametrization. Most recently, the root structure of the CC polynomial was studied by Csirik and Laestadius using the topological degree theory~\cite{csirik2023coupled}.

\subsection{Contributions}
In this article, we leverage the framework of modern {\it algebro-computational} methods with the goal to bound the number of roots of the CC equations. We complement our results with numerical simulations using state-of-the-art software libraries such as polymake~\cite{polymake_drawing:2010,polymake_XML:ICMS_2016,polymake:2000, polymake:2017, polymake:FPSAC_2009, polymake:ICMS_2006}. We numerically solve tractable examples as well as analyze the CC root structure using mixed volume computations. The latter is used to establish a Bernstein--Khovanskii--Kushnirenko (BKK) inspired bound \cite{Bernshtein1975,Kouchnirenko1976,cox2006using}.
Moreover, we perform calculations based on homotopy continuation methods, using the Bertini~\cite{Bertini} software, to understand the quality of the approximations provided by truncated and untruncated couple cluster parameterizations. 
Other algebraic methods such as univariate rational approximation~\cite{Xiao2021} and multi-resultant techniques~\cite{SYAM2004417} are not considered in this article.
As a first step towards a thorough understanding of this specific family of polynomial systems, we establish a vertex description of surrogate Newton polytopes and prove an alternative description of the polynomials as cubic forms on some affine algebraic variety. 
Due to the exponential scaling of the polynomial system with respect to the number of electrons, already small systems become computationally challenging and might provide an interesting benchmark problem for future algebro-computational techniques.

This article consists of two main parts, the first part,~\cref{sec:ElecSchroed,sec:exponentialparam} describes the underlying structures used in CC theory, i.e., the formulation of the electronic Schr\"odinger equation, excitation operators, and the exponential parametrization, from an algebraic and geometric perspective. 
The second part,~\cref{sec:CCEq,sec:twoInFour,sec:threeInSix}, contains the derivation of the CC equations, and the numerical simulation of two electrons in four spin-orbitals and three electrons in six spin-orbitals. 
Since the proofs of the novel results in this manuscript are rather straightforward but often tedious, we provide them in their entirety in the appendix~\ref{sec:appendix}; for proofs of employed elementary results, we provide the reader with direct references.

\section{Brief Introduction to the Electronic Schr\"odinger Equation}
\label{sec:ElecSchroed}

The electronic Schr\"odinger equation of $N$ electrons in the vicinity of $N_{\rm nuc}$ (clamped) nuclei is an eigenvalue equation
\begin{equation}
H\Psi_* = E_* \Psi_*,
\end{equation}
where the operator $H$, also called the Coulomb Hamiltonian, is given by 
\begin{equation}
\label{eq:ElecSE}
H
=
-\frac{1}{2}\sum_{i=1}^N\Delta_{r_i}
-\sum_{i=1}^N\sum_{j=1}^{N_{\rm nuc}}\frac{Z_j}{\vert r_i-R_j\vert} 
+\sum_{i=1}^N\sum_{j>i}^N\frac{1}{\vert r_i-r_j\vert},
\end{equation}
and $\Psi_*$, the eigenfunction, is an anti-symmetric function in $N$ variables. 
We here denote $\Delta_{r_i} = \sum_{j=1}^3 \partial^2_{r_i^{(j)}}$ with $r_i= (r_i^{(1)},r_i^{(2)},r_i^{(3)})$, $Z_j\in \mathbb{R}$ denotes the j$th$ nuclear charge, $R_j\in\mathbb{R}^3$ denotes the position of the j$th$ nucleus and $r_i\in\mathbb{R}^3$ denotes the position of the i$th$ electron. 
The energetically lowest eigenpair, subsequently denoted $(\Psi_0, E_0)$, is of particular importance to the computational chemistry community~\cite{helgaker2014molecular,yserentant2010regularity}. 
The anti-symmetry is important since we consider electrons, i.e., fermions, which must obey {\it Pauli's exclusion principle}.\footnote{For the physical interpretation of Pauli's exclusion principle, we refer the interested reader to standard quantum mechanics textbooks, e.g., see~\cite{messiah1990albert,helgaker2014molecular}.} 
For the content of this article, it suffices to accept that Pauli's exclusion principle is equivalent to the anti-symmetry property of the eigenfunction and is interpreted as a quantum mechanical axiom~\cite{yserentant2010regularity}.
Since the Hamiltonian in~\cref{eq:ElecSE} is a continuous operator, we require a discretization to be able to perform numerical simulations.
There are in principle two ways of doing this: The first is to use a {\it Galerkin discretization}~\cite{ern2004theory} of the operator in~\cref{eq:ElecSE} which is known as first quantization. The second, and in practice more common approach, is to use the fermionic Fock space~\cite{messiah1990albert,szalay2015tensor}.
This approach is called second quantization and shall be further elaborated on in this article. For a more well-defined introduction to the electronic Schr\"odinger equation, we refer the interested reader to~\cref{app:FuncAna}. 

Either discretization approach results in a linear eigenvalue problem. The challenge of this {\it quantum-many body problem} is the scaling of the discretization matrix with respect to the number of electrons. In the parlance of high-dimensional numerical simulations, this is also referred to as the {\it curse of dimensionality}, in other mathematical fields this might be known as the {\it combinatorial explosion}.  
Over the past century, there have been a great many approaches of different scalings and accuracies to circumvent this curse of dimensionality making the electronic Schr\"odinger equation numerically tractable. 
One of the most widely used high-accuracy approaches is the single reference CC approach, which is the subject of this article. 
Despite the great progress made in the field, we emphasize that the quantum many-body problem still poses one of today's most challenging problems in computational and applied sciences.

\section{Exponential parametrization}
\label{sec:exponentialparam}

The success of CC methods can be traced back to the exponential parametrization of the eigenfunction.
Mathematically, this can be derived in several ways~\cite{schneider2009analysis,rohwedder2013continuous,csirik2023coupled1}; here, we introduce an alternative way following algebraic considerations. 

The above-mentioned first quantization is a useful framework for analyzing the Schr\"odinger equation from a functional analysis perspective~\cite{schneider2009analysis, rohwedder2013continuous,laestadius2018analysis,faulstich2019analysis}. However, for practical implementations as well as for the algebraic and geometric analysis presented in this article, the second quantization is more suitable. 
We will not give a complete characterization of the second quantization, the interested reader is referred to~\cite{messiah1990albert}, instead, we will focus only on the mathematical constructions and objects that are important for the subsequently introduced CC theory. 

\subsection{Creation and annihilation operators}

Let $\mathcal{B}$ be an orthonormal basis of a finite-dimensional Hilbert space $h\subset H^1(X)$, where $X=\mathbb{R}^3 \times \{ \pm 1/2 \}$ and $|\mathcal{B}|=K\gg N$; we call $h$ the single-particle Hilbert space and its basis functions spin-orbitals.
The $M$-particle Hilbert space, denoted $\mathcal{H}^{(M)}$, is then the $M$-th exterior power of $h$ equipped with the natural orthonormal basis, denoted $\mathfrak{B}^{(M)}$, coming from wedge-products of elements of $\mathcal{B}$, which are called Slater determinants. We define the fermionic Fock space as the direct sum of the $M$-particle Hilbert spaces $\mathcal{H}^{(M)}$, i.e., 
\begin{equation}
\label{eq:FockSpace}
\mathcal{F} = \bigoplus_{M=0}^K \mathcal{H}^{(M)},
\end{equation}
which is the Grassmann algebra on $h=\mathcal{H}^{(1)}$. All vector spaces defined above are real vector spaces. Recall that by Pauli's exclusion principle, a spin-orbital can either be occupied or unoccupied, hence, any element in $\mathcal{F}$ can be expressed by means of an occupation vector, i.e., a $K$-tuple of zeros and ones determining the occupation of an individual spin-orbital~\cite{helgaker2014molecular}. Said occupation vectors can be expressed in the compact bra-ket notation, i.e., $\phi = |k_1, k_2, \dots, k_K\rangle$ where $k_i\in\{0,1\}$ denotes the occupation number of the $i$-th spin-orbital. 
We can now define the creation operator $a_i^\dagger:\mathcal{F} \to \mathcal{F}$ as taking the wedge product with the $i$-th basis vector of $h$ in $\mathcal{B}$, and the annihilation operator $a_i:\mathcal{F} \to \mathcal{F}$ as contracting with $i$-th basis element.

Given the fermionic creation and annihilation operators, the Coulomb Hamiltonian can be reformulated in second quantized form~\cite{messiah1990albert}; it reads
\begin{equation}
H = \sum_{p,q = 1}^K h_{p,q}a_p^\dagger a_q + \frac{1}{2}\sum_{p,q,r,s=1}^K v_{p,q,r,s} a_p^\dagger a_r^\dagger a_s a_q,
\end{equation}
where 
\begin{equation}
h_{p,q} = \int_X \chi_p^*(x_1) \left(
-\frac{\Delta}{2} - \sum_{j}\frac{Z_j}{|r_1 - R_j|}
\right)\chi_q(x_1)dx_1
\end{equation}
and 
\begin{equation}
v_{p,q,r,s} 
=
\int_{X \times X}
\frac{
\chi_p^*(x_1) \chi_q(x_1) \chi_r^*(x_2) \chi_s(x_2)
}{|r_1 - r_2|} dx_1 dx_2,
\end{equation}
with $x_i = (r_i,s_i)$, i.e. tuples of space and spin coordinates.
Arguably, the most important property of the creation and annihilation operators is the canonical anti-commutation relation (CAR).

\begin{lemma}[CAR]
\label{Thm:CAR}
The fermionic creation and annihilation operators fulfill the canonical anti-commutation relations, i.e., 
\begin{equation}
[a_i^\dagger, a_j]_+ 
=
\delta_{i,j}
\quad 
\text{and}
\quad 
[a_i^\dagger, a_j^\dagger]_+
=
[a_i, a_j]_+ 
=0.
\end{equation}
\end{lemma}
The proof can be found in e.g.~Ref~\cite{helgaker2014molecular}. We can now define the vector space $\mathfrak{a} = {\rm Span}\{a_i\}_i\subseteq {\rm End}(\mathcal{F})$ and the corresponding coordinate map
\begin{equation}
a:\mathbb{R}^K \to \mathfrak{a}~;~ v = (v_1,...,v_K) \mapsto v_1a_1+...+v_Ka_K.
\end{equation}
The canonical anti-commutation relations from Theorem~\ref{Thm:CAR} then generalize on the vector space $\mathfrak{a}$, which is commonly referred to as the representation of the CARs of $K$ degrees of freedom.  

\begin{remark}
The above construction is the starting point for a description of fermionic ladder operators in the more general framework of $C^*$-algebras~\cite{arveson2012invitation}.
\end{remark}

\subsection{Excitation operators}

Having the basic objects at hand we can now define the excitation operators which are essential in CC theory. 
The $i$-$a$ (single) excitation operator is then given by 
\begin{equation}
X_i^a = a^\dagger_a a_i.
\end{equation}
Since the excitation operators are defined by means of the creation and annihilation operators, their action on Slater determinants is simply given by contracting with $i$-th basis
element in $\mathcal{B}$ followed by taking the wedge product with the $a$-th basis element in $\mathcal{B}$. 
Higher-order excitation operators correspond to the normal-ordered distinct products of such single excitation operators, e.g.,
\begin{equation}
X_{ij}^{ab}= \{X_i^aX_j^b\} 
= a^\dagger_aa^\dagger_b a_i a_j,
\end{equation}
where $\{\cdot\}$ highlights the normal ordering (see e.g.~\cite{messiah1990albert,helgaker2014molecular} for the definition of normal ordered operators). We furthermore introduce a canonical order on the indices, i.e., $i<j$ and $a<b$, introducing an additional phase, i.e., a sign.  For a general excitation operator $X^{a_1,...,a_k}_{i_1,...,i_k}$ we define the excitation multiindex
\begin{equation}
\mu=
\begin{pmatrix}
a_1,...,a_k\\ i_1,...,i_k
\end{pmatrix},
\end{equation}
and denote $X_\mu = X^{a_1,...,a_k}_{i_1,...,i_k}$---we call $|\mu| = k$ the excitation rank.
We subsequently restrict the considerations to excitations with $i_p\in\{1,...,N\}$ and $a_p \in \{ N+1,...K\}$, and refer to $\Phi_0=\Phi[1,...,N]$ as the reference Slater determinant. 
Moreover, the orbitals $\{1,...,N\}$ are called {\it occupied orbitals}, $\{ N+1,...K\}$ and called {\it virtual orbitals}. 
The excitation indices $\mu$ that excite from the occupied into the virtual orbitals define the multi-index set $\mathcal{I}$.  
Since this set of excitations corresponds to simply replacing indices in the string $[1,..., N]$ with indices in the string $[N+1,..., K]$ (plus some additional permutation), we may deduce that there is a one-to-one relation between excitation operators and Slater determinants except for the reference Slater determinant $\Phi_0$. 
In other words, the excitation operators map the reference Slater determinant $\Phi_0$ to all other Slater determinants. Hence, the basis can be recovered as
\begin{equation}
\mathfrak{B} = \mathfrak{B}^{(N)} = \{\Phi_0\}\cup \{\Phi_\mu=X_\mu\Phi_0 ~|~ \mu \in \mathcal{I} \}.
\end{equation}
This allows us to introduce an alternative representation of Slater determinants which is very useful for practical considerations. Let $\mu$ be a multi-index of rank $k$, then we may write the Slater determinant $\Phi_\mu$ as 
\begin{equation}
\Phi_\mu
=
X_\mu \Phi_0
=
X_{i_1,...,i_k}^{a_1,...,a_k}\Phi_0
=
\Phi_{i_1,...,i_k}^{a_1,...,a_k}.
\end{equation}
Moreover, assuming that $\mathfrak{B}$ is ordered by the Slater determinant's excitation rank, the excitation operators correspond to strictly lower-triangular block matrices and are consequently nilpotent. 
\begin{lemma}
\label{thm:ExcitationOp}
The excitation operators fulfill the commutation relations 
\begin{equation}
    \begin{aligned}
    \,[X_\mu, X_\nu]  = X_\mu X_\nu -  X_\nu X_\mu &= 0 \\
    [X_\mu^\dagger, X_\nu^\dagger] = X_\mu^\dagger X_\nu^\dagger -  X_\nu^\dagger X_\mu^\dagger &= 0
    \end{aligned}
    \end{equation}
\end{lemma}
The proof can be found in Ref.~\cite{helgaker2014molecular,csirik2023coupled1}. We emphasize that $[X_\mu, X^\dagger_\nu]$ is not necessarily zero. However, the following statement can be made.
\begin{lemma}
\label{lem:Commutation}
Let $X_\mu, X_\nu$ be excitation operators with
\begin{equation}
\mu=
\begin{pmatrix}
a^{(\mu)}_1,...,a^{(\mu)}_{k^{(\mu)}}\\ i^{(\mu)}_1,...,i^{(\mu)}_{k^{(\mu)}}
\end{pmatrix},
\qquad
\nu=
\begin{pmatrix}
a^{(\nu)}_1,...,a^{(\nu)}_{k^{(\nu)}}\\ i^{(\nu)}_1,...,i^{(\nu)}_{k^{(\nu)}}
\end{pmatrix},
\end{equation}
and of rank $k^{(\mu)}, k^{(\nu)}$, respectively.
Moreover, let $a^{(\mu)}_j \neq a^{(\nu)}_\ell$, $i^{(\mu)}_j \neq i^{(\nu)}_\ell$ for all pairs $(j,\ell)$, then
\begin{equation}
[X_\mu, X^\dagger_\nu]=0.
\end{equation}
\end{lemma}

For the $N$-particle Hilbert space $\mathcal{H}$, we can define the set of all creation operators defined on $\mathcal{H}$ denoted by $\mathfrak{E}(\mathcal{H})$. 

\begin{theorem}
\label{th:Basis<->operator}
There exists a one-to-one relation between the $N$-particle basis functions $\mathfrak{B}^{(N)}$ and $\mathfrak{E}(\mathcal{H})\cup\{I\}$.
\end{theorem}
Then, we define the vector space
\begin{equation}
\mathfrak{b} = {\rm Span}\big( \mathfrak{E}(\mathcal{H})\big),
\end{equation} 
which is known in the quantum-chemistry community as the space of cluster operators, i.e., 
\begin{equation}
\mathfrak{b} = \left\lbrace T = \sum_{k=1}^NT^{(k)}~\bigg|~T^{(k)} = \sum_{|\mu| = k} t_\mu X_\mu \right\rbrace.
\end{equation}
This vector space takes an essential role in CC theory~\cite{helgaker2014molecular,faulstich2020mathematical}. 

\subsection{Coupled cluster parametrization}

From a mathematical perspective, it is natural to approach the characterization of the fermionic many-body problem from a Ritz--Galerkin perspective~\cite{ern2004theory}, i.e., 
\begin{equation}
\label{eq:RitzGalerkin}
E_0 = \min_{\Psi \in \mathcal{H}}\frac{\langle \Psi , H \Psi \rangle}{\langle \Psi , \Psi \rangle }.
\end{equation}
Note that the above expression is invariant w.r.t~a normalization factor of $\Psi$.  
We denote $\mathcal{H}_{\rm int}$ the affine subspace of $\mathcal{H}$ consisting of all functions that are intermediately normalized w.r.t.~the reference state $\Phi_0$, i.e., $\langle \Psi, \Phi_0 \rangle = 1$.
The construction of excitation operators allows us to transfer the degrees of freedom from the function space to the {\it wave operator}. 
We define the wave-operator map
\begin{equation}
\Omega : \mathfrak{b} \to \mathcal{G}~;~ C\mapsto Id + C,
\end{equation}
where 
\begin{equation}
\mathcal{G} = \{ Id + C~|~C\in\mathfrak{b} \}.
\end{equation}
We show the one-to-one correspondence between intermediately normalized functions $\Psi\in \mathcal{H}_{\rm int} \subset \mathcal{H}$, and cluster operators $c\in\mathfrak{b}$. 
To that end, we require the following two lemmata:

\begin{lemma}
\label{lemma:1}
Let $\Psi \in\mathcal{H}_{\rm int}$. There exists a unique element $(Id +C)\in\mathcal{G}$, s.t., 
\begin{equation}
\Psi = (Id +C)\Phi_0.
\end{equation}
\end{lemma}

\begin{lemma}
\label{lemma:2}
The wave-operator map $\Omega$ is bijective. 
\end{lemma}
Combining these two lemmata, we can show the desired one-to-one correspondence. 
\begin{theorem}
\label{th:fci}
Let $\Psi \in\mathcal{H}_{\rm int}$. There exists a unique element $C\in\mathfrak{b}$, s.t., 
\begin{equation}
\Psi = \Omega(C)\Phi_0.
\end{equation}
\end{theorem}
Although we have restricted the above theorem to intermediately normalized functions (the reason will become apparent shortly),~\cref{th:fci} is in fact the core of the (full) configuration interaction expansion~\cite{helgaker2014molecular}. 
Note moreover, since $\mathfrak{b}$ is nilpotent the exponential series $\exp(T)$ for any element $T\in\mathfrak{b}$ terminates after at most $N$ terms. 
We can therefore define the set
\begin{equation}
\tilde{\mathcal{G}}
=
\left\lbrace
\exp(T)~|~T \in \mathfrak{b}
\right\rbrace,
\end{equation}
where $\exp(T) = Id + P(T)$ and $P$ is a polynomial of at most degree $N$.
Since $\mathcal{G}$ is a unipotent algebraic group, the exponential map from $\mathfrak{b}$ to $\mathcal{G}$ is a bijection.
Hence, any function that is intermediately normalized can be uniquely expressed by means of an element in $\mathcal{G}$ through the exponential of a cluster operator $T\in\mathfrak{b}$. 
This aligns with the known functional analytic results~\cite{schneider2009analysis,rohwedder2013continuous,laestadius2018analysis,faulstich2019analysis}, and is known in the quantum-chemistry community as the equivalence of Full Interaction Configuration (FCI) and Full CC (FCC). 

\begin{remark}
In order to make practical computations numerically tractable, the cluster operator is truncated in some way. 
This could be either by the level of excitation, i.e., CCSD, etc.~\cite{helgaker2014molecular}, by means of additional system-specific information~\cite{gururangan2021high, chakraborty2022benchmarking}, or active space considerations~\cite{magoulas2022addressing,Piecuch2010,veis2016coupled,faulstich2019numerical,evangelista2011alternative,hino2006tailored}. 
\end{remark}
\section{Coupled Cluster Equations}\label{sec:CCEq}
In order to approximate solutions to the Schr\"odinger equation using CC theory, we need to derive working equations that can be solved numerically. 
There is a variety of ways the exponential ansatz can be used to derive numerical schemes that approximately solve the Schr\"odinger equation. 
We shall focus here on the most common application, i.e., the projected single-reference CC approach. 
In this approach, we define the CC energy as 
\begin{equation}
\label{eq:cc-energy}
\mathcal{E}_{\rm CC} = \langle\Phi_0, e^{-T}He^T\Phi_0\rangle,
\end{equation}
where $T\in \mathfrak{b}$ fulfills the CC equations
\begin{equation}
\label{cc:amplitude}
\langle\Phi_\mu, e^{-T}He^T\Phi_0\rangle = 0, \quad \forall \Phi_\mu \in \mathfrak{B} \setminus \{ \Phi_0\}.
\end{equation}

Note that $T$ fulfilling Eq.~\eqref{cc:amplitude} is not unique. 
In fact, different cluster matrices describe different eigenfunctions and by taking the energetically lowest solution, this ansatz yields the FCI ground state energy (i.e. the Galerkin energy~\cite{laestadius2019coupled,csirik2023coupled1}).
Put differently, if all eigenfunctions of $H$ have a non-zero projection onto $\Phi_0$, then, the untruncated CC ansatz can describe the spectrum of $H$.

If the CC ansatz is truncated, i.e., $T\in\bar{\mathfrak{b}} \subseteq \mathfrak{b}$, the above is not necessarily true (see e.g.~\cite{helgaker2014molecular}).
Moreover, it is worth noting that the CC equations (Eq.~\eqref{cc:amplitude}) do not arise from the Rayleigh quotient, hence, $\mathcal{E}_{\rm CC}$ is not necessarily bounded from below by the FCI ground state energy, i.e, {\it the CC method is not variational}.

That being said, the exponential parameterization has a major advantage over the linear parameterization, namely, for non-interacting systems it factorizes.  
This is a highly desired property since it yields physical quantities (like the system's energy) to be correctly described in certain limits; in quantum chemistry parlance, this is called size consistency~\cite{helgaker2014molecular}. 
Hence, the CC ansatz yields a physically accurate description of the system which makes truncated CC methods more favorable than a truncated linear approximation.

Since each operator $T\in\mathfrak{b}$ can be uniquely described by its expansion coefficients $t = (t_\mu)_{\mu}$, the CC equations can be understood as amplitude equations, i.e., equations that characterize the amplitudes $t$ instead of the operator $T$. 
From this perspective, the CC amplitudes correspond to the roots of the CC function
\begin{equation}
f_\mu(t) = \langle\Phi_\mu, e^{-T(t)}He^{T(t)}\Phi_0\rangle,\quad \forall\mu\in\mathcal{I} ,
\end{equation}
which can be either real or complex. In this article, we consider the complex-valued amplitudes.
A characteristic of the projected CC approach is that the expansion of the similarity-transformed Hamiltonian truncates after at most four terms~\cite{helgaker2014molecular}, i.e., 
\begin{equation}
e^{-T}H e^T 
=H + [H,T] + \frac{1}{2}[[H,T],T] + \frac{1}{6}[[[H,T],T],T] + \frac{1}{24}[[[H,T],T],T],T] 
\end{equation}
and, 
hence, every polynomial $f_\mu(t)$ is at most of degree four. 
Since we have as many polynomial equations as we have excited Slater determinants, i.e., we exclude the reference state, B\'ezout's theorem yields the rough upper bound of
\begin{equation}
\mathcal{N} \leq 4^{\mathcal{K}-1},
\end{equation}
where $\mathcal{N}$ is the number of solutions to the set of CC equations and $\mathcal{K}-1$ is the number of excited Slater determinants. 
Note, that already in the case of four electrons in eight spin-orbitals, for example for LiH in the minimal basis description, this na\"ive B\'ezout bound for CC truncated at the single and double level of excitation (CCSD) is $2^{104}\approx 2\cdot 10^{31}$. Any homotopy algorithm based on this B\'ezout bound is therefore infeasible. We emphasize that this bound grossly overestimates the number of roots since no further structure of the CC equations is taken into account. Clearly, more information on the number of roots of the projected CC equations would be beneficial for various numerical approaches that characterize the roots of the CC function.

In order to derive more accurate bounds to the number of roots, we require a more explicit description of the working equations. Since a general description of these equations for arbitrary truncations is not feasible, we subsequently focus on the most popular truncation scheme, i.e., the CCSD method.
This results in a restriction of the projections onto merely singly and doubly excited Slater determinants. 
Although this seems like an aggressive truncation, the expressions of the working equations are still highly complicated. 
Therefore, we split the description of $f_\mu$ into the case when $\mu$ describes singly excited Slater determinants (single-excitation CCSD equations), and the case when $\mu$ describes doubly excited Slater determinants (double-excitation CCSD equations).
We note that the subsequent derivation of the CCSD amplitude equations can be made more compact when advanced diagrammatic techniques based on Feynman-style (e.g., Hugenholtz and Brandow) diagrams are used~\cite{shavitt2009many}.

We investigate the different terms of the Baker--Campbell--Hausdorff expansion individually. 
It holds
\begin{align*}
    [H,T] &= HT-TH,\\
    [[H,T],T] 
    &= HT^2-2THT+T^2H,\\
    [[[H,T],T],T]
    &= HT^3-3THT^2+3T^2HT-T^3H,\\
    [[[[H,T],T],T],T]
    &= HT^4-4THT^3+6T^2HT^2-4T^3HT+T^4H.
\end{align*}
Some very important rules of quantum chemistry simplify the higher order terms~\cite{helgaker2014molecular}: Since we are only interested in the CCSD case all terms of the form $T^nHT^m$ with $m$ arbitrary and $n>2$ vanish since we are projecting onto at most doubly excited Slater determinants. Similarly, we find that $$\langle\Phi_i^a,T^2HT^m\rangle = 0.$$
Furthermore, the Slater--Condon rules~\cite{messiah1990albert} imply that 
\begin{align*}
    \langle \Phi_0, H T^n\Phi_0\rangle &= 0, \quad \text{for }n>2\\
    \langle \Phi_i^a, H T^n\Phi_0\rangle &= 0, \quad \text{for }n>3\\
    \langle \Phi_{ij}^{ab}, H T^n\Phi_0\rangle &= 0, \quad \text{for }n>4
\end{align*}
and by de-excitation rules
\begin{align*}
    \langle \Phi_i^a, THT^m\Phi_0\rangle &= t_i^a\langle\Phi_0, HT^m\Phi_0\rangle.
\end{align*}

These rules may now be applied to expand the CCSD equations. To that end, we investigate the individual projections of the respective commutator terms. 
We exemplify this procedure with the standard commutator expression projected onto a singly excited Slater determinant. This yields    
\begin{align*}
\langle \Phi_i^a, [H,T]\Phi_0\rangle  
&=  
\langle \Phi_i^a, HT\Phi_0\rangle
- 
\langle \Phi_i^a, TH\Phi_0\rangle\\
&=
\sum_{j,b}t_j^b\langle\Phi_i^a,H\Phi_j^b\rangle
+
\sum_{\substack{j<k\\b<c}}
t_{jk}^{bc}\langle\Phi_i^a,H\Phi_{jk}^{bc}\rangle
-
t_i^a\langle\Phi_0, H\Phi_0\rangle,
\end{align*}
i.e., a linear polynomial in the cluster amplitudes $t_j^b, t_{ij}^{ab}$.
We emphasize that this procedure rapidly gets more involved, and refer the interested reader for a complete derivation of the CC equations to~\cref{sec:ccWorkingEQ}.  
Alternatively, these (nested) commutators of the Hamiltonian with the elementary excitation operators, which generate excited determinants, can
be directly calculated using more advanced techniques like Wick’s theorem or other contraction theorems~\cite{helgaker2014molecular}.
Note that due to the Slater--Condon rules, there exists no quadruple commutator contribution in the single-projection equations.
Hence, the resulting polynomials that arise from projecting onto singly excited Slater determinants are of the form
\begin{equation}
\label{eq:singleprojections}
\begin{aligned}
f_i^a(t)
&=
\langle
\Phi_i^a, H\Phi_0 
\rangle
+
\langle
\Phi_i^a, [H,T]\Phi_0 
\rangle
+\frac{1}{2}
\langle
\Phi_i^a, [[H,T],T]\Phi_0 
\rangle\\
&\quad+\frac{1}{6}
\langle
\Phi_i^a, [[[H,T],T],T]\Phi_0 
\rangle\\
&=\langle
\Phi_i^a, H\Phi_0 
\rangle
+
\langle
\Phi_i^a, [H,T]\Phi_0 
\rangle
+\frac{1}{2}
\langle
\Phi_i^a, HT^2\Phi_0 
\rangle-\langle
\Phi_i^a, THT\Phi_0 
\rangle\\
&\quad+\frac{1}{6}
\langle
\Phi_i^a, HT^3\Phi_0 
\rangle-\frac 1 2 \langle
\Phi_i^a, THT^2\Phi_0 
\rangle.
\end{aligned}
\end{equation}
The expansion of the individual terms yields a polynomial of the form 
\begin{equation}
\label{eq:singles-ccsdEQ}
\begin{aligned}
f_i^a(t)
&=
\mathcal{C}(i,a)
+\sum_{k,c} \mathcal{C}(k,c)\,t_{ik}^{ac}
+\frac{1}{2}\sum_{k,c,d} \mathcal{C}(a,k,c,d)\, t_{ik}^{cd}
-\frac{1}{2}\sum_{k,l,c} \mathcal{C}(k,l,i,d)\, t_{kl}^{ac}
\\&\quad 
+\sum_c \mathcal{C}(a,c)\, t_i^c
-\sum_k  \mathcal{C}(k,i)\,  t_k^a
+\sum_{k,c}  \mathcal{C}(a,k,i,c)\, t_{k}^{c}
\\&\quad 
-\frac{1}{2}\sum_{k,l,c,d}  \mathcal{C}(k,l,c,d)\,\rangle (t_i^c  t_{kl}^{ad}+ t_k^a  t_{il}^{cd}- t_k^c  t_{li}^{da})
-\sum_{k,c}\mathcal{C}(k,c) t_i^c  t_k^a
\\&\quad 
+\sum_{k,c,d} \mathcal{C}(a,k,c,d) \, t_i^c  t_k^d
- \sum_{k,l,c}\mathcal{C}(k,l,i,c) \, t_k^a  t_l^c
- \sum_{k,l,c,d} \mathcal{C}(k,l,c,d) \, \rangle t_i^c  t_k^a  t_l^d ,
\end{aligned}
\end{equation}
where $\mathcal{C}$ denotes a coefficient map that can be assembled by the coefficients that are derived in~\cref{sec:ccWorkingEQ}. We note that $\mathcal{C}$ obeys some structural rules which are partially exploited in the subsequent paragraphs, see~\cref{sec:ccWorkingEQ} for more details.
We see from~\cref{eq:singles-ccsdEQ} that the single-excitation CCSD equations consists of $n_{\rm s}$ equations of at most degree three, where $
n_{\rm s}
=
N(K-N)$.

In a similar fashion, we can derive the polynomials that arise from projecting onto doubly excited Slater determinants. 
Here we find
\begin{equation}
\label{eq:doubleprojections}
\begin{aligned}
f_{ij}^{ab}(t)
&=
\langle
\Phi_{ij}^{ab}, H\Phi_0 
\rangle
+
\langle
\Phi_{ij}^{ab}, [H,T]\Phi_0 
\rangle
+\frac{1}{2}
\langle
\Phi_{ij}^{ab}, [[H,T],T]\Phi_0 
\rangle\\
&\quad+\frac{1}{6}
\langle
\Phi_{ij}^{ab}, [[[H,T],T],T]\Phi_0 
\rangle 
+\frac{1}{24}
\langle
\Phi_{ij}^{ab}, [[[[H,T],T],T],T]\Phi_0 
\rangle 
.
\end{aligned}
\end{equation}
The full derivation of the terms is again performed in~\cref{sec:ccWorkingEQ}. Without going into full detail of the coefficient distribution, we present the general form of the polynomial that arises from the projections onto the doubly excited Slater determinants
\begin{equation}
\label{eq:doubles-ccsdEQ}
\begin{aligned}
f_{i,j}^{a,b}(t)
&=
\mathcal{C}(a,b,i,j)
+\sum_k \mathcal{C}(k,j)\, t_{ik}^{ab}
+\frac{1}{2}\sum_{cd} \mathcal{C}(a,b,c,d)\, t_{ij}^{cd}
\\&\quad
+\frac{1}{2}\sum_{kl} \mathcal{C}(k,l,i,j)\, t_{kl}^{ab}
+ \sum_{kc} \mathcal{C}(k,b,c,j)\, t_{ik}^{ac}
\\&\quad
+\frac{1}{4} \sum_{klcd} \mathcal{C}(k,l,c,d)\, \rangle t_{ij}^{cd} t_{kl}^{ab}
+\sum_{klcd} \mathcal{C}(k,l,c,d)\, t_{ij}^{cd} t_{ik}^{ac} t_{jl}^{bd}
\\&\quad
-\frac{1}{2}  \sum_{k l c d}\mathcal{C}(k,l,c,d)\, t_{i k}^{d c} t_{l j}^{a b}
-\frac{1}{2} \sum_{k l c d}\mathcal{C}(k,l,c,d)\, t_{l k}^{a c} t_{i j}^{d b}
\\&\quad
+ \sum_{c} \mathcal{C}(a,b,c,j)\,  t_{i}^{c}
-\sum_{k} \mathcal{C}(k,b,i,j)\, t_{k}^{a}
-\sum_{k c} \mathcal{C}(k,c)\, t_{i}^{c} t_{k j}^{a b} 
\\&\quad
-\sum_{k c} \mathcal{C}(k,c) t_{k}^{a} t_{i j}^{c b}
+ \sum_{k c d}\mathcal{C}(a,k,c,d) t_{i}^{c} t_{k j}^{d b} 
-\sum_{k l c}\mathcal{C}(k,l,i,c) t_{k}^{a} t_{l j}^{c b}
\\&\quad
-\frac{1}{2} \sum_{k c d}\mathcal{C}(k,b,c,d) t_{k}^{a} t_{i j}^{c d}
+\frac{1}{2} \sum_{k l c}\mathcal{C}(k,l,c,j) t_{i}^{c} t_{k l}^{a b}
\\&\quad
+\sum_{k c d}\mathcal{C}(k,a,c,d) t_{k}^{c} t_{i j}^{d b}
- \sum_{k l c}\mathcal{C}(k,l,c,i) t_{k}^{c} t_{l j}^{a b}
\\&\quad
+\sum_{c d}\mathcal{C}(a,b,c,d)  t_{i}^{c} t_{j}^{d}
+\sum_{k l}\mathcal{C}(k,l,i,j)  t_{k}^{a} t_{l}^{b} 
\\&\quad
-\sum_{k c}\mathcal{C}(k,b,c,j)  t_{i}^{c} t_{k}^{a}
+\frac{1}{2} \sum_{k l c d}\mathcal{C}(k,l,c,d)  t_{i}^{c} t_{j}^{d} t_{k l}^{a b}
\\&\quad
+\frac{1}{2} \sum_{k l c d}\mathcal{C}(k,l,c,d)  t_{k}^{a} t_{l}^{b} t_{i j}^{c d} 
- \sum_{k l c d}\mathcal{C}(k,l,c,d) t_{i}^{c} t_{k}^{a} t_{l j}^{d b}
\\&\quad
- \sum_{k l c d}\mathcal{C}(k,l,c,d) t_{k}^{c} t_{i}^{d} t_{l j}^{a b} 
- \sum_{k l c d}\mathcal{C}(k,l,c,d) t_{k}^{c} t_{l}^{a} t_{i j}^{d b}
\\&\quad
+ \sum_{k c d}\mathcal{C}(k,d,c,d) t_{i}^{c} t_{k}^{a} t_{j}^{d} 
+ \sum_{k l c}\mathcal{C}(k,l,c,j) t_{i}^{c} t_{k}^{a} t_{l}^{b}
\\&\quad
+\sum_{k l c d}\mathcal{C}(k,l,c,d) t_{i}^{c} t_{j}^{d} t_{k}^{a} t_{l}^{b}.
\end{aligned}
\end{equation}

From Eq.~\eqref{eq:doubles-ccsdEQ}, we see that we have $n_{\rm d}$ equations of at most degree four, where 
\begin{equation}
n_{\rm d} 
=
\frac{N(N-1)}{2}\frac{(K-N)(K-N-1)}{2}.
\end{equation}
This yields a refined B\'ezout bound of $$\mathcal{N}\leq 3^{n_{\rm s}}4^{n_{\rm d}},$$
which has been known and reported in the quantum chemistry community~\cite{piecuch2000search}.

\begin{remark}
The CC equations yield a significant computational advantage: 
Instead of tracking the actual matrix $T\in\mathbb{C}^{\mathcal{K}}$, where $\mathcal{K}= {{K}\choose{N}}$, and then exponentiating this huge matrix to obtain $\exp(T)$, using the CC equations allows to merely track the CC amplitudes where the singles amplitudes are a matrix of size $N \times (K-N)$ and the double amplitudes are a fourth-order tensor of size  $N \times N \times (K-N) \times (K-N)$. Moreover, introducing beneficial factorizations in the CC equations~\cite{helgaker2014molecular} yields a complexity of $\mathcal{O}((K-N)^4N^2)$ to perform all contractions. This is numerically much more feasible than computing $\exp(-T)H\exp(T)$ at the matrix level, which results in $2\cdot N+3$ matrix products that (na\"ivly) scale as $\mathcal{O}(\mathcal{K}^3)$.
\end{remark}

\subsection{Newton Polytopes} 
In order to advance towards a BKK-type bound, we provide the vertex description of two polytopes $N_S$ and $N_D$. The polytope $N_S$ is a superset for the Newton polytopes corresponding to the equations coming from projections to singly-excited Slater determinants. Similarly, $N_D$ is a superset for the Newton polytopes corresponding to the equations coming from projections to doubly-excited Slater determinants.
In principle, these supersets can be used to bound the mixed volume of the actual Newton polytopes. Yet, for this to be of practical relevance, we require a numerically tractable way to compute the mixed volume of these surrogate polytopes, which is left for future investigations. 

We start by characterizing the sets of all possible vertices corresponding to monomials in the singly-excited projections. We define the sets $N_{\rm s} = \{1,\dots,n_{\rm s}\}$ and $N_{\rm d} = \{1+n_{\rm s},\dots,n_{\rm d}+n_{\rm s}\}$.
All vertices describing the single-excited CCSD equations are then described by the following sets $\mathcal{S}_i$:
\begin{equation*}
\begin{aligned}
\mathcal{S}_1
&=\{v~|~ \exists i\in N_{\rm s}  \;s.t.\; v_i=1\;and\; v_j=0\ \forall j\in N_{\rm s}\setminus \{i\} \text{ or } j\in N_{\rm d} \}\\
\mathcal{S}_2
&=\{v~|~ \exists i\in N_{\rm d}  \;s.t.\; v_i=1\;and\; v_j=0\ \forall j\in N_{\rm d}\setminus \{i\} \text{ or } j\in N_{\rm s} \}\\
\mathcal{S}_3
&=\{v~|~ \exists i\neq j\in N_{\rm s}  \;s.t.\; v_i =v_j=1\;and\; v_k=0\ \forall k\in N_{\rm s}\setminus \{i,j\} \text{ or } k\in N_{\rm d} \}\\
\mathcal{S}_4
&=\{v~|~ \exists i\in N_{\rm s},\exists j\in N_{\rm d}  \;s.t.\; v_i =v_j=1\;and\; v_k=0\ \forall k\in N_{\rm s}\setminus \{i\} \text{ or } k\in N_{\rm d}\setminus \{j\} \}\\
\mathcal{S}_5
&=\{v~|~ \exists i\neq j\neq k\in N_{\rm s}  \;s.t.\; v_i =v_j= v_k =1\;and\; v_\ell=0\ \forall \ell\in N_{\rm s}\setminus \{i,j,k\} \text{ or } \ell\in N_{\rm d} \}\\
\mathcal{S}_6
&=\{v~|~ \exists i\in N_{\rm s}  \;s.t.\; v_i=2\;and\; v_j=0\ \forall j\in N_{\rm s}\setminus \{i\} \text{ or } j\in N_{\rm d} \}\\
\mathcal{S}_7
&=\{v~|~ \exists i\neq j\in N_{\rm s} \;s.t.\; v_i=2 \;and\; v_j=1 \;and \; v_k=0\ \forall k\in N_{\rm s}\setminus \{i,j\} \text{ or } k\in N_{\rm d}  \}
\end{aligned}
\end{equation*}
\begin{remark}
A careful analysis of the coefficient map $\mathcal C$ reveals that terms of the form $\mathcal{C}(k,l,c,d) \,  t_i^c  t_k^a  t_l^d$ can only give rise to monomials of the form $xyz$ or $x^2y$, where $x,y,z$ are some distinct single excitation amplitudes. Thus, $\mathcal S_7$ and $\mathcal S_5$ are the only contributions of cubic terms.
\end{remark}
The Newton polytopes corresponding to the monomials in the singly-excited projections lie in the convex hull of the derived vertices $\mathcal{S}_i$. 
In order to find the extremal points of said convex hull we establish the following lemmata of relations between the individual sets, the proofs are again provided in the appendix.

\begin{lemma}\label{lem:set_inclus_sing_Newton_Pol}
The following set of inclusions hold:
\begin{align}
\mathcal{S}_3 &\subset {\rm conv}(\mathcal{S}_6)&
\mathcal{S}_1 &\subset {\rm conv}(\mathcal{S}_6\cup \{0\}).
\end{align}
\end{lemma}
\begin{lemma}\label{lem:S_5_not_extr}
The vertices described by $\mathcal{S}_5$ are not extremal points of the convex hull of $\mathcal S_1,\dots,\mathcal S_7$.
\end{lemma}

\begin{lemma}\label{lem:extrm_vert_sing}
The vertices described by $\{0\},\mathcal{S}_2, \mathcal{S}_4, \mathcal{S}_6, \mathcal{S}_7$
are extremal vertices of the convex hull of $\mathcal S_1,\dots,\mathcal S_7$.
\end{lemma}
We will then define the surrogate polytope $N_S$ as follows,$$N_S\colon= {\rm conv}(\{0\}\cup\mathcal{S}_2\cup\mathcal{S}_4\cup\mathcal{S}_6\cup\mathcal{S}_7).$$
Similarly to the projections onto singly excited Slater determinants, we characterize the extremal vertices of the monomials arising in the projection to doubly-excited Slater determinants. All vertices describing double-excited CCSD equations are described by the following sets $\mathcal{D}_i$:

\begin{equation*}
\begin{aligned}
\mathcal{D}_1
&=
\{v~|~ \exists i\in N_{\rm d}  \;s.t.\; v_i=1\;and\; v_j=0\ \forall j\in N_{\rm d}\setminus \{i\} \text{ or } j\in N_{\rm s} \}\\
\mathcal{D}_2
&=\{v~|~ \exists i\neq j\in N_{\rm d}  \;s.t.\; v_i =v_j=1\;and\; v_k=0\ \forall k\in N_{\rm d}\setminus \{i,j\} \text{ or } k\in N_{\rm s} \}\\
\mathcal{D}_3
&=
\{v~|~ \exists i\in N_{\rm d}  \;s.t.\; v_i=2\;and\; v_j=0\ \forall j\in N_{\rm d}\setminus \{i\} \text{ or } j\in N_{\rm s} \}\\
\mathcal{D}_4
&=
\{v~|~ \exists i\in N_{\rm s}  \;s.t.\; v_i=1\;and\; v_j=0\ \forall j\in N_{\rm s}\setminus \{i\} \text{ or } j\in N_{\rm d} \}\\
\mathcal{D}_5
&= \{v~|~ \exists i\neq j\in N_{\rm s}  \;s.t.\; v_i =v_j=1\;and\; v_k=0\ \forall k\in N_{\rm s}\setminus \{i,j\} \text{ or } k\in N_{\rm d} \}\\
\mathcal{D}_6
&=\{v~|~ \exists i\in N_{\rm s},\exists j\in N_{\rm d}  \;s.t.\; v_i =v_j=1\;and\; v_k=0\ \forall k\in N_{\rm s}\setminus \{i\} \text{ or } k\in N_{\rm d}\setminus \{j\} \}\\
\mathcal{D}_7
&=\{v~|~ \exists i\neq j\in N_{\rm s},\exists k\in N_{\rm d}  \;s.t.\; v_i =v_j=v_k=1\;and\; v_\ell=0\ \forall \ell\in N_{\rm s}\setminus \{i,j\} \text{ or } \ell\in N_{\rm d}\setminus \{k\} \}\\
\mathcal{D}_8
&=
\{v~|~ \exists i\neq j\neq k\in N_{\rm s}  \;s.t.\; v_i =v_j= v_k =1\;and\; v_\ell=0\ \forall \ell\in N_{\rm s}\setminus \{i,j,k\} \text{ or } \ell\in N_{\rm d} \}\\
\mathcal{D}_9
&=\{v~|~ \exists i\neq j\neq k\neq \ell\in N_{\rm s}  \;s.t.\; v_i =v_j= v_k =v_\ell =1\;and\; v_m=0\ \forall m\in N_{\rm s}\setminus \{i,j,k,\ell\} \text{ or } m\in N_{\rm d} \}\\
\mathcal{D}_{10}
&=\{v~|~ \exists i\in N_{\rm s}  \;s.t.\; v_i=2\;and\; v_j=0\ \forall j\in N_{\rm s}\setminus \{i\} \text{ or } j\in N_{\rm d} \}\\
\mathcal{D}_{11}
&=\{v~|~ \exists i\in N_{\rm s},\exists j\in N_{\rm d}  \;s.t.\; v_i =2 \; and \; v_j=1\;and\; v_k=0\ \forall k\in N_{\rm s}\setminus \{i\} \text{ or } k\in N_{\rm d}\setminus \{j\} \}\\
\mathcal{D}_{12}
&=\{v~|~ \exists i\neq j\in N_{\rm s}  \;s.t.\; v_i =v_j=2\;and\; v_k=0\ \forall k\in N_{\rm s}\setminus \{i,j\} \text{ or } k\in N_{\rm d} \}\\
\end{aligned}
\end{equation*}
\begin{remark}
Again, a careful analysis of the coefficient map $\mathcal C$ reveals that terms of the form $\mathcal{C}(k,l,m,c,d,e) \,  t_i^c  t_k^a  t_l^d t_m^e$ can only give rise to monomials of the form $wxyz$ or $x^2y^2$, where $w,x,y,z$ are distinct single excitation amplitudes. Thus, $\mathcal D_9$ and $\mathcal D_{12}$ are the only contributions of quartic terms.
\end{remark}
In order to find the extremal vertices, we use the following lemmata.
\begin{lemma}
The following set of inclusions hold
\begin{align*}
\mathcal{D}_1 &\subset {\rm conv} (\mathcal{D}_3 \cup \{0\}) &
\mathcal{D}_2 &\subset {\rm conv} (\mathcal{D}_3 )&
\mathcal{D}_4 &\subset {\rm conv} (\mathcal{D}_{10} \cup \{0\}\\
\mathcal{D}_5 &\subset {\rm conv} (\mathcal{D}_{10})&
\mathcal{D}_6 &\subset {\rm conv} (\mathcal{D}_3 \cup \mathcal{D}_{10})&
\mathcal{D}_7 &\subset {\rm conv} (\mathcal{D}_{3} \cup \mathcal{D}_{12})\\
\mathcal{D}_8 &\subset {\rm conv} (\mathcal{D}_{10} \cup \mathcal{D}_{12})&
\mathcal{D}_9 &\subset {\rm conv} (\mathcal{D}_{12})&
\mathcal{S}_2&\subset {\rm conv} (\mathcal{D}_3 \cup \{0\}) \\
\mathcal{S}_4&\subset {\rm conv} (\mathcal{D}_3 \cup \mathcal{D}_{10}) &
\mathcal{S}_6 &= \mathcal{D}_{10} &
\mathcal{S}_7&\subset {\rm conv} (\mathcal{D}_{10} \cup \mathcal{D}_{12}). 
\end{align*}
\end{lemma}
Therefore, we can summarize the vertex description.
\begin{theorem}
The vertices described by $\{0\},\mathcal D_3, \mathcal D_{10},\mathcal{D}_{11}, \mathcal D_{12}  $ are the extremal vertices of the convex hull of $\mathcal D_1,\dots, \mathcal D_{12}$.
\end{theorem}
We will then define the surrogate polytope $N_D$ as follows,$$N_D\colon= {\rm conv}(\{0\}\cup\mathcal{D}_3\cup\mathcal{D}_{10}\cup\mathcal{D}_{11}\cup\mathcal{D}_{12}).$$
By construction, the polytopes $N_S$ and $N_D$ contain all Newton polytopes of the single and double projections, respectively. 
Then, the mixed volume of the Newton polytopes can be bounded from above and below as
\begin{align*}
    \frac{2^n}{n!}&= {\rm Vol} (S_n(2))\leq {\rm mixed Volume}(N_{S_i},N_{D_i})\leq {\rm mixed Volume}(\underbrace{N_S}_{\rm n_s-times},\underbrace{N_D}_{\rm n_d-times})\\
    &\leq
{\rm Vol}({\rm conv} \{\{0\},\mathcal{D}_3, \mathcal{D}_{10}, \mathcal{D}_{11}, \mathcal{D}_{12}\})
\leq {\rm Vol} (S_n(4))= \frac{4^n}{n!}
\end{align*}

\subsection{The Coupled Cluster equations as Cubic Polynomial}

The asymmetry rules allow us to exploit some symmetries of the coefficient tensor $\mathcal C$, i.e. all terms of the form $\sum C(i,j,c,d)t_i^ct_j^d$ can be written as $\sum C(i,j,c,d)(t_i^ct_j^d-t_i^dt_j^c)$. Note we have as many disconnected double excitations $t_i^at_j^b$ as doubles. Formally, we define an index map $\iota$ that flattens the tuple $(i,j,a,b)$. Observe, that $\iota$ is invertible. As an abuse of notation, we use the inverse of $\iota$ also for single excitations. 
Note that since $i$, $j$ are occupied indices, and $a$, $b$ are virtual indices, there are exactly $n_d$ auxiliary indices that are obtained by $\iota$.  Thus, we can introduce new variables $$x_k =\begin{cases}
    t_i^a\quad \text{with } i,a = \iota^{-1}(k) \text{ for } k\leq n_s\\
    t_{ij}^{ab}\quad \text{with } i,a,j,b = \iota^{-1}(k) \text{ for } n_s<k\leq n_d\\
    t_i^at_j^b-t_i^bt_j^a\quad \text{with } i,a,j,b = \iota^{-1}(k-n_d) \text{ for } n_d<k.
\end{cases} $$
We then define the variety $\mathcal{A}$ as
\begin{equation*}
\mathcal{A} = \{x\in \mathbb F^{n_s+2n_d}| x_k - x_i^ax_j^b+x_i^bx_j^a= 0,\, \forall k = \iota (i,j,a,b)+n_s+n_d\}.
\end{equation*}
Using again the Slater--Condon rules and the antisymmetry of the Slater determinants we get
\begin{equation}
\begin{aligned}
\langle \Phi_{ij}^{ab}, HT^4\Phi_0\rangle 
&=
\sum_{klmn,cdef} t_k^c t_l^d t_m^e t_n^f
\langle \Phi_{ij}^{ab}, H X_k^c X_l^d X_m^e X_n^f \Phi_0\rangle \\
&= 
\sum_{k<l,c<d}(t_i^at_j^b-t_i^bt_j^a)(t_k^ct_l^d-t_k^dt_l^c)\langle \Phi_{ij}^{ab}, H \Phi_{ijkl}^{abcd}\rangle
\end{aligned}
\end{equation}
and, thus, on this variety, every projected equation $f_\mu$ can be reduced to a cubic polynomial.
\begin{theorem}[CC reduction theorem]
\label{th:quadQe}
The CCSD equations can be written as a sparse system of cubic equations in $n_s+2n_d$ variables, i.e.,
\begin{equation}
f_\mu(\Pi {\bf x}) 
= 
\sum_{q,r} h^\mu_{q,r}x_qx_r
+
\sum_{p,q,r} \eta^\mu_{pq,r}x_px_qx_r
\end{equation}
on the variety $\mathcal{A}$.
The matrix $h^\mu$ and the tensor $\eta^\mu$ together with their sparsity patterns are derived in the appendix. Furthermore, $\eta^\mu=0$ for single excitation indices $\mu.$
\end{theorem}
The proof of~\cref{th:quadQe} is performed in~\cref{sec:ccWorkingEQ}. 
We explicitly work out the expression of $h^\mu$ for the example of two electrons in four spin orbitals in Eq.~\eqref{eq_h_matrix}.
With this cubic form of the CCSD equations at hand, we can establish the following bound for the roots.
\begin{corollary}
\label{cor:LowerBezout}
The number of roots to the CCSD equations can be bounded by
\begin{equation}
\mathcal{N}\leq 2^{n_s+n_d}3^{n_d}.
\end{equation}
\end{corollary}
\begin{remark}
Note, that this bound improves the na\"ive Bezout bound only for $N=2$. The more important observations that can be made in Theorem \ref{th:quadQe} lie within the sparsity pattern. For example, the 3-tensor $\eta^\mu$ has only four non-zeros slices, each of which is sparse itself. Another interesting observation is that the coefficients of the double excitation amplitudes and the corresponding ``extended'' doubles are the same. Exploiting these structures numerically has the potential to accelerate root calculations tremendously.
\end{remark}

\section{Two Electrons in Four Spin-Orbitals}
\label{sec:twoInFour}

In this section, we scrutinize the algebraic and geometric properties and structures of the two electrons in four spin-orbitals case and perform benchmark calculations revealing the highly intricate structures of CC theory. 
In this model system, we find that the $N$-particle Hilbert space (here the two-particle Hilbert space), is spanned by six Slater determinants. 
Without loss of generality, we denote $\Phi_0 = \Phi[1,2]$ the reference state, and the respectively excited Slater determinants are
\begin{equation}
\left\lbrace
\Phi_2^3, \Phi_2^4, \Phi_1^3, \Phi_1^4, \Phi_{12}^{34}
\right\rbrace.
\end{equation}

The aim of this section is a mathematical exposition of the CC equations and their root structure. 
We therefore do not include spin-symmetries here, e.g., spin-restricted or spin-unrestricted formulations.  
When modeling physical systems these symmetries have to be discussed which results in a reduction of the dimensionality.
Therefore, the model presented here describes a mathematically abstracted problem without a direct physical interpretation.

\subsection{Algebraic Description}
We exemplify equation \eqref{eq:singles-ccsdEQ} for the polynomial corresponding to the projection onto $\Phi_2^3$ in the minimal working example of two electrons in four spin-orbitals and derive the Newton polytope on the fly.

The corresponding polynomial is given by
\begin{equation}
\langle \Phi_1^3, e^{-T}H e^T \Phi_0  \rangle.
\end{equation}

We shall now expand the similarity transformed Hamiltonian and inspect the individual terms; we will mainly focus on the degree of the individual monomials and as we are seeking a bound to $\mathcal N$ that is based on the BKK theorem, i.e., we require to characterize the Newton polytope.
The first part is given by the constant contribution 
\begin{equation}
\begin{aligned}
\langle \Phi_1^3, H \Phi_0  \rangle = C,
\end{aligned}
\end{equation}
which corresponds to the origin $(0,0,0,0,0)$.
The second part is given by the singly nested commutator contribution, i.e. 
\begin{equation*}
\begin{aligned}
\langle \Phi_1^3, [H,T] \Phi_0  \rangle &=
 t_2^3 \langle \Phi_1^3, H \Phi_2^3  \rangle
+t_2^4 \langle \Phi_1^3, H \Phi_2^4  \rangle
+t_1^3 \langle \Phi_1^3, H \Phi_1^3  \rangle\\
&\quad+t_1^4 \langle \Phi_1^3, H \Phi_1^4  \rangle
+t_{12}^{34} \langle \Phi_1^3, H \Phi_{12}^{34}  \rangle
-t_1^3 \langle \Phi_0, H \Phi_0  \rangle,
\end{aligned}
\end{equation*}
which corresponds to the vertices 
\begin{equation*}
\{
(1,0,0,0,0),
(0,1,0,0,0),
(0,0,1,0,0),
(0,0,0,1,0),
(0,0,0,0,1)
\}.
\end{equation*}
The third part is given by the doubly nested commutator contribution, i.e. 
\begin{equation*}
\begin{aligned}
\langle \Phi_1^3, [[H,T],T] \Phi_0  \rangle &=
  \langle \Phi_1^3, HT^2 \Phi_0 \rangle
-2\langle \Phi_1^3, THT \Phi_0 \rangle
+ \langle \Phi_1^3, T^2H \Phi_0 \rangle\\
&=t_2^3t_1^4  \langle \Phi_1^3, H \Phi_{12}^{34}  \rangle 
- t_1^3t_2^4 \langle \Phi_1^3, H \Phi_{12}^{34}  \rangle
-2t_1^3 \langle \Phi_0, HT \Phi_0  \rangle,
\end{aligned}
\end{equation*}
which corresponds to the vertices 
\begin{equation*}
\{
(1,0,0,1,0),
(0,1,1,0,0),
(1,0,1,0,0),
(0,1,1,0,0),
(0,0,2,0,0),
(0,0,1,1,0),
(0,0,1,0,1)
\}.
\end{equation*}
The fourth part is given by the triply nested commutator contribution, i.e.
\begin{equation*}
\begin{aligned}
\langle \Phi_1^3, [[[H,T],T],T] \Phi_0  \rangle 
&=
 \langle \Phi_1^3, HT^3 \Phi_0 \rangle
-3\langle \Phi_1^3, (THT^2 + T^2HT) \Phi_0 \rangle
+ \langle \Phi_1^3, T^3H \Phi_0 \rangle\\
&=-3t_1^3  \langle \Phi_0, HT^2 \Phi_0  \rangle\\
&=-3t_1^3(t_2^3t_1^4  \langle \Phi_0, H \Phi_{12}^{34}  \rangle 
- t_1^3t_2^4 \langle \Phi_0,  H \Phi_{12}^{34}  \rangle),
\end{aligned}
\end{equation*}
which corresponds to the vertices 
\begin{equation}
\{
(1,0,1,1,0),
(0,1,2,0,0)
\}.
\end{equation}
The quadruple commutator gives no contributions due to the Slater--Condon rules. The first polynomial can then be summarized as
\begin{equation*}
\begin{aligned}
f_1(t) &= 
\langle \Phi_1^3, e^{-T}H e^T \Phi_0  \rangle\\
&=
C + 
t_2^3  \langle \Phi_1^3, H \Phi_2^3  \rangle
+t_2^4  \langle \Phi_1^3, H \Phi_2^4  \rangle 
+t_1^3  \langle \Phi_1^3, H \Phi_1^3  \rangle 
+t_1^4  \langle \Phi_1^3, H \Phi_1^4  \rangle 
-t_1^3 \langle \Phi_0, H \Phi0 \rangle \\
&\quad 
-t_1^3 \big( 
 t_2^3 \langle \Phi_0, H \Phi_2^3  \rangle 
+t_2^4 \langle \Phi_0, H \Phi_2^4  \rangle
+t_1^3 \langle \Phi_0, H \Phi_1^3  \rangle 
+t_1^4 \langle \Phi_0, H \Phi_1^4  \rangle 
+t_{12}^{34}  \langle \Phi_0, H \Phi_{12}^{34}  \rangle \big)
\\
&\quad 
+\big(t_{12}^{34} 
-t_2^3t_1^4 
+t_1^3t_2^4 \big) \langle \Phi_1^3, H \Phi_{12}^{34}  \rangle 
-t_1^3 \big( t_1^3t_2^4   
- t_2^3t_1^4 \big) \langle \Phi_0, H \Phi_{12}^{34}  \rangle 
\end{aligned}
\end{equation*}
and the corresponding Newton polytope is given by
\begin{equation}
\begin{aligned}
{\rm New}_1^3
={\rm conv} \big(&
  (0, 0, 0, 0, 0),~
  (1, 0, 0, 0, 0),~
  (0, 1, 0, 0, 0),~
  (0, 0, 0, 1, 0),~\\
  &(0, 0, 0, 0, 1),~
  (0, 1, 1, 0, 0),~
  (1, 0, 1, 0, 0),~
  (0, 1, 1, 0, 0),~\\
  &(0, 0, 2, 0, 0),~
  (0, 0, 1, 0, 1),~
  (1, 0, 1, 1, 0),~
  (0, 1, 2, 0, 0)  
\big).
\end{aligned}
\end{equation}
The polynomials and polytopes corresponding to the remaining projections onto singly and doubly excited Slater determinants can be obtained correspondingly, see~\cref{app:2in4}. 

Using polymake, we calculated the BKK bound on the number of roots via the mixed volume which is equal to $50$. Comparing this with the different B\'ezout bounds derived above puts their quality into perspective
\begin{equation}
50 < 2^53 = 96 < 3^4 4 = 324 .
\end{equation}
Using polymake, we can also obtain the description of the facets of ${\rm New}_1^3$ in homogeneous coordinates:

\begin{align*}
\left\lbrace
\begin{pmatrix}
0\\ 0\\ 0\\ 0\\ 0\\ 1
\end{pmatrix},
\begin{pmatrix}
0\\ 0\\ 0\\ 0\\ 1\\ 0 
\end{pmatrix},
\begin{pmatrix}
0\\ 0\\ 0\\ 1\\ 0\\ 0 
\end{pmatrix},
\begin{pmatrix}
0\\ 0\\ 1\\ 0\\ 0\\ 0
\end{pmatrix},
\begin{pmatrix}
0\\ 1\\ 0\\ 0\\ 0\\ 0
\end{pmatrix},
\begin{pmatrix}
1\\ 0\\ -1\\ 0\\ -1\\ -1
\end{pmatrix},
\begin{pmatrix}
1\\ -1\\ -1\\ 0\\ 0\\ -1 
\end{pmatrix},
\begin{pmatrix}
2\\ -1\\ 0\\ -1\\ 0\\ -1
\end{pmatrix},
\begin{pmatrix}
2\\ 0\\ 0\\ -1\\ -1\\ -1
\end{pmatrix}
\right\rbrace 
\end{align*}
This description reflects Pauli's exclusion principle in the Newton polytope. More precisely, in the vectors that do not correspond to canonical basis vectors, the distributions of the `$-1$' entries follow a particular selection procedure. We first note that the last entry which corresponds to the double excitation $t_{12}^{34}$ is set to $-1$. Then, we see that the remaining entries are chosen corresponding to the physically allowed double excitations that correspond to $t_{12}^{34}$, i.e., pairs $t_i^a,t_j^b$ that form $t_{12}^{34}$ have one entry equal to $-1$ and one entry equal to $0$. This procedure results in the particular distribution of $-1$ and $0$ depicted in the facet description of ${\rm New}_1^3$ in homogeneous coordinates.

\begin{figure}[h!]
    \centering
    \includegraphics[width=0.5\textwidth]{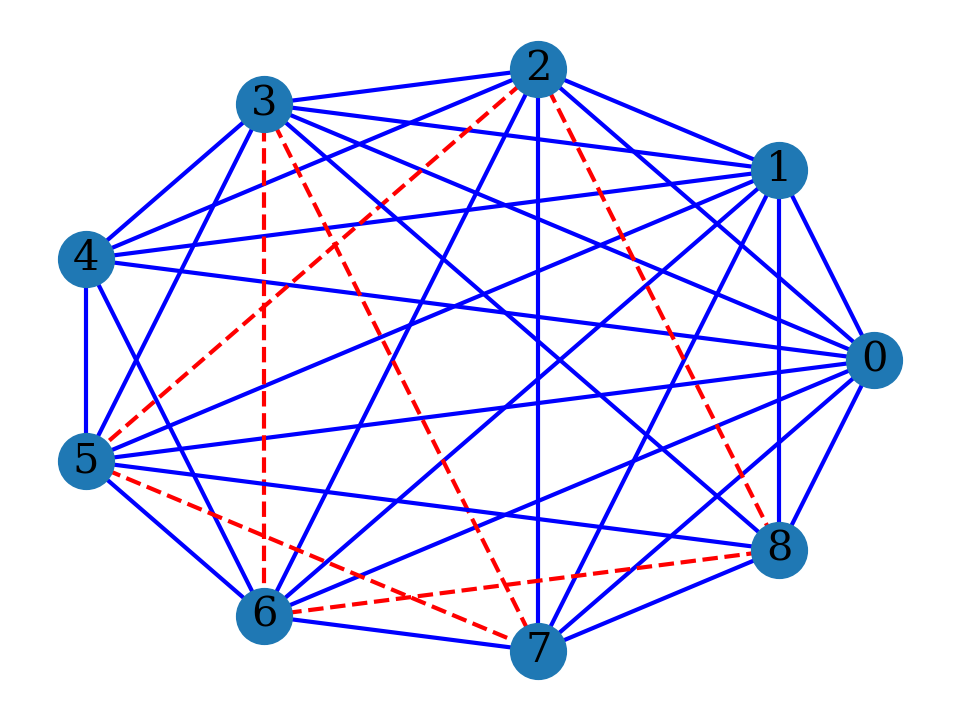}
    \caption{Graph that shows the dimensionality of the intersection of facets of the third Newton polytope. If two facets have a four-dimensional intersection they are connected with a solid blue line, and if two facets have a three-dimensional intersection they are connected with a dotted red line.}
    \label{fig:graph_polytope_1}
\end{figure}

Polymake also allows computing the f-vector of the Newton polytope, for New$_1^3$ which yields $(12, 33, 42, 28, 9 )$. 
In order to obtain a more geometric interpretation of the Newton polytopes we generate a graph that relates the individual facets. To that end, we enumerate the facets from zero to seven and connect two facets if they intersect. It is worth noticing that in the case of two electrons in four spin-orbitals, the intersection of two facets is either four- or three-dimensional, see~\cref{fig:graph_polytope_1}. Moreover, we see that the majority of facets intersect in four dimensions. 

Since CCSD is equivalent to FCI in this particular example, we expect to find exactly six roots, which numerical simulations do. This drastic reduction in the number of roots, compared to the generic case, might originate from the highly structured coefficient tensor of the polynomial system.
As shown above, the polynomial system can be written as a quadratic polynomial in $$1,t_2^3,t_2^4, t_1^3,t_1^4,t_{12}^{34}, t_1^3t_2^4-t_2^3t_1^4.$$ 
Due to the system's size, this can be confirmed by direct computation for the two electrons in four spin-orbitals case. 
Indeed, the CC equations can be written as 
\begin{align*}
p_i = 
v^T h^i v
\end{align*}
where $v\in A$ and
\begin{equation}
A = \{(x_0,x_1,x_2,x_3,x_4,x_5,x_6):x_0=1 \text{ and } x_6-x_2x_3+x_1x_4=0\}.
\end{equation}
We exemplify the coefficient matrix for the index $i=1$, here, 
\begin{equation}
\label{eq_h_matrix}
h^1  = \begin{pmatrix}
    h^1_{1,0}&h^1_{1,1}&h^1_{1,2}&h^1_{1,3}&h^1_{1,4}&h^1_{1,5}&-h^1_{1,5}\\
    -h^1_{0,0}&-h^1_{0,1}&-h^1_{0,2}&-h^1_{0,3}&-h^1_{0,4}&-h^1_{0,5}&h^1_{0,5}\\
    0&0&0&0&0&0&0\\
    0&0&0&0&0&0&0\\
    0&0&0&0&0&0&0\\
    0&0&0&0&0&0&0\\
    0&0&0&0&0&0&0
    \end{pmatrix}.
\end{equation}
The other coefficient matrices are listed in the~\cref{app:2in4}. 

We emphasize that the B\'ezout bound that is closest to the BKK bound arises from this cubic representation, i.e., $\mathcal{N}\leq 2^53 = 96$. 
A deeper investigation regarding the precise structures of the coefficient matrices $h^i$ can be found in the appendix. 

Since CCSD, under certain assumptions, corresponds to the exact diagonalization of the underlying chemical system, we shall now consider the truncated version of CCSD which is merely CCS. 
One could alternatively consider CCD, however, for two electrons in four spin orbitals, the CCD system consists of only one equation. 
We note that CCS corresponds to orbital rotations (see Thouless theorem~\cite{thouless1960stability}), yet, from the algebraic perspective this experiment reveals the highly complicated structure of the roots and the corresponding non-trivial effect on the computed energies. 
To that end, we perform the following experiment. We define an arbitrary Hamiltonian through
\begin{equation}
H(\varepsilon)
=
V(\varepsilon)^\top {\rm diag}(-\frac{11}{12}, -\frac{9}{12}, -\frac{4}{6}, -\frac{3}{6}, -\frac{2}{6}, -\frac{1}{6}) V(\varepsilon),
\end{equation}
where $V(\varepsilon)$ is the orthonormalized part (via QR-decomposition) of $\tilde{V}(\varepsilon)$ which is defined as 
\begin{equation}
\tilde{V}(\varepsilon)
=
\begin{pmatrix}
1 & 1+\varepsilon & 1+\varepsilon & 1+\varepsilon & 1+\varepsilon & 1+\varepsilon\\
0.05 & 1 & 0 & 0 & 0 & 0\\
0.05 & 0 & 1 & 0 & 0 & 0\\
0.05 & 0 & 0 & 1 & 0 & 0\\
0.05 & 0 & 0 & 0 & 1 & 0\\
0 & \varepsilon & \varepsilon & \varepsilon & \varepsilon & \varepsilon\\
\end{pmatrix}.
\end{equation}

\begin{figure}
    \centering
    \includegraphics[width=0.9\textwidth]{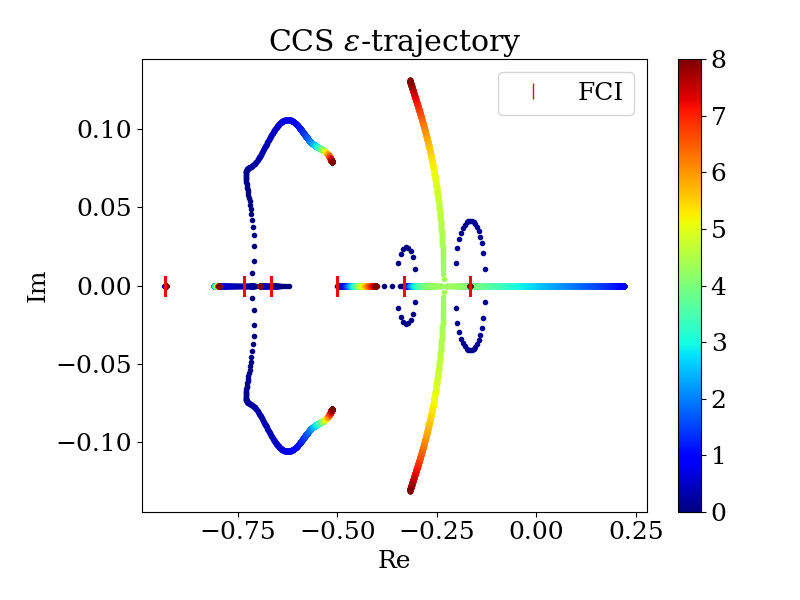}
    \caption{The $\varepsilon$-energy trajectory of CCS, where $\varepsilon$ was varied between zero and eight. The value of $\varepsilon$ is indicated through different colors, the six FCI energies are marked with red vertical lines.}
    \label{fig:ccs_trajectory}
\end{figure}

The parameter $\varepsilon$ affects the overlap of the eigenstates with the reference state, which is $(1,0,0,0,0,0)^\top$. 
Tuning the parameter $\varepsilon$ allows us to track the roots of the polynomial system numerically. 
The parameter $\varepsilon$ steers the overlap of the eigenstates with the doubly excited Slater determinant, except for the lowest-lying state.
Hence, the lowest lying state does not change with the parameter $\varepsilon$, see~\cref{fig:ccs_trajectory}.
However, this is not the case for the remaining eigenstates. As we change the parameter $\varepsilon$ the overlap of the eigenstates with the reference state changes, yielding poor representability in the exponential form. This poor representation has severe effects on the computed energies, see~\cref{fig:ccs_trajectory}.
In fact, not only are the computed CC energies inaccurate but also complex-valued, which means that they are unphysical.

Moreover, the number of roots changes as a function of $\varepsilon$. 
In Table~\ref{tab:Roots} we provide the number of roots for selected values of $\varepsilon$. 
We emphasize that the minimal number of roots is 5, and the maximal number of roots is 9.
\begin{table}[h!]
    \centering
    \begin{tabular}{c|ccccccc}
         $\varepsilon$ &  0 & 0.05  & 0.1   & 0.2   & 1.5   & 3.0   & 6.5\\
         \hline
         $\mathcal{N}$ &  5 & 9     & 8     & 9     & 9     & 9     & 8
    \end{tabular}
    \caption{Number of roots as for different values of $\varepsilon$.}
    \label{tab:Roots}
\end{table}

We emphasize that employing CCSD, in this case, recovers all eigenstates correctly unless one of the eigenstates is orthogonal to the chosen reference state.
Noting that the overlap with the reference state is deeply rooted in the local analysis of CC theory, this experiment indicates that the assumptions made in~\cite{schneider2009analysis,laestadius2019coupled}, also provide sufficient conditions for the CC solution to provide a physical description of the system.
Hence, this experiment is an alternative confirmation that these assumptions can be used to derive a useful {\it a posteriori} diagnostic for single reference CC calculations~\cite{faulstich2023s}.

Algebraically, it is remarkable that the number of roots to the polynomial system (see Table~\ref{tab:Roots}) is dramatically lower than the BKK number (i.e. 50 for the considered system), which turns out to be sharp for a generic polynomial system of this size.  
We, therefore, believe that in order to obtain improved insight into the number and structure of roots of the CCSD equations, new algebraic techniques have to be developed that incorporate the eminent structure of the quantum many-body problem. 

We note that the above example is a mere model problem used to illustrate the root behavior of CC equations in a general setting. To extend our study to truly physical systems (as has been considered in e.g.~\cite{kowalski1998towards,piecuch2000search,kowalski2000complete,kowalski2000complete2}), a number of careful considerations have to be made. However, these considerations are based on physical concepts that go beyond the scope of this manuscript and we, therefore, refrain from including them here.

\section{Three Electrons in Six Spin-Orbitals}
\label{sec:threeInSix}

The rapid increase in the size of the considered polynomial system is well demonstrated by the fact that for the next bigger system, i.e., three electrons in six spin-orbitals, the system already supersedes the abilities of state-of-the-art algebraic geometry software when used in an out-of-the-box fashion.
In order to numerically investigate the problem further, we perform a root search using Newton's method for the energetically lowest state (i.e. the ground state). 
This procedure closely resembles the practical use of CC theory, where the cluster amplitudes corresponding to the ground state solution are approximated using a quasi-Newton method. 
An important unanswered question that arises in this regard is the convergence behavior depending on the initial guess. 
In practical quantum chemical simulations, it is common to use an initial guess that stems from second-order M\o{}ller–Plesset perturbation theory (MP2)~\cite{helgaker2014molecular}. 

In order to understand the behavior of Newton's methods and the root structure around the lowest lying eigenstate, we choose from a sequence of randomly generated, symmetric Hamiltonians fulfilling $\langle \Phi_0, H\Phi_{123}^{456}\rangle =0$ an instance, such that for one of the eigenvectors $v$ we have $\langle \Phi_0, v \rangle \approx 1$ and the contribution of the triple CC amplitude $t_{123}^{456}\ll 1$. Note that for the three electrons in six spin-orbitals case, CCSDT is exact and there is exactly one triple CC amplitude $t_{123}^{456}$, i.e. there is a bijective map between the FCI and CCSDT description. We want to recover the CC amplitudes of this eigenvector $v$.

\begin{figure}[tbhp]
\centering
\subfloat[]{\label{fig:newtonfit_01}\includegraphics[width=0.45\textwidth]{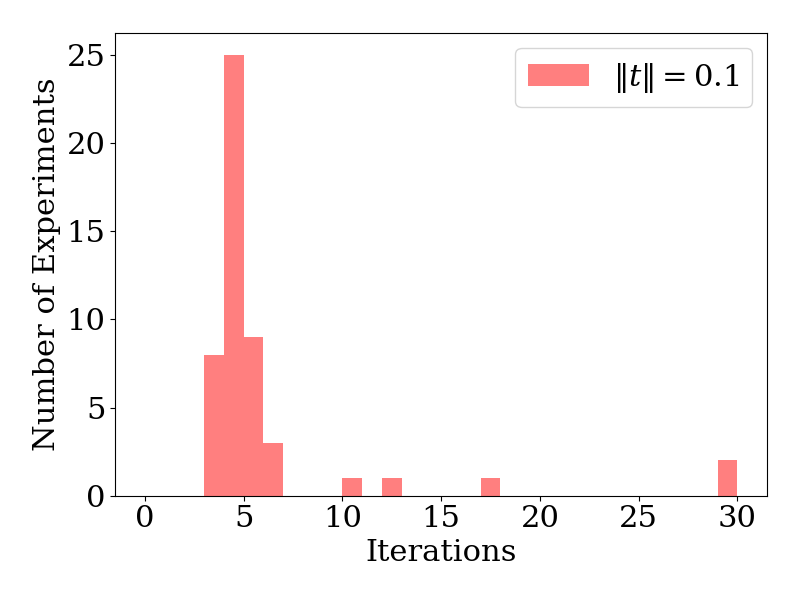}}
\subfloat[]{\label{fig:newtonfit_06}\includegraphics[width=0.45\textwidth]{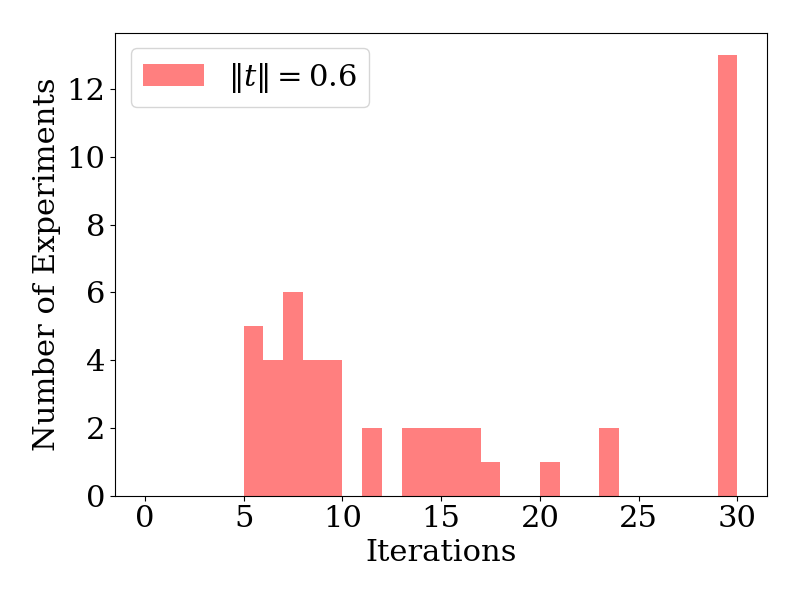}}
\linebreak
\subfloat[]{\label{fig:newtonfit_1}\includegraphics[width=0.45\textwidth]{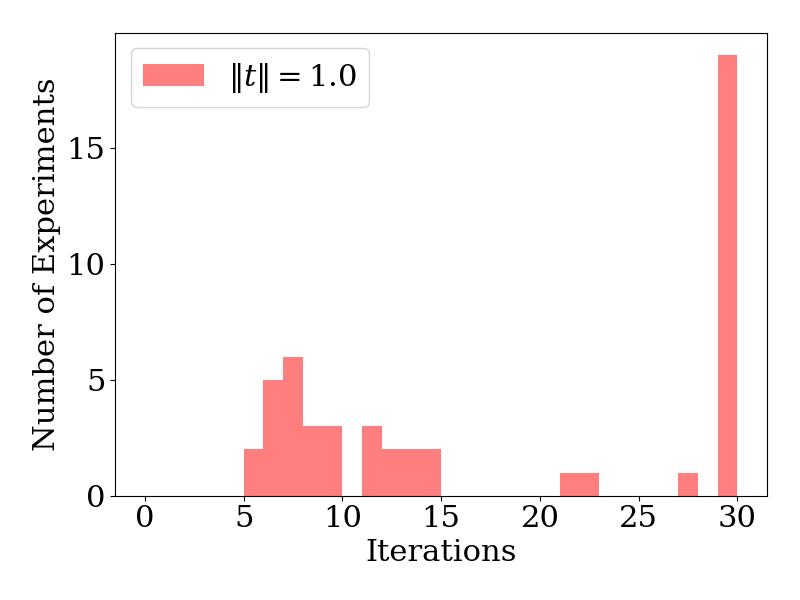}}
\subfloat[]{\label{fig:newtonfit_12}\includegraphics[width=0.45\textwidth]{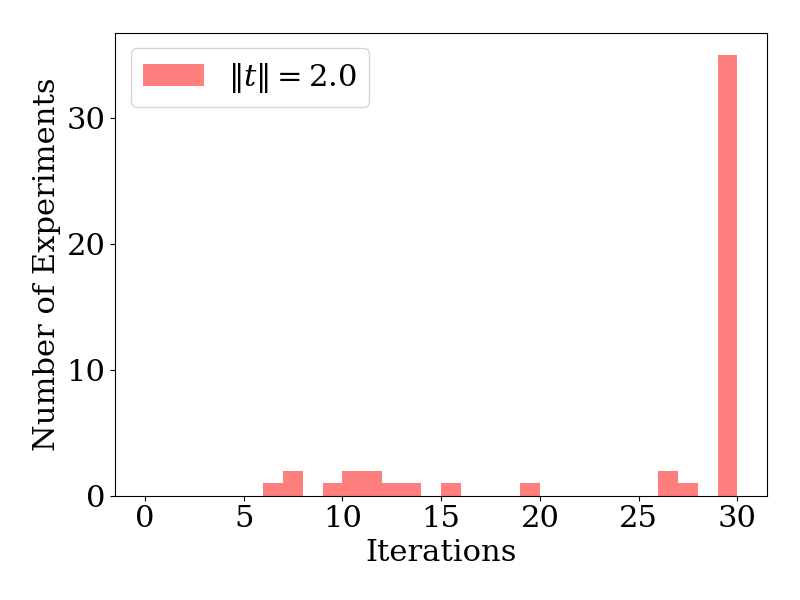}}
\caption{
Number of Newton iterations for randomly perturbed initial guess. The imposed perturbations are of size $0.1$, $0.6$, $1.0$, $2.0$ in~\ref{fig:newtonfit_01},~\ref{fig:newtonfit_06},~\ref{fig:newtonfit_1},~\ref{fig:newtonfit_12}, respectively.}

\label{fig:newton}
\end{figure}

In order to investigate the sensitivity of the Newton method in CC theory we perform Newton optimizations for a set of random perturbations of the exact CC amplitudes as the initial value. For each experiment, we draw 50 normally distributed perturbations and normalize them to be of the desired size. We then perform a Newton optimization with at most 30 iterations, which is sufficient for the considered experimental setup.
It turns out that in CC theory, the sensitivity of Newton's method is rather severe, see~\cref{fig:newton}. We emphasize that the CCSD amplitude vector in this example consists of 18 values, i.e., in a perturbation that has norm 1 each entry is on average of size $1/18\approx 0.056$. The \cref{fig:newton} clearly shows the decrease of successfully converged Newton optimizations as we increase the norm of the perturbation (at 30 iterations we deem the optimization to have failed). 

\begin{figure}[tbhp]
    \centering
    \includegraphics[width=0.6\textwidth]{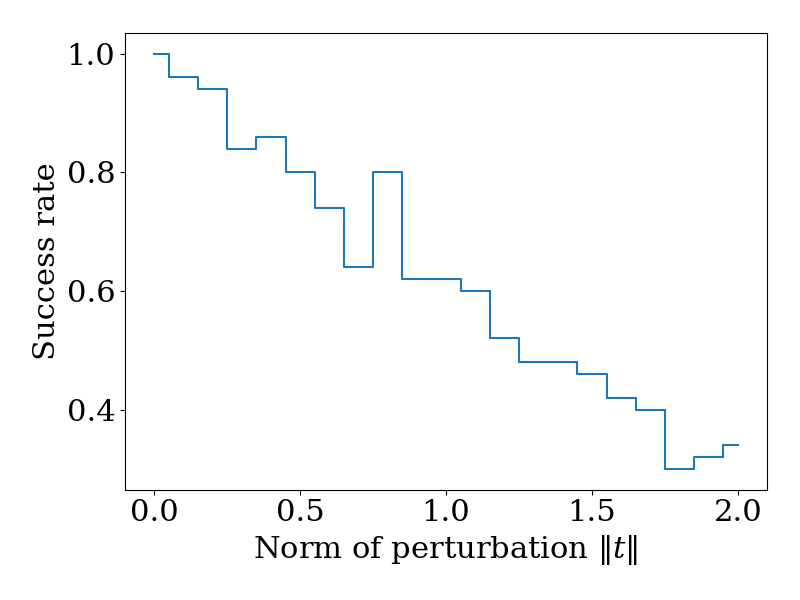}
    \caption{Success rate of Newton's method as a function of the norm of the perturbation $\Vert t\Vert$. }
    \label{fig:success_rate}
\end{figure}

In~\cref{fig:success_rate} we show the success rate as a function of the norm of the perturbation. As expected, this clearly shows a decreasing trend. We emphasize that in practical simulations, the initial guess is not chosen at random but corresponds to the MP2 result. Yet, our experiments show that for the Newton optimization to be on average successful, we require the initial guess to be of high quality. This again supports the assumptions made in~\cite{schneider2009analysis}.
We moreover note that if Newton's method converges successfully, it mostly converges towards the originally perturbed solution. However, in a few cases, we observe that Newton's method successfully converges but towards a different solution. It is worth mentioning that this solution was not physical, i.e., the solution yielded an energy value of 51 which is outside the spectrum of the considered Hamiltonian (spec(H)$= [\![ 1,20 ]\!]$).

\section{Conclusion and Outlook}
In this article, we investigated the coupled cluster (CC) formalism from an algebraic geometry perspective. 
We moreover provided a derivation of the CC equations that correspond to the CC method with excitation rank truncation at the level of single and double excitations, called CCSD.
This particular method is of utmost importance to the computational chemistry community since it is a widely used high-accuracy method, and as such is commonly used for benchmarking purposes.\\

We approached the study of the CCSD equations, which is a system of polynomial equations, from an algebraic perspective.
Since the number of roots is of general interest we investigated different approaches. 
By mere inspection of the CCSD equations, we find that the most naive B\'ezout bound of $\mathcal{O}(4^{(n_s+n_d)})$ can be straightforwardly improved to $\mathcal{O}(3^{n_s}4^{n_d})$. 
In order to establish a more profound bound, we investigated the use of the BKK theorem. 
To that end, we analyzed the Newton polytopes that correspond to the system of polynomial equations.  
Although the polynomials are at most of degree four, they have an intricate structure that impedes volume estimations from the vertex description. We establish surrogate Newton polytopes that show a simpler structure, which may simplify the path to a bound based on the BKK theorem.\\
\indent
We moreover proved a novel description of the CCSD equations as a cubic system.  
This description provides a fundamentally new perspective to the CCSD equations, and yields a further improvement to the bound of zeros, namely, $\mathcal{O}(2^{n_s+2n_d})$.\\

We then performed numerical investigations for two non-trivial systems. The first is the two electrons in four spin-orbitals system.
Here we calculated all Newton polytopes in the vertex description in order to provide an insight into their underlying structure. 
Using the Software package polymake we were able to compute the corresponding facet description. This revealed that Pauli's exclusion principle is deeply embedded into the structure of the polytopes as certain occupation rules of the orbitals have to be obeyed. 
We moreover investigated the geometrical properties of the five-dimensional polytope using a graph structure to reveal the orientation of the facets relative to each other. 
We found that if facets of this polytope have a non-empty intersection, this intersection is mostly of dimension four, very few have a three-dimensional intersection, but lower dimensional intersections did not occur. 
In order to understand the algebraic structures, especially in the case of truncated CC theory, we constructed a model system that can be steered to produce unphysical solutions in the truncated CC approach.
This experiment moreover related our algebraic investigations to previous functional analysis results and revealed that the assumptions made in~\cite{schneider2009analysis,laestadius2019coupled} not only apply to the quality of the approximation provided by CC theory but also related to the physicality of the described solution.  \\

To conclude, this article shows that the polynomial system corresponding to the CCSD equations is highly structured and attains significantly fewer roots than estimated with state-of-the-art techniques in algebraic geometry. 
Given the impact of computational chemistry and its contemporary relevance, it appears worthwhile to develop new algebraic techniques that take the apparent structures into account in order to provide more detailed information on the number and structure of roots.   
We expect that the representation of the polynomials as sparse cubic forms on a quadratic algebraic variety together with the banded structure of the coefficient matrix opens new lines of thought on tackling highly structured polynomial systems.

\section{Acknowledgements}
This work was partially supported by the Air Force Office of Scientific Research under award
number FA9550-18-1-0095 and by the Simons Targeted Grants in Mathematics and Physical Sciences on Moir\'e Materials Magic (F.M.F.). We also thank Jan Draisma for his initial comments that have improved the manuscript's accessibility and Prof.~Bernd Sturmfels, and his graduate students Yelena Mandelshtam and Yulia Alexandr as well as Prof. Timo de Wolff and his group for fruitful discussions. 
Last but certainly not least, we wish to acknowledge that this work was predominantly developed during a {\it Research in Pairs} stay at the Oberwolfach Research Institute, and we thank its members for providing an outstanding, unique, and inspiring atmosphere.

\bibliographystyle{IEEEtranS}
\bibliography{lib.bib}

\newpage
\section{Appendix}
\label{sec:appendix}

\subsection{Functional analytic framework of the electronic Schr\"odinger equation}
\label{app:FuncAna}

The governing equation for the motion of electrons in a molecular system employing the adiabatic---or Born--Oppenheimer---approximation~\cite{born1927quantentheorie,Born1954} is the {\it electronic Schr\"odinger equation}~\cite{helgaker2014molecular}.
For a system of $N$ electrons and $M$ nuclei, each of charge $Z_j$, the electronic Hamiltonian in atomic units is formally defined as 
\begin{equation}
\label{eq:ElecSE_app}
H
=
-\frac{1}{2}\sum_{i=1}^N\Delta_{r_i}
-\sum_{i=1}^N\sum_{j=1}^M\frac{Z_j}{\vert r_i-R_j\vert} 
+\sum_{i=1}^N\sum_{j>i}^N\frac{1}{\vert r_i-r_j\vert},
\end{equation}
where $\{r_i\}_{i = 1}^N$ and $\{R_i\}_{i=1}^M$ describe the spatial coordinates of the electrons and the clamped nuclei, respectively, and $\Delta_{r_i}$ is the Laplacian with respect to the $r_i$-spatial coordinate.  
Central to the quantum-mechanical study of such electronic systems are the eigenfunctions of the Hamiltonian in Eq.~\eqref{eq:ElecSE_app}.
We define wavefunctions through the following mapping 
\begin{equation}
\label{eq:Wavefunction_app}
\Psi\colon(\mathbb{R}^3)^N \times \{ \pm 1/2\}^N\to \mathbb{R}
\colon (r,s) \mapsto \Psi(r,s),
\end{equation}
where $s$ denotes the {\it electron spin} (or simply {\it spin}).
Although spin does not appear directly in the electronic Hamiltonian in Eq.~\eqref{eq:ElecSE_app}, it influences the structure of atoms and molecules~\cite{yserentant2010regularity}. 
On the fundamental level of measure theory, which is necessary to characterize the wave function in Eq.~\eqref{eq:Wavefunction_app}, it is straightforward to incorporate the spin variable; in the single-particle case we simply consider the product measure space 
\begin{equation*}
(X,\mathfrak{B}(X),\lambda_{X})
=
\Big(
\mathbb{R}^{3}\times\{\pm 1/2 \},\mathfrak{B}(\mathbb{R}^{3})\otimes \mathfrak{B}(\{\pm 1/2 \}),\lambda\times c
\Big),
\end{equation*}
where $\mathfrak{B}$ denotes the {\it Borel $\sigma$-algebra}, $\lambda$ is the corresponding {\it Lebesgue measure} and $c$ is the {\it counting measure}.
This naturally extends to the N-body case, where we define
\begin{equation*}
(X^N,\mathfrak{B}(X^N),\lambda_{X^N})
=
\Big(
\mathbb{R}^{3N}\times\{\pm 1/2 \}^N,\mathfrak{B}(\mathbb{R}^{3N})\otimes \mathfrak{B}(\{\pm 1/2 \}^N),\lambda\times c
\Big).
\end{equation*}
We will use the short-hand notation and write $\dd\lambda_{X^N}(x) = \dd x$ and similar for the single-particle case. 

In order to provide a precise description of the considered Coulomb Hamiltonian we require a characterization of the considered domain. 
To that end, we consider its weak formulation~\cite{ern2004theory}. 
Following the standard procedure, we can derive the {\it form domain} of the Coulomb Hamiltonian, i.e.,   
\begin{equation*}
\mathcal{D}(H)= H^1(X^N) \cap \bigwedge_{i=1}^NL^2(X).
\end{equation*}
For a more detailed description of this procedure, we refer the interested reader to~\cite{laestadius2019coupled,faulstich2019analysis,yserentant2010regularity,reed1980functional}

\subsection{Exponential Parametrization}
In this section of the appendix, we list proofs and additional results related to excitation operators and exponential parametrization.

\subsubsection{Proofs of~\cref{sec:exponentialparam}}

\begin{proof}{\it [\cref{th:Basis<->operator}]}
Since excitation operators are defined w.r.t.~a chosen reference Slater determinant $\Phi_0$ it follows immediately that $\Phi_0= I \Phi_0$. 
Subsequently, we assume without loss of generality that $\Phi_0=\Phi [1,...,N]$. 
Let $\Phi [p_1,...,p_N]\in \mathcal{B}$ be arbitrary. 
Comparing $[1,...,N]$ to $[p_1,...,p_N]$ we can identify a multi-index $\mu$ describing the indices that have to be changed in $[1,...,N]$ to obtain $[p_1,...,p_N]$. Due to the canonical ordering, the multi-index $\mu$ is unique. 
Then, by definition we obtain $\Phi [p_1,...,p_N]= X_\mu \Phi_0$, hence it follows the claim.
\end{proof}

\begin{proof}{\it [\cref{lemma:1}]}
Note that $\mathcal{H}_{\rm int} \subset \mathcal{H}_N$ is convex since it is an affine subspace, and that the functions $\mathcal{B}_{\rm int}= \{Id + \Phi_\mu~|~ \Phi_\mu \in \mathcal{B}\}$ form an affine-independent set generating $\mathcal{H}_{\rm int} = {\rm conv}(\mathcal{B}_{\rm int})$. It is therefore sufficient to show the theorem for $\Psi \in \mathcal{B}_{\rm int}$. In this case, however, it follows trivially from the fact that there exists a one-to-one relation between $\Phi_\mu$ and $X_\mu$, see Theorem~\ref{th:Basis<->operator}.
\end{proof}

\begin{proof}{\it [\cref{lemma:2}]}
By definition of $ \mathcal{G}$, the map $\Omega$ is trivially surjective. 
Moreover, since ${\rm dim}(\mathfrak{b})= {\rm dim}(\mathcal{G})$, it follows that $\Omega$ is a bijection.
\end{proof}

\begin{proof} {\it [\cref{th:fci}]}
The proof follows directly from Lemma~\ref{lemma:1} and Lemma~\ref{lemma:2}.
\end{proof}

\subsubsection{Additional Results}
\begin{theorem}
\label{th:AbelianLieAlg}
The vector space $\mathfrak{b}$ together with the commutator 
\begin{equation}
[\cdot, \cdot ] : \mathfrak{b} \times \mathfrak{b} ~;~
(S,T) \mapsto S T - TS
\end{equation}
forms an Abelian Lie algebra.
\end{theorem}
\begin{proof} {\it [\cref{th:AbelianLieAlg}] }
It is sufficient to show that the commutator fulfills the Jacobi identity, i.e., 
\begin{equation}
[S,[T,U]]+[T,[U,S]]+[U,[S,T]] = 0,
\end{equation}
for $S,T,U\in \mathfrak{b}$, which is a direct consequence from Lemma~\ref{thm:ExcitationOp} and the bilinearity of the commutator.
\end{proof}

\begin{corollary}
\label{Cor:AbelianLieAlg}
~
\begin{itemize}
    \item[i)] The Lie algebra $\mathfrak{b}$ is an Abelian subalgebra of a nilalgebra, and hence a nilpotent Lie algebra itself.
    \item[ii)] The Lie algebra $\mathfrak{b}$ is solvable.
    \item[iii)] The center $Z(\mathfrak{b}) = \{T\in\mathfrak{b}~|~ [T,S] = 0~\forall S\in\mathfrak{b}\} = \mathfrak{b}$
    \item[iv)] The derivative algebra of $\mathfrak{b}$, denoted $[\mathfrak{b}\mathfrak{b}]$, is equal to zero. 
\end{itemize}
\end{corollary}
\begin{proof} {\it [\cref{Cor:AbelianLieAlg}]}
Property i) is a direct consequence from Engel's theorem, and properties $ii)-iv)$ follow immediately as $\mathfrak{b}$ is Abelian.
\end{proof}

\begin{lemma}
\label{lemma:5}
The set $\tilde{\mathcal{G}} $ is equal to $ \mathcal{G}$.
\end{lemma}
\begin{proof}{\it [\cref{lemma:5}]}
Let $T \in \mathfrak{b}$. Then by definition
\begin{equation}
\exp(T) = Id + P(T) \in \tilde{\mathcal{G}}
\end{equation}
where $P(T)\in \mathfrak{b}$. Since $\mathfrak{b}$ is a vector space, we get $P(T)\in \mathfrak{b}$, and therewith $\tilde{\mathcal{G}} \subseteq \mathcal{G}$.\\ Conversely, let $C\in \mathcal{G}$. Then $C - Id \in \mathfrak{b}$, which implies that 
\begin{equation}
\log (C) = \sum_{n=0}^\infty \frac{(-1)^n}{n+1} (C-Id)^{n+1} 
\end{equation}
terminates after a final number of elements. Hence ${\rm log}(C)$ is an element in $\mathfrak{b}$ and therewith $ \mathcal{G}\subseteq \tilde{\mathcal{G}} $.
\end{proof}

\begin{lemma} 
\label{lemma:3}
$(\mathfrak{b}, *)$ is a closed Abelian Lie group 
\end{lemma}
\begin{proof}{\it [\cref{lemma:3}]}
$(\mathfrak{b}, *)$ is a group: $0\in\mathfrak{b}$ is the identity for $*$, the $*$-inverse of $T\in\mathfrak{b}$ is $-T \in\mathfrak{b}$, and the associative of $*$ follows from the simple form that the BCH formula takes on Abelian algebras. As $\mathfrak{b}$ is a vector space and $*$ corresponds to the addition, the group law and the inverse function mapping are continuous, hence, the group is furthermore closed w.r.t. the euclidean topology. 
\end{proof}

\begin{lemma}
\label{Lemma:group}
$(\mathcal{G}, \odot)$ is a connected Abelian Lie group.
\end{lemma}
\begin{proof} {\it [\cref{Lemma:group}]}
The Abelian group property follows directly from Lemma~\ref{lemma:3}. 
Let $X\in\mathcal{G}$, then the Jordan-Chevalley decomposition is $X = Id + S$ where $S\in \mathfrak{b}$, and $X$ is path connected to the identity via the linear segment $\gamma(t) =  Id + t S$ where $0\leq t\leq 1$. Hence, $(\mathcal{G}, \odot)$ is connected. 
\end{proof}

\begin{theorem}
\label{Th:ExpMapSurj}
Given the Lie group $\mathcal{G}$ with Lie algebra $\mathfrak{b}$. The exponential map $\exp :\mathfrak{b} \to \mathcal{G}$ is surjective. 
\end{theorem}

\begin{proof}{\it [\cref{Th:ExpMapSurj}]}
Let $X \in \mathcal{G}$. The inverse of the exponential series is given by 
\begin{equation}
\log (X) = \sum_{n=0}^\infty \frac{(-1)^n}{n+1} (X-Id)^{n+1} .
\end{equation}
This series terminates (see proof of Lemma~\ref{lemma:5}), and $\log (X) = P(C)\in \mathfrak{b}$, which concludes the proof.
\end{proof}

Note that this theorem can be generalized to the any nilpotent Lie algebra. However, the proof shows (see~\cref{sec:appendix}) that the inverse of the exponential is in this particular case well-defined, which proves the following theorem. 
\begin{theorem}
\label{th:ExpMapInv}
The exponential map from $\mathfrak{b}$ to $\mathcal{G}$ is bijective. 
\end{theorem}

\begin{corollary}
\label{Corr:ClosedLieGroup}
The group $(\mathcal{G}, \odot)$ is closed Lie subgroup of ${\rm GL}_{\mathcal{K}}(\mathbb{R})$ where $\mathcal{K}$ is the dimension of the FCI space, i.e.,
\begin{equation*}
\mathcal{K}= {{K}\choose{N}}.
\end{equation*}
\end{corollary}
\begin{proof}{\it [\cref{Corr:ClosedLieGroup}]}
Since $(\mathfrak{b}, *)$ is closed w.r.t.~the euclidean topology. Together with Theorem~\ref{th:ExpMapInv}, this shows the claim. 
\end{proof}

Similar to the untruncated case, we define the set 
\begin{equation}
\bar{\mathcal{G}} = \{ \exp (T)~|~T\in\bar{\mathfrak{b}}\},
\end{equation}
and find that $(\bar{\mathcal{G}}, \odot)$ is a Lie group. 

\begin{definition}
We define the inclusion map $\iota :  (\bar{\mathcal{G}}, \odot) \to  (\mathcal{G}, \odot)$ with 
\begin{equation}
\iota: \exp(T) \mapsto \exp(\bar{T})
\end{equation}
where 
\begin{equation}
\bar{T} = T + 0,
\end{equation}
i.e., trivially extending $T$ with zeros.
\end{definition}

For the inclusion map $\iota$ we can show the following two lemmata.
 
\begin{lemma}
\label{lemma:iotaHomo}  
The inclusion map $\iota: (\bar{\mathcal{G}}, \odot) \to  (\mathcal{G}, \odot)$ is a group homomorphism.
\end{lemma}
\begin{proof}{\it [\cref{lemma:iotaHomo} ]}
Let $\exp(S),\exp(T)\in\bar{\mathcal{G}} $. Then
\begin{equation}
\begin{aligned}
\iota(\exp(S) \odot \exp(T)) 
&= 
\iota(\exp(S + T))
=
\exp(S + T)
=
\exp(S)) \odot \iota(\exp(T)\\
&=
\iota(\exp(S)) \odot \iota(\exp(T)).
\end{aligned}
\end{equation}
\end{proof}

\begin{lemma}
\label{lemma:iotaImmer}
The inclusion map $\iota: (\bar{\mathcal{G}}, \odot) \to  (\mathcal{G}, \odot)$ is an injective immersion.
\end{lemma}
\begin{proof}{\it [\cref{lemma:iotaImmer}]}
As $\iota$ is taken to be the identity on the cluster operator, we inherit the injectivity from the exponential map. However, as the derivative of the exponential map is again the exponential map, this also implies that the $\iota$ is an immersion.  
\end{proof}

With these properties at hand, we can show the closed connected Abelian subgroup property of $(\bar{\mathcal{G}}, \odot)$. 

\begin{theorem}
\label{Th:Subgroup}
$(\bar{\mathcal{G}}, \odot)$ is a closed connected Abelian subgroup of $(\mathcal{G}, \odot)$. 
\end{theorem}
\begin{proof}{\it [\cref{Th:Subgroup}]}
The theorem follows from Lemma \ref{lemma:iotaHomo} and \ref{lemma:iotaImmer}, i.e., $\iota$ is an injective immersion that is a group homomorphism, which shows the subgroup property. This subgroup is furthermore Abelian. Moreover, analogously to the proof of Lemma~\ref{Lemma:group}, each element $X\in \bar{\mathcal{G}}$ is path connected to the identity. The subgroup is indeed closed w.r.t.~the euclidean topology, as any Cauchy sequence in $\bar{\mathcal{G}}$ clearly converges in $\bar{\mathcal{G}}$. 
\end{proof}

\begin{corollary}
\label{Corr:EmbeddedLieSubGroup}
$\bar{\mathcal{G}}$ is an embedded Lie subgroup.
\end{corollary}

\begin{proof}{\it [\cref{Corr:EmbeddedLieSubGroup}]}
This is a direct application of Cartan's theorem~\cite{lee2013smooth}.
\end{proof}

\subsection{Coupled Cluster Equation}
\label{sec:ccWorkingEQ}

We shall begin this section with the complete derivation of the coupled cluster working equations.
Similarly to the example in the main text, we get by using above rules
\begin{align*}
\langle \Phi_i^a, [[H,T],T]\Phi_0\rangle 
&= 
\langle\Phi_i^a, HT^2\Phi_0\rangle
-2 \langle \Phi_i^a,THT\Phi_0\rangle\\
&=
\sum_{\substack{j\neq k \\ b\neq c}} 
t_j^bt_k^c \langle \Phi_i^a,H\Phi_{jk}^{bc}\rangle
+
\sum_{\substack{j\neq k \\ b\neq c}} 
t_j^bt_{ik}^{ac} \langle\Phi_i^a,H\Phi_{ijk}^{abc}\rangle
\\
&\quad
-2t_i^a\sum_{\substack{j\neq k \\ b\neq c}} 
t_{jk}^{bc} \langle\Phi_0,H\Phi_{jk}^{bc}\rangle
-
2t_i^a\sum_{j,b}t_j^b\langle\Phi_0,H\Phi_j^b\rangle
\end{align*}

and
\begin{align*}
\langle \Phi_i^a, [[[H,T],T],T]\Phi_0\rangle 
&= 
\langle\Phi_i^a, HT^3\Phi_0\rangle
- 
3\langle \Phi_i^a,THT^2\Phi_0\rangle\\
&=
\sum_{\substack{j\neq k \neq l \\ b\neq c\neq d }} t_j^bt_k^ct_l^d \langle\Phi_i^a,H\Phi_{jkl}^{bcd}\rangle\\
&\quad-
3t_i^a
\left(
\sum_{\substack{j\neq k\\\ b\neq c}}
t_j^bt_k^c\langle \Phi_0, H\Phi_{jk}^{cd}\rangle
+
\sum_{\substack{j\neq k \\ b\neq c }} 
t_j^bt_{ik}^{ac} \langle\Phi_i^a,H\Phi_{ijk}^{abc}\rangle
+
\sum_{\substack{ j\neq k \\ b\neq c }} 
t_i^at_{jk}^{bc} \langle\Phi_i^a,H\Phi_{jk}^{bc}\rangle
\right)
\end{align*}

Due to the Slater--Condon rules, there exists no quadruple commutator contribution in the single-projection equations.
Hence, the resulting polynomials that arise from projecting onto singly excited Slater determinants are of the form
\begin{equation}
f_i^a(t)
=
\langle
\Phi_i^a, e^{-T} H e^T \Phi_0 
\rangle
\end{equation}

Next we shall derive the polynomials that arise from projecting onto doubly excited Slater determinants. 
For the regular commutator term, we find
\begin{align*}
\langle \Phi_{ij}^{ab}, [H,T]\Phi_0\rangle  
&=  
\langle \Phi_{ij}^{ab}, HT\Phi_0\rangle
- 
\langle \Phi_{ij}^{ab}, TH\Phi_0\rangle\\
&=
\sum_{k,c} t_k^c 
\langle \Phi_{ij}^{ab},H\Phi_k^c\rangle
+
\sum_{\substack{k<l\\c<d}}
t_{kl}^{cd}\langle\Phi_{ij}^{ab},H\Phi_{kl}^{cd}\rangle\\
&\quad -t_{ij}^{ab}\langle\Phi_0, H\Phi_0\rangle
-
\sum_{
\substack{
k\in\{i,j\},\ l\in \{i,j\}\setminus k \\
c\in\{a,b\},\ d\in \{a,b\}\setminus c}
} t_k^c\langle\Phi_l^d,H\Phi_0\rangle
\end{align*}
Note that the canonical ordering, which is imposed here, results in further restrictions to the possible index pairs when de-exciting with a single excitation operator, see e.g. last the term in the equation above. 
For the double commutator we obtain
\begin{align*}
\langle \Phi_{ij}^{ab}, [[H,T],T]\Phi_0\rangle  
&=  
\langle \Phi_{ij}^{ab}, HT^2\Phi_0\rangle
- 2
\langle \Phi_{ij}^{ab}, THT\Phi_0\rangle
+
\langle\Phi_{ij}^{ab}T^2H\Phi_0\rangle\\
&=
\sum_{\substack{k\neq l\\ c\neq d}}
t_{k}^ct_{l}^{d}\langle\Phi_{ij}^{ab},H\Phi_{kl}^{cd}\rangle
+
\sum_{\substack{k\neq l<m \\ c\neq d<e}}
t_{k}^{c}t_{lm}^{de}\langle\Phi_{ij}^{ab},H\Phi_{klm}^{cde}\rangle\\
&\quad 
+ \sum_{\substack{k\neq l \neq i,j\\ c\neq d \neq a,b}}
t_{ij}^{ab}t_{kl}^{ce}\langle
\Phi_{ij}^{ab},
H\Phi_{ijkl}^{abcd}
\rangle
-2
t_{ij}^{ab}
\sum_{k,c}
t_k^c\langle\Phi_0, H\Phi_k^c\rangle \\
&\quad
-2
\sum_{
\substack{
k\in\{i,j\},\ l\in \{i,j\}\setminus k \\
c\in\{a,b\},\ d\in \{a,b\}\setminus c}
} 
\sum_{m, e}
t_k^c
t_m^e
\langle\Phi_l^d,H\Phi_m^e\rangle\\
&\quad +
\sum_{
\substack{
k\in\{i,j\},\ l\in \{i,j\}\setminus k \\
c\in\{a,b\},\ d\in \{a,b\}\setminus c}
} 
t_k^c
t_l^d
\langle\Phi_0,H\Phi_0\rangle
\end{align*}
In a similar spirit, we find for the triple commutator term
\begin{align*}
\langle \Phi_{ij}^{ab}, [[[H,T],T],T]\Phi_0\rangle  
&=
\langle \Phi_{ij}^{ab}, HT^3\Phi_0\rangle 
-3 \langle \Phi_{ij}^{ab}, THT^2\Phi_0\rangle
+3 \langle \Phi_{ij}^{ab}, T^2HT\Phi_0\rangle\\ 
&=
\sum_{\substack{k\neq l \neq m \\ c\neq d \neq e}}
t_k^c t_l^d t_m^e
\langle \Phi_{ij}^{ab}, H\Phi_{klm}^{cde}\rangle 
+
\sum_{\substack{k\neq l \neq m \neq n\\ c\neq d \neq e \neq f}}
t_k^ct_l^dt_{m,n}^{ef} 
\langle \Phi_{ij}^{ab}, H\Phi_{klmn}^{cdef}\rangle
\\
&\quad 
-3
\sum_{
\substack{
k\in\{i,j\},\ l\in \{i,j\}\setminus k \\
c\in\{a,b\},\ d\in \{a,b\}\setminus c}
} 
\sum_{\substack{m\neq n \\ e \neq f}}
t_k^ct_m^et_n^f\langle\Phi_l^d,H\Phi_{mn}^{ef}\rangle\\
&\quad
-3
\sum_{
\substack{
k\in\{i,j\},\ l\in \{i,j\}\setminus k \\
c\in\{a,b\},\ d\in \{a,b\}\setminus c}
} 
\sum_{\substack{m\neq n<o \\ e\neq f<g}}
t_k^c
t_m^et_{no}^{fg}
\langle\Phi_l^d,H\Phi_{mno}^{efg}\rangle\\
&\quad +3
\sum_{
\substack{
k\in\{i,j\},\ l\in \{i,j\}\setminus k \\
c\in\{a,b\},\ d\in \{a,b\}\setminus c}
} 
\sum_{m, e}
t_k^c
t_l^d
t_m^e
\langle\Phi_0,H\Phi_m^e\rangle
\end{align*}

For the quadruple commutator, we obtain

\begin{align*}
\langle \Phi_{ij}^{ab}, [[[[H,T],T],T],T]\Phi_0\rangle  
&=
\langle \Phi_{ij}^{ab}, HT^4\Phi_0\rangle
-4
\langle \Phi_{ij}^{ab}, THT^3\Phi_0\rangle
+6
\langle \Phi_{ij}^{ab},T^2HT^2\Phi_0\rangle\\
&=
\sum_{\substack{k\neq l \neq m \neq n\\ c\neq d \neq e \neq f}}
t_k^c t_l^d t_m^e t_n^f
\langle \Phi_{ij}^{ab}, H \Phi_{klmn}^{cdef} \rangle\\
&\quad-4
\sum_{
\substack{
k\in\{i,j\},\ l\in \{i,j\}\setminus k \\
c\in\{a,b\},\ d\in \{a,b\}\setminus c}
} 
\sum_{\substack{m\neq n\neq o \\ e \neq f \neq g}}
t_k^c t_m^e t_n^f t_o^g \langle\Phi_l^d,H\Phi_{mnp}^{efg}\rangle\\
&\quad+6
\sum_{
\substack{
k\in\{i,j\},\ l\in \{i,j\}\setminus k \\
c\in\{a,b\},\ d\in \{a,b\}\setminus c}
} 
\sum_{
\substack{
m\neq n\\ e\neq f
}
}
t_k^c
t_l^d
t_m^e
t_n^f
\langle\Phi_0,H\Phi_{mn}^{ef}\rangle
\end{align*}

\begin{proof}{\it [\cref{lem:S_5_not_extr}]}
For any $v\in \mathcal{S}_5$ with non-zero entries in $i,j,k$ there exists two vectors $v_1, v_2 \in \mathcal{S}_7$ where $(v_1)_i = (v_2)_k = 2$ and  $(v_1)_j = (v_2)_j = 1$, then, $v=0.5 v_1 + 0.5 v_2$. Therefore $\mathcal{S}_5$ are not extremal points of the maximal polytope.
\end{proof}

\begin{proof}{\it [\cref{lem:extrm_vert_sing}]}
We show that $\nexists v_i \in \mathcal{S}_g$ s.t. $\exists (\alpha_j)$ with 
\[
v_i = \sum_{i\neq j} \alpha_j v_j 
\]
and $1\geq \alpha_j \geq 0$ and $\sum_j \alpha_j = 1$ for $g\in\{2,4,6,7\}$ and $v_j \in \mathcal{S}_g$.
We investigate the individual cases:
\paragraph{\underline{$v_i \in \mathcal{S}_2$}}

Let $v_i \in \mathcal{S}_2$. Elements in $\mathcal{S}_2$ have exactly one entry; since convex combinations cannot remove entries, this set has to be extremal.

\paragraph{\underline{$v_i \in \mathcal{S}_4$}}

Let $v_i \in \mathcal{S}_4$. 
Then, since convex combinations cannot remove entries $\alpha_k$ has to be zero for vectors in $\mathcal{S}_6$. 
Moreover, $\alpha_k$ has to be zero for vertices in $\mathcal{S}_4$; otherwise we would obtain vertices with at least three entries.
Also, $\alpha_k$ has to be zero for vertices in $\mathcal{S}_7$, since $\mathcal{S}_7$ has two single entries.
Hence, the only possible linear combination is 
\[
v_i =\alpha_1 v_1 + \alpha_2 v_2
\]
where $v_1 \in \mathcal{S}_2$ with the same non-zero double entry as $v_i$, and $v_2 \in \mathcal{S}_6$ with the same non-zero single entry as $v_i$. However, then $\alpha_2= 1/2$ and $\alpha_1 =1$ which is not convex.

\paragraph{\underline{$v_i \in \mathcal{S}_6$}}

Let $v_i \in \mathcal{S}_6$. Elements in $\mathcal{S}_6$ have exactly one entry which differs in position from the vectors in $\mathcal{S}_2$. Again, since convex combinations cannot remove entries, this set has to be extremal.

\paragraph{\underline{$v_i \in \mathcal{S}_7$}}
Let $v_i \in \mathcal{S}_7$. Then, $\alpha_k$ has to be zero for vectors in $\mathcal{S}_2$ and $\mathcal{S}_4$ since those vectors have non-zero entries in the double component.
Moreover, $\alpha_k$ has to be zero for vectors in $\mathcal{S}_5$ since vectors in $\mathcal{S}_5$ have three non-zero components.
Hence, the only possible linear combination is 
\[
v_i =\alpha_1 v_1 + \alpha_2 v_2
\]
where $v_1, v_2 \in \mathcal{S}_6$ with the same non-zero single entries as $v_i$. However, then $\alpha_2= 1/2$ and $\alpha_1 =1$ which is not convex.
\end{proof}

\subsection{Proof of~\cref{th:quadQe}}

The main idea is to show, that in each coupled cluster equation every term is a product of $1,t_i^a,t_{ij}^{ab}$ and $t_i^at_j^b-t_i^bt_j^a$. Then we introduce new variables $\tilde t_{i,j,a,b} = t_i^at_j^b-t_i^bt_j^a$ which are used in Eq.~\eqref{eq:double_proj_with_explicit_dummy},~\eqref{eq:double_proj_with_explicit_dummy_2} and~\eqref{eq:double_proj_with_explicit_dummy_3}. There are as many $\tilde t_{i,j,a,b}$ as double excitations. The old variables in canonical ordering together with the new ones will be called $x$.
\subsubsection{Singles projections}
We first look at the equations that arise from single projections, recall
\begin{equation}
\begin{aligned}
f_i^a
&=
\langle
\phi_i^a, H \phi_0
\rangle
+
\langle
\phi_i^a, [H,T] \phi_0
\rangle
+
\frac{1}{2}
\langle
\phi_i^a, [[H,T],T] \phi_0
\rangle
\\
&\quad +
\frac{1}{6}
\langle
\phi_i^a, [[[H,T],T],T] \phi_0
\rangle\\
&=
\langle
\phi_i^a, H \phi_0
\rangle
+
\langle
\phi_i^a, HT \phi_0
\rangle
-
\langle
\phi_i^a, TH \phi_0
\rangle+
\frac{1}{2}
\langle
\phi_i^a, HT^2 \phi_0
\rangle\\
&\quad
-
\langle
\phi_i^a, THT \phi_0
\rangle
+
\frac{1}{6}
\langle
\phi_i^a, HT^3 \phi_0
\rangle
-
\frac{1}{2}
\langle
\phi_i^a, THT^2 \phi_0
\rangle\\
\end{aligned}
\end{equation}

\noindent
We begin with the terms 
\begin{equation}
\langle
\phi_i^a, H \phi_0
\rangle
+
\langle
\phi_i^a, HT \phi_0
\rangle
+
\frac{1}{2}
\langle
\phi_i^a, HT_s^2 \phi_0
\rangle,
\end{equation}
where $\langle
\phi_i^a, HT^2 \phi_0
\rangle = \langle
\phi_i^a, HT_s^2 \phi_0
\rangle$ due to the Slater--Condon rules,  and find
\begin{equation}
\label{eq:Single_projection_part1}
\begin{aligned}
\langle
\phi_i^a, H \phi_0
\rangle
&+
\langle
\phi_i^a, HT \phi_0
\rangle
+
\frac{1}{2}
\langle
\phi_i^a, HT_s^2 \phi_0
\rangle
\\
&=
\langle
\phi_i^a, H \phi_0
\rangle
+
\sum_{j}\sum_{b}
t_j^b
\langle
\phi_i^a, H \phi_j^b 
\rangle
+
\sum_{j<k}\sum_{b<c}
t_{jk}^{bc}
\langle
\phi_i^a, H \phi_{jk}^{bc} 
\rangle\\
&\quad
+
\sum_{j\neq k }\sum_{b\neq c }
t_j^b t_k^c
\langle
\phi_i^a, H \phi_{jk}^{bc}
\rangle
\\
&=
\langle
\phi_i^a, H \phi_0
\rangle
+
\sum_{j}\sum_{b}
t_j^b
\langle
\phi_i^a, H \phi_j^b 
\rangle
+
\sum_{j<k}\sum_{b<c}
t_{jk}^{bc}
\langle
\phi_i^a, H \phi_{jk}^{bc} 
\rangle\\
&\quad
+
\sum_{j< k}\sum_{b< c }
(t_j^b t_k^c - t_j^c t_k^b)
\langle
\phi_i^a, H \phi_{jk}^{bc}
\rangle\\
\end{aligned}
\end{equation}

\noindent
Next, we take a look at
\begin{equation}
-
\langle
\phi_i^a, TH \phi_0
\rangle
-
\langle
\phi_i^a, THT \phi_0
\rangle
-
\frac{1}{2}
\langle
\phi_i^a, THT_s^2 \phi_0
\rangle,
\end{equation}

\noindent
where $\langle
\phi_i^a, THT^2 \phi_0
\rangle = 
\langle
\phi_i^a, THT_s^2 \phi_0
\rangle$ due to the Slater--Condon rules, and find

\begin{equation}
\label{eq:Single_projection_part2}
\begin{aligned}
-
&\langle
\phi_i^a, TH \phi_0
\rangle
-
\langle
\phi_i^a, THT \phi_0
\rangle
-
\frac{1}{2}
\langle
\phi_i^a, THT_s^2 \phi_0
\rangle\\
&=
-
t_i^a
\langle
\phi_0, H \phi_0
\rangle
-
t_i^a
\sum_{j}\sum_{b}
t_j^b
\langle
\phi_0, H \phi_j^b
\rangle\\
&\quad-
t_i^a
\sum_{j<k}\sum_{b<c}
t_{jk}^{bc}
\langle
\phi_0, H \phi_{jk}^{bc} 
\rangle
-
t_i^a
\sum_{j\neq k}\sum_{b\neq c}
t_j^bt_k^c
\langle
\phi_0, H \phi_{jk}^{bc}
\rangle\\
&=
-
t_i^a
\langle
\phi_0, H \phi_0
\rangle
-
t_i^a
\sum_{j}\sum_{b}
t_j^b
\langle
\phi_0, H \phi_j^b
\rangle\\
&\quad-
t_i^a
\sum_{j<k}\sum_{b<c}
t_{jk}^{bc}
\langle
\phi_0, H \phi_{jk}^{bc} 
\rangle
-
t_i^a
\sum_{j< k }\sum_{b< c}
(t_j^b t_k^c - t_j^c t_k^b)
\langle
\phi_0, H \phi_{jk}^{bc}
\rangle
\end{aligned}
\end{equation}

\noindent
Introducing a matrix notation we may write Eqs~\eqref{eq:Single_projection_part1} and~\eqref{eq:Single_projection_part2} 
\begin{equation}
\begin{aligned}
&\phi_i^a, H \phi_0
\rangle
+
\langle
\phi_i^a, HT \phi_0
\rangle
+
\frac{1}{2}
\langle
\phi_i^a, HT_s^2 \phi_0
\rangle\\
&\quad 
-
\langle
\phi_i^a, TH \phi_0
\rangle
-
\langle
\phi_i^a, THT \phi_0
\rangle
-
\frac{1}{2}
\langle
\phi_i^a, THT_s^2 \phi_0
\rangle\\
&=
\begin{bmatrix}
1 & 
\dots  t_i^a \dots  
& {\bf t_d} & {\bf t_{x}}
\end{bmatrix}
\begin{bmatrix} 
H_{(i,a),0} & H_{(i,a),s} & H_{(i,a),d} & H_{(i,a),d}\\
\vdots&\vdots&\vdots&\vdots\\
-H_{0,0} & -H_{0,s} & -H_{0,d} & -H_{0,d}\\
\vdots&\vdots&\vdots&\vdots\\ 
\end{bmatrix}
\begin{bmatrix}
1\\
{\bf t_s}\\ 
{\bf  t_d}\\
{\bf t_{x}}
\end{bmatrix}
\end{aligned}
\end{equation}
where the ``\dots''-entries in the matrix correspond to zero entries.
Next, we take a look at the terms
\begin{equation}
\frac{1}{6}
\langle
\phi_i^a, HT^3 \phi_0
\rangle
+
\frac{1}{2}
\langle
\phi_i^a, HT_sT_d \phi_0
\rangle,
\end{equation}
where $\langle
\phi_i^a, HT^2 \phi_0
\rangle = \langle
\phi_i^a, HT_sT_d \phi_0
\rangle$
due to the Slater--Condon rules, and find 
\begin{equation}
\begin{aligned}
\frac{1}{6}
&
\langle
\phi_i^a, HT^3 \phi_0
\rangle
+
\frac{1}{2}
\langle
\phi_i^a, HT_sT_d \phi_0
\rangle\\
&=
\frac{1}{6}
\sum_{j\neq k\neq l }\sum_{b\neq c\neq d }
t_j^bt_k^ct_l^d
\langle
\phi_i^a, H \phi_{jkl}^{bcd}
\rangle
+
\frac{1}{2}
\sum_{j\neq k\neq l }\sum_{b\neq c\neq d }
t_j^b t_{kl}^{cd} 
\langle
\phi_i^a, H\phi_{jkl}^{bcd}
\rangle\\
&=
\sum_{j\neq k< l }\sum_{b\neq c< d }
t_j^b(t_k^ct_l^d-t_k^dt_l^c)
\langle
\phi_i^a, H \phi_{c.o.(jkl)}^{c.o.(bcd)}
\rangle\\
&\quad 
+
\sum_{\substack{j\neq k\neq l \\ k<l}} \sum_{\substack{b\neq c\neq d \\ c<d}}
t_j^b t_{kl}^{cd} 
\langle
\phi_i^a, H\phi_{c.o.(jkl)}^{c.o.(bcd)}
\rangle\\
&=
\begin{bmatrix}
1 & {\bf t_s} & {\bf t_d} & {\bf t_{x}}
\end{bmatrix}
\begin{bmatrix} 0&0&0&0\\
0 & 0 & {\bf H}_{(i,a),tr} & {\bf H}_{(i,a),tr}\\
0 & 0 & 0 & 0\\  0&0&0&0
\end{bmatrix}
\begin{bmatrix}
1\\
{\bf t_s}\\ 

{\bf  t_d}\\
{\bf t_{x}}
\end{bmatrix}
\end{aligned}
\end{equation}

\noindent
We here introduced the sub-matrix ${\bf H}_{(i,a),tr}$ which also has a sparsity pattern due to the Slater--Condon rules.
Elementwise we define 
\begin{equation}
\begin{footnotesize}
({\bf H_{(i,a),tr}})_{(j,b),(kl,cd)} = \begin{cases}
0& \text{ if } (i,a)\notin\ \{(jkl,bcd)\}\\
\mathrm{sign}(\sigma)H_{(i,a),C.O.((ikl,acd))}& \text{ if } (i,a)=(j,b)\\
\mathrm{sign}(\sigma)H_{(i,a),C.O.((ikl,bad))}& \text{ if } (i,a)=(j,c)\\
\mathrm{sign}(\sigma)H_{(i,a),C.O.((ikl,bca))}& \text{ if } (i,a)=(j,d)\\
\mathrm{sign}(\sigma)H_{(i,a),C.O.((jil,acd))}& \text{ if } (i,a)=(k,b)\\
\mathrm{sign}(\sigma)H_{(i,a),C.O.((jil,bad))}& \text{ if } (i,a)=(k,c)\\
\mathrm{sign}(\sigma)H_{(i,a),C.O.((jil,bca))}& \text{ if } (i,a)=(k,d)\\
\mathrm{sign}(\sigma)H_{(i,a), C.O.(\{(jki,acd)\})}& \text{ if } (i,a)=(l,b)\\
\mathrm{sign}(\sigma)H_{(i,a), C.O.(\{(jki,bad)\})}& \text{ if } (i,a)=(l,c)\\
\mathrm{sign}(\sigma)H_{(i,a), C.O.(\{(jki,bca)\})}& \text{ if } (i,a)=(l,d),
\end{cases}
\end{footnotesize}
\end{equation}

\noindent
where $\sigma $ is the permutation to achieve canonical ordering and $C.O.$ is the canonical ordering of the indices. The only dense row in this sub-matrix appears if $(i,a)=(j,b)$. Otherwise, we only have as many non-zero entries per row as single amplitudes. 
\\

\noindent
In summary, we conclude that 
\begin{equation}
\begin{footnotesize}
\begin{aligned}
f_i^a
&=
\begin{bmatrix}
1 & 
\dots  t_i^a \dots  
& {\bf t_d} & {\bf t_{x}}
\end{bmatrix}
\begin{bmatrix} 
H_{(i,a),0} & H_{(i,a),s} & H_{(i,a),d} & H_{(i,a),d}\\
\vdots&\vdots&\vdots&\vdots\\
-H_{0,0} & -H_{0,s} & -H_{0,d} & -H_{0,d}\\
\vdots&\vdots&\vdots&\vdots\\ 
\end{bmatrix}
\begin{bmatrix}
1\\
{\bf t_s}\\ 
{\bf  t_d}\\
{\bf t_{x}}
\end{bmatrix}\\
&\quad 
+
\begin{bmatrix}
1 & {\bf t_s} & {\bf t_d} & {\bf t_{x}}
\end{bmatrix}
\begin{bmatrix} 0&0&0&0\\
0 & 0 & {\bf H}_{(i,a),tr} & {\bf H}_{(i,a),tr}\\
0 & 0 & 0 & 0\\  0&0&0&0
\end{bmatrix}
\begin{bmatrix}
1\\
{\bf t_s}\\ 

{\bf  t_d}\\
{\bf t_{x}}
\end{bmatrix}
\end{aligned}
\end{footnotesize}
\end{equation}
which summarizes the sparsity pattern of the CCSD single's equations. 

\subsubsection{doubles projections}
Next, we shall derive the polynomials that arise from projecting onto doubly excited Slater determinants. 
Recall that 
\begin{equation}
\begin{aligned}
f_{ij}^{ab}
&=
\langle
\phi_{ij}^{ab}, H \phi_0
\rangle
+
\langle \phi_{ij}^{ab}, [H,T]\phi_0\rangle  
+
\frac{1}{2}
\langle \phi_{ij}^{ab}, [[H,T],T]\phi_0\rangle \\
&\quad
+
\frac{1}{6}
\langle \phi_{ij}^{ab}, [[[H,T],T],T]\phi_0\rangle  
+
\frac{1}{24}
\langle \phi_{ij}^{ab}, [[[[H,T],T],T],T]\phi_0\rangle\\
&=
\langle
\phi_{ij}^{ab}, H \phi_0
\rangle
+
\langle \phi_{ij}^{ab}, HT\phi_0\rangle
- 
\langle \phi_{ij}^{ab}, TH\phi_0\rangle
+
\frac{1}{2}
\langle \phi_{ij}^{ab}, HT^2\phi_0\rangle\\
&\quad
- 
\langle \phi_{ij}^{ab}, THT\phi_0\rangle
+
\frac{1}{2}
\langle\phi_{ij}^{ab},T^2H\phi_0\rangle
+
\frac{1}{6}
\langle \phi_{ij}^{ab}, HT^3\phi_0\rangle \\
&\quad
-\frac{1}{2}
\langle \phi_{ij}^{ab}, THT^2\phi_0\rangle
+\frac{1}{2}
\langle \phi_{ij}^{ab}, T^2HT\phi_0\rangle
+
\frac{1}{24}
\langle \phi_{ij}^{ab}, HT^4\phi_0\rangle\\
&\quad 
-\frac{1}{6}
\langle \phi_{ij}^{ab}, THT^3\phi_0\rangle
+\frac{1}{4}
\langle \phi_{ij}^{ab},T^2HT^2\phi_0\rangle
\end{aligned}
\end{equation}

\noindent
We begin with the terms 
\begin{equation}
\langle
\phi_{ij}^{ab}, H \phi_0
\rangle
+
\langle
\phi_{ij}^{ab}, HT \phi_0
\rangle
+
\frac{1}{2}
\langle
\phi_{ij}^{ab}, HT_s^2 \phi_0
\rangle,
\end{equation}

\noindent
where we separated the term 
\begin{equation}
\langle
\phi_{ij}^{ab}, HT^2 \phi_0
\rangle
=
\langle
\phi_{ij}^{ab}, HT_s^2 \phi_0
\rangle
+\langle
\phi_{ij}^{ab}, HT_sT_d \phi_0
\rangle
+\langle
\phi_{ij}^{ab}, HT_d^2 \phi_0
\rangle
\end{equation}
and included only the first term, the remaining terms are discussed in Eqs.~\eqref{eq:T2_splie_part_2} and~\eqref{eq:T2_splie_part_3}.
We find
\begin{equation}
\label{eq:doubles_first_part}
\begin{aligned}  
&\langle
\phi_{ij}^{ab}, H \phi_0
\rangle
+
\sum_{k,c} t_k^c 
\langle \phi_{ij}^{ab},H\phi_k^c\rangle
+
\sum_{\substack{k<l\\c<d}}
t_{kl}^{cd}\langle\phi_{ij}^{ab},H\phi_{kl}^{cd}\rangle\\
&\quad
+
\frac 1 2\sum_{\substack{k\neq l\\ c\neq d}}
t_{k}^ct_{l}^{d}\langle\phi_{ij}^{ab},H\phi_{kl}^{cd}\rangle\\
&=
\langle
\phi_{ij}^{ab}, H \phi_0
\rangle
+
\sum_{k,c} t_k^c  \langle \phi_{ij}^{ab},H\phi_k^c\rangle
+
\sum_{\substack{k<l\\c<d}}
t_{kl}^{cd}\langle\phi_{ij}^{ab},H\phi_{kl}^{cd}\rangle\\
&\quad
+\sum_{\substack{k< l\\ c< d}}
(t_{k}^ct_{l}^{d} -t_{k}^dt_{l}^{c})\langle\phi_{ij}^{ab},H\phi_{kl}^{cd}\rangle\\
=&\begin{bmatrix}
1 & {\bf t_s} & {\bf t_d} & {\bf t_{x}}
\end{bmatrix}
\begin{bmatrix}
H_{(ij,ab),0} & H_{(ij,ab),s} & H_{(ij,ab),d} & H_{(ij,ab),d}\\
0 & 0 & 0 & 0\\ \vdots&\vdots&\vdots&\vdots\\ 0&0&0&0
\end{bmatrix}
\begin{bmatrix}
1\\
{\bf t_s}\\ 
{\bf t_d}\\
{\bf t_{x}}
\end{bmatrix}
\end{aligned}
\end{equation}

\noindent
Next, we take a look at
\begin{equation}
\label{eq:doubles_second_part}
\begin{aligned}
& -
\langle
\phi_{ij}^{ab}, TH \phi_0
\rangle
-
\langle
\phi_{ij}^{ab}, THT \phi_0
\rangle
-
\frac{1}{2}
\langle
\phi_{ij}^{ab}, THT_s^2 \phi_0
\rangle\\
&=
F_{ij,ab}^d + F_{i,a}^s + F_{j,b}^s - F_{i,b}^s - F_{j,a}^s
\end{aligned}
\end{equation}
where the expressions $F_{ij,ab}^d$ and $F_{i,a}^s$ are given by 
\begin{equation}
\begin{aligned}
F_{ij,ab}^d &= 
-
t_{ij}^{ab}\langle \phi_0,H\phi_0\rangle 
-
t_{ij}^{ab} \sum_{(k,c)}t_k^c \langle \phi_0, H\phi_{k}^{c}\rangle \\
&\quad-
t_{ij}^{ab}\sum_{\substack{k<l\\c<d}} t_{kl}^{cd}\langle\phi_0,H\phi_{kl}^{cd}\rangle
-
t_{ij}^{ab}\sum_{\substack{k< l\\ c< d}}
(t_{k}^ct_{l}^{d}-t_{k}^dt_{l}^{c})\langle\phi_0,H\phi_{kl}^{cd}\rangle
\end{aligned}
\end{equation}

\noindent
and
\begin{equation}
\begin{aligned}
F_{i,a}^s &=
- t_i^a \langle \phi_j^b,H\phi_0\rangle 
- 
t_i^a \sum_{(k,c)}t_k^c \langle \phi_j^b, H\phi_{k}^{c}\rangle \\
&\quad
- 
t_i^a \sum_{\substack{k<l\\c<d}}
t_{kl}^{cd}\langle\phi_j^b,H\phi_{kl}^{cd}\rangle
- 
t_i^a \sum_{\substack{k< l\\ c< d}}
(t_{k}^ct_{l}^{d}-t_{k}^dt_{l}^{c}) \langle\phi_j^b,H\phi_{kl}^{cd}\rangle
\end{aligned}
\end{equation}

\noindent
We may combine Eqs.~\eqref{eq:doubles_first_part} and~\eqref{eq:doubles_second_part} and find

\begin{equation}
\begin{aligned}
\langle
\phi_{ij}^{ab}, H \phi_0
\rangle
&+
\langle
\phi_{ij}^{ab}, HT \phi_0
\rangle
+
\frac{1}{2}
\langle
\phi_{ij}^{ab}, HT_s^2 \phi_0
\rangle\\
&\quad
-
\langle
\phi_{ij}^{ab}, TH \phi_0
\rangle
-
\langle
\phi_{ij}^{ab}, THT \phi_0
\rangle
-
\frac{1}{2}
\langle
\phi_{ij}^{ab}, THT_s^2 \phi_0
\rangle\\
=&
\begin{bmatrix}
1 \\ 
\vdots\\
t_i^a\\
\vdots \\
t_i^b \\
\vdots \\ 
t_j^a \\
\vdots \\ 
t_j^b \\
\vdots\\ 
\vdots \\t_{ij}^{ab} \\
\vdots \\
{\bf t_{x}}
\end{bmatrix}^\top
\begin{bmatrix}
H_{(ij,ab),0} & H_{(ij,ab),s} & H_{(ij,ab),d} & H_{(ij,ab),d}\\
\vdots&\vdots&\vdots&\vdots\\
H_{(j,b),0} & H_{(j,b),s} & H_{(j,b),d} & H_{(j,b),d}\\
\vdots&\vdots&\vdots&\vdots\\
H_{(j,a),0} & H_{(j,a),s} & H_{(j,a),d} & H_{(j,a),d}\\
\vdots&\vdots&\vdots&\vdots\\
H_{(i,b),0} & H_{(i,b),s} & H_{(i,b),d} & H_{(i,b),d}\\
\vdots&\vdots&\vdots&\vdots\\
H_{(i,a),0} & H_{(i,a),s} & H_{(i,a),d} & H_{(i,a),d}\\
\vdots&\vdots&\vdots&\vdots\\
H_{0,0} & H_{0,s} & H_{0,d} & H_{0,d}\\
\vdots&\vdots&\vdots&\vdots\\ 
0&0&0&0
\end{bmatrix}
\begin{bmatrix}
1\\
{\bf t_s}\\ 
{\bf t_d}\\
{\bf t_{x}}
\end{bmatrix}
\end{aligned}
\end{equation}
Note that the ``$\cdots$'' in the matrix are zero entries, and we only list the non-zeros lines. We moreover assumed an ordering where $t_i^a <t_i^b < t_j^a <t_j^b$.

Next, we take a look at
\begin{equation}
\label{eq:T2_splie_part_2}
\begin{aligned}
\frac{1}{2}
\langle
\phi_{ij}^{ab}, H T_sT_d \phi_0
\rangle
&=
\sum_{k,c}\sum_{\substack{\ell < m \\ d < e}}
t_{k}^ct_{\ell m}^{de}
\langle
\phi_{ij}^{ab}, H \phi_{c.o.(k\ell m)}^{c.o.(cde)}\rangle
\\
&=\begin{bmatrix}
1 & {\bf t_s} & {\bf t_d} & {\bf t_{x}}
\end{bmatrix}
\begin{bmatrix} 0&0&0&0\\
0 & 0 & {\bf H}_{(ij,ab),tr} & 0\\
0 & 0 & 0 & 0\\  0&0&0&0
\end{bmatrix}
\begin{bmatrix}
1\\
{\bf t_s}\\ 
{\bf t_d}\\
{\bf t_{x}}
\end{bmatrix}
\end{aligned}
\end{equation}

\noindent 
where we define
\begin{equation}
\begin{aligned}
&({\bf H_{(ij,ab),tr}})_{(k,c),(\ell m,d e)} \\
&\qquad\qquad= \begin{cases}
\mathrm{sign}(\sigma)H_{(ij,ab),C.O.(k\ell m,cde)}& \text{if } |\{i,j\}\cap \{k,l,m\}|\geq 1 \\&\& |\{a,b\}\cap\{c,d,e\}|\geq 1\\
0& \text{else}
\end{cases}
\end{aligned}
\end{equation}

\noindent
Next, we look at the term
\begin{equation}
\label{eq:T2_splie_part_3}
\begin{aligned}
\frac{1}{2}
\langle
\phi_{ij}^{ab}, H T_dT_d \phi_0
\rangle
&= 
\sum_{\substack{k < \ell \neq m < n \\ c< d\neq e < f}} t_{k\ell}^{cd}t_{mn}^{ef}\langle \phi_{ij}^{ab},H \phi_{k\ell m n}^{cdef}\rangle\\
&=\begin{bmatrix}
1 & {\bf t_s} & {\bf t_d} & {\bf t_{x}}
\end{bmatrix}
\begin{bmatrix} 0&0&0&0\\
0 & 0 & 0 & 0\\
0 & 0 & {\bf H}_{(ij,ab),qd}& 0\\  0&0&0&0
\end{bmatrix}
\begin{bmatrix}
1\\
{\bf t_s}\\ 
{\bf t_d}\\
{\bf t_{x}}
\end{bmatrix}
\end{aligned}
\end{equation}
where we define
\begin{equation}
\begin{aligned}
&({\bf H_{(ij,ab),qd}})_{(k\ell,cd),(m n, e f)}\\
&\qquad \qquad
= \begin{cases}
\mathrm{sign}(\sigma)H_{(ij,ab),C.O.((k\ell m n,cdef))}& \text{if } \{i,j\}\in \{k,l,m,n\} \\&\& \{a,b\}\in \{c,d,e,f\}\\
0& \text{else}
\end{cases}
\end{aligned}
\end{equation}

\noindent
Next, we investigate the term $\langle \phi_{ij}^{ab}, HT^3\phi_0\rangle$, where we find 
\begin{equation}
\langle \phi_{ij}^{ab}, HT^3\phi_0\rangle = \langle \phi_{ij}^{ab}, HT_s^3\phi_0\rangle +\langle \phi_{ij}^{ab}, HT_s^2T_d\phi_0\rangle 
\end{equation}
due to the Slater--Condon rules. 
Therefore
\begin{equation}
\label{eq:block_matrix_part_3}
\begin{aligned}
\frac{1}{6}
\langle \phi_{ij}^{ab}, HT^3\phi_0\rangle 
&= 
\frac{1}{6}\langle \phi_{ij}^{ab}, HT_s^3\phi_0\rangle 
+
\frac{1}{6}\langle \phi_{ij}^{ab}, HT_s^2T_d\phi_0\rangle\\
&= 
\frac{1}{6}
\sum_{\substack{k\neq\ell\neq m\\ c\neq d\neq e}} t_k^ct_{\ell}^dt_m^e\langle \phi_{ij}^{ab}, H\phi_{k\ell m}^{cde}\rangle\\
&\quad
+
\frac{1}{6}
\sum_{\substack{k\neq\ell\neq m<n\\ c\neq d\neq e<f}} t_k^ct_{\ell}^dt_{mn}^{ef}\langle \phi_{ij}^{ab}, H\phi_{k\ell mn}^{cdef}\rangle\\
&=
\sum_{\substack{k\neq\ell< m\\ c\neq d< e}} t_k^c(t_{\ell}^dt_m^e-t_{\ell}^et_m^d)\langle \phi_{ij}^{ab}, H\phi_{c.o.(k\ell m)}^{c.o.(cde)}\rangle\\
&\quad+
\sum_{\substack{k<\ell\neq m<n\\ c< d\neq e<f}} (t_k^ct_{\ell}^d-t_k^dt_{\ell}^c)t_{mn}^{ef}\langle \phi_{ij}^{ab}, H\phi_{c.o.(k\ell mn)}^{c.o.(cdef)}\rangle\\
&=
\begin{bmatrix}
1 & {\bf t_s} & {\bf t_d} & {\bf t_{x}}
\end{bmatrix}
\begin{bmatrix} 0&0&0&0\\
0 & 0 & 0 & {\bf H}_{(ij,ab),tr}\\
0 & 0 & 0& {\bf H}_{(ij,ab),qd}\\  0&0&0&0
\end{bmatrix}
\begin{bmatrix}
1\\
{\bf t_s}\\ 
{\bf t_d}\\
{\bf t_{x}}
\end{bmatrix}
\end{aligned}
\end{equation}

\noindent
Next, we consider the term
\begin{equation}
\label{eq:block_matrix_part_4}
\begin{aligned}
\frac{1}{24}
\langle \phi_{ij}^{ab},HT^4\phi_0\rangle 
=
&\frac{1}{24}\langle \phi_{ij}^{ab},HT_s^4\phi_0\rangle\\
=&\frac{1}{24}\sum_{\substack{k\neq l \neq m \neq n\\ c\neq d \neq e \neq f}}
t_k^c t_l^d t_m^e t_n^f
\langle \phi_{ij}^{ab}, H \phi_{klmn}^{cdef} \rangle\\
=&\sum_{\substack{k< l \neq m < n\\ c< d \neq e < f}}
(t_k^c t_l^d-t_k^d t_l^c)( t_m^e t_n^f-t_m^f t_n^e)
\langle \phi_{ij}^{ab}, H \phi_{c.o.(klmn)}^{c.o.(cdef)} \rangle\\
=&\begin{bmatrix}
1 & {\bf t_s} & {\bf t_d} & {\bf t_{x}}
\end{bmatrix}
\begin{bmatrix} 0&0&0&0\\
0 & 0 & 0 &  0\\
0 & 0 & 0& 0\\  0&0&0&{\bf H}_{(ij,ab),qd}
\end{bmatrix}
\begin{bmatrix}
1\\
{\bf t_s}\\ 
{\bf t_d}\\
{\bf t_{x}}
\end{bmatrix}
\end{aligned}
\end{equation}

We may combine Eqs.~\eqref{eq:T2_splie_part_2},\eqref{eq:T2_splie_part_3},\eqref{eq:block_matrix_part_3} and \eqref{eq:block_matrix_part_4} to
\begin{equation}
\begin{aligned}
&\frac{1}{2}
\langle \phi_{ij}^{ab}, H T_sT_d \phi_0 \rangle
+
\frac{1}{2}
\langle \phi_{ij}^{ab}, H T_dT_d \phi_0 \rangle
+
\frac{1}{6}
\langle \phi_{ij}^{ab}, HT_s^3\phi_0\rangle \\
&\quad
+
\frac{1}{6}
\langle \phi_{ij}^{ab}, HT_s^2T_d\phi_0\rangle
+
\frac{1}{24}
\langle \phi_{ij}^{ab},HT_s^4\phi_0\rangle\\
&=
\begin{bmatrix}
1 & {\bf t_s} & {\bf t_d} & {\bf t_{x}}
\end{bmatrix}
\begin{bmatrix} 
0&0&0&0\\
0 & 0 & {\bf H}_{(ij,ab),tr} & {\bf H}_{(ij,ab),tr}\\
0 & 0 & {\bf H}_{(ij,ab),qd}& {\bf H}_{(ij,ab),qd}\\ 
0&0&0&{\bf H}_{(ij,ab),qd}
\end{bmatrix}
\begin{bmatrix}
1\\
{\bf t_s}\\ 
{\bf t_d}\\
{\bf t_{x}}
\end{bmatrix}
\end{aligned}
\end{equation}

\noindent
Next, we investigate the term
\begin{equation}
\label{eq:double_proj_with_explicit_dummy}
\begin{aligned}
&
\frac{1}{2}
\langle \phi_{ij}^{ab}, T^2HT\phi_0\rangle \\
&=
\frac{1}{2}
\langle \phi_{ij}^{ab}, T_s^2HT\phi_0\rangle   \\
&= \sum_{k,c} (t_i^at_j^b-t_i^bt_j^a)t_k^c\langle\phi_0,H\phi_k^c\rangle+
\sum_{\substack{k<m\\c<d}} (t_i^at_j^b-t_i^bt_j^a)t_{km}^{cd}\langle\phi_0,H\phi_{km}^{cd}\rangle\\
&=\begin{bmatrix}
1 & {\bf t_s} & {\bf t_d} & \dots \tilde{t}_{i,j,a,b} \dots
\end{bmatrix}
\begin{bmatrix} 
0&0&0&0\\
0 & 0 & 0 &  0\\
0 & 0 & 0& 0\\
\vdots & \vdots & \vdots & \vdots\\
0&H_{(ij,ab),s}&H_{(ij,ab),d}&0\\
\vdots & \vdots & \vdots & \vdots\\
\end{bmatrix}
\begin{bmatrix}
1\\
{\bf t_s}\\ 
{\bf t_d}\\
{\bf t_{x}}
\end{bmatrix}
\end{aligned}
\end{equation}

\noindent
Note that $\langle \phi_{ij}^{ab}, T^3H T^n\phi_0\rangle=0$ for all $n$.\\

\noindent
Next, we look at the term
\begin{equation}
\label{eq:double_proj_with_explicit_dummy_2}
\begin{aligned}
\frac{1}{2}
\langle
\phi_{ij}^{ab}, T_sT_s H \phi_0
\rangle
&=
(t_i^a t_j^b - t_i^b t_j^a) \langle \phi_0,  H \phi_0 \rangle\\
&=
\begin{bmatrix}
1 & {\bf t_s} & {\bf t_d} & \dots \tilde{t}_{i,j,a,b} \dots
\end{bmatrix}
\begin{bmatrix} 
0&0&0&0\\
0 & 0 & 0 & 0\\
0 & 0 & 0 & 0\\  
\vdots & \vdots & \vdots & \vdots\\
H_{0,0}&0&0&0\\
\vdots & \vdots & \vdots & \vdots\\
\end{bmatrix}
\begin{bmatrix}
1\\
{\bf t_s}\\ 
{\bf t_d}\\
{\bf t_{x}}
\end{bmatrix}
\end{aligned}
\end{equation}

\noindent
For the term $\langle \phi_{ij}^{ab}, T^2HT^2\phi_0\rangle$ we may again apply the Slater--Condon rules which yield
\begin{equation}
\langle \phi_{ij}^{ab}, T^2HT^2\phi_0\rangle = \langle \phi_{ij}^{ab}, T_s^2HT_s^2\phi_0\rangle.
\end{equation}
Therefore
\begin{equation}
\label{eq:double_proj_with_explicit_dummy_3}
\begin{aligned}
\frac{1}{4}
\langle \phi_{ij}^{ab}, T^2 & HT^2 \phi_0\rangle\\
=& 
\frac{1}{4}
\langle \phi_{ij}^{ab}, T_s^2HT_s^2\phi_0\rangle\\
=& 
\frac{1}{2}\sum_{k\neq\ell, c\neq d}(t_i^at_j^b-t_i^bt_j^a)t_k^ct_\ell^d\langle\phi_0,H\phi_{k\ell}^{cd}\rangle\\
=& 
\sum_{k< \ell, c < d}(t_i^at_j^b-t_i^bt_j^a)(t_k^ct_\ell^d-t_k^dt_\ell^c)\langle\phi_0,H\phi_{k\ell}^{cd}\rangle\\
=&
\begin{bmatrix}
1 & {\bf t_s} & {\bf t_d} & \dots \tilde{t}_{i,j,a,b} \dots
\end{bmatrix}
\begin{bmatrix} 0&0&0&0\\
0 & 0 & 0 &  0\\
0 & 0 & 0& 0\\
\vdots & \vdots & \vdots & \vdots\\
0&0&0&H_{(ij,ab),d}\\
\vdots & \vdots & \vdots & \vdots\\
\end{bmatrix}
\begin{bmatrix}
1\\
{\bf t_s}\\ 
{\bf t_d}\\
{\bf t_{x}}
\end{bmatrix}
\end{aligned}
\end{equation}

\noindent
We may combine Eq.~\eqref{eq:double_proj_with_explicit_dummy},~\eqref{eq:double_proj_with_explicit_dummy_2} and~\eqref{eq:double_proj_with_explicit_dummy_3} to
\begin{equation}
\begin{aligned}
&
\frac{1}{2}
\langle \phi_{ij}^{ab}, T_sT_s H \phi_0 \rangle
+
\frac{1}{2}
\langle \phi_{ij}^{ab}, T_s^2HT\phi_0\rangle 
+
\frac{1}{4}
\langle \phi_{ij}^{ab}, T^2  HT^2 \phi_0\rangle\\
&=
\begin{footnotesize}
\begin{bmatrix}
1 & {\bf t_s} & {\bf t_d} & \dots \tilde{t}_{i,j,a,b} \dots
\end{bmatrix}
\begin{bmatrix} 
0&0&0&0\\
0 & 0 & 0 &  0\\
0 & 0 & 0& 0\\  
\vdots & \vdots & \vdots & \vdots\\
H_{0,0}&H_{(ij,ab),s}&H_{(ij,ab),d}&H_{(ij,ab),d}\\
\vdots & \vdots & \vdots & \vdots\\
\end{bmatrix}
\begin{bmatrix}
1\\
{\bf t_s}\\ 
{\bf t_d}\\
{\bf t_{x}}
\end{bmatrix}
\end{footnotesize}
\end{aligned}
\end{equation}

\noindent
Note that 
$\langle
\phi_{ij}^{ab}, T_d H T_sT_d \phi_0
\rangle =0$
due to the Slater--Condon rules. 
Similarly, we find
\begin{equation}
\begin{aligned}
\frac{1}{2}
\langle
\phi_{ij}^{ab}, T_s H T_sT_d \phi_0
\rangle
&=
t_i^a\sum_{k,c}\sum_{\substack{\ell < m \\ d < e}}
t_{k}^ct_{\ell m}^{de}
\langle
\phi_{j}^{b}, H \phi_{c.o.(k\ell m)}^{c.o.(cde)}\rangle\\
&\quad
+
t_j^b\sum_{k,c}\sum_{\substack{\ell < m \\ d < e}}
t_{k}^ct_{\ell m}^{de}
\langle
\phi_{i}^{a}, H \phi_{c.o.(k\ell m)}^{c.o.(cde)}\rangle\\
&\quad 
-
t_i^b\sum_{k,c}\sum_{\substack{\ell < m \\ d < e}}
t_{k}^ct_{\ell m}^{de}
\langle
\phi_{j}^{a}, H \phi_{c.o.(k\ell m)}^{c.o.(cde)}\rangle\\
&\quad
-
t_j^a\sum_{k,c}\sum_{\substack{\ell < m \\ d < e}}
t_{k}^ct_{\ell m}^{de}
\langle
\phi_{i}^{b}, H \phi_{c.o.(k\ell m)}^{c.o.(cde)}\rangle\\
=&
t_i^a \begin{bmatrix}
1 \\
{\bf t_{s}}\\
{\bf t_{d}}\\
{\bf t_{x}}
\end{bmatrix}^\top
\begin{bmatrix} 0&0&0&0\\
0 & 0 & {\bf H_{(j,b),tr}} &  0\\
0 & 0 & 0& 0\\
0&0&0&0\\
\end{bmatrix}
\begin{bmatrix}
1\\
{\bf t_s}\\ 
{\bf t_d}\\
{\bf t_{x}}
\end{bmatrix}\\
&\quad +t_i^b \begin{bmatrix}
1 \\
{\bf t_{s}}\\
{\bf t_{d}}\\
{\bf t_{x}}
\end{bmatrix}^\top
\begin{bmatrix} 0&0&0&0\\
0 & 0 & {\bf -H_{(j,a),tr}} &  0\\
0 & 0 & 0& 0\\
0&0&0&0\\
\end{bmatrix}
\begin{bmatrix}
1\\
{\bf t_s}\\ 
{\bf t_d}\\
{\bf t_{x}}
\end{bmatrix}\\
&\quad+t_j^a \begin{bmatrix}
1 \\
{\bf t_{s}}\\
{\bf t_{d}}\\
{\bf t_{x}}
\end{bmatrix}^\top
\begin{bmatrix} 0&0&0&0\\
0 & 0 & {\bf -H_{(i,b),tr}} &  0\\
0 & 0 & 0& 0\\
0&0&0&0\\
\end{bmatrix}
\begin{bmatrix}
1\\
{\bf t_s}\\ 
{\bf t_d}\\
{\bf t_{x}}
\end{bmatrix}\\
&\quad+ t_j^b\begin{bmatrix}
1 \\
{\bf t_{s}}\\
{\bf t_{d}}\\
{\bf t_{x}}
\end{bmatrix}^\top
\begin{bmatrix} 0&0&0&0\\
0 & 0 & {\bf H_{(i,a),tr}} &  0\\
0 & 0 & 0& 0\\
0&0&0&0\\
\end{bmatrix}
\begin{bmatrix}
1\\
{\bf t_s}\\ 
{\bf t_d}\\
{\bf t_{x}}
\end{bmatrix} \\
& = \underline{\bf H_{(ij,ab),tr,s}^1}\times_1 \begin{bmatrix}
1\\
{\bf t_s}\\ 
{\bf t_d}\\
{\bf t_{x}}
\end{bmatrix}\times_2 \begin{bmatrix}
1\\
{\bf t_s}\\ 
{\bf t_d}\\
{\bf t_{x}}
\end{bmatrix}\times_3 \begin{bmatrix}
1\\
{\bf t_s}\\ 
{\bf t_d}\\
{\bf t_{x}}
\end{bmatrix}
\end{aligned}
\end{equation}

\noindent
where $\underline{\bf H_{(ij,ab),tr,s}^1}$ is the 3-tensor
\begin{equation}
\begin{footnotesize}
\begin{aligned}
    &\underline{\bf H_{(ij,ab),tr,s}}=\\
   &\begin{bmatrix}  ... ,    \begin{bmatrix} 0&0&0&0\\
0 & 0 & {\bf H_{(j,b),tr}} &  0\\
0 & 0 & 0& 0\\
0&0&0&0\\
\end{bmatrix},...,    \begin{bmatrix} 0&0&0&0\\
0 & 0 & {\bf -H_{(j,a),tr}} &  0\\
0 & 0 & 0& 0\\
0&0&0&0\\
\end{bmatrix},... ,   \begin{bmatrix} 0&0&0&0\\
0 & 0 & {\bf -H_{(i,b),tr}} &  0\\
0 & 0 & 0& 0\\
0&0&0&0\\
\end{bmatrix},
\cdots,\begin{bmatrix} 0&0&0&0\\
0 & 0 & {\bf H_{(i,a),tr}} &  0\\
0 & 0 & 0& 0\\
0&0&0&0\\
\end{bmatrix},... \end{bmatrix}   
\end{aligned}
\end{footnotesize}
\end{equation}

\noindent
Note that the Slater--Condon rules imply that 
\begin{equation}
\langle \phi_{ij}^{ab}, THT^2_d\phi_0\rangle=0
\end{equation}
since we have at least one deexcitation applied to $\phi_{ij}^{ab}$ and at least 4 excitations to $\phi_0$.\\

\noindent
Similarly, we get that 
\begin{equation}
\langle \phi_{ij}^{ab}, THT^3\phi_0\rangle = \langle \phi_{ij}^{ab}, T_sHT_s^3\phi_0\rangle
\end{equation}
due to the Slater-Condon rules, and find
\begin{equation}
\begin{aligned}
&
\frac{1}{6}
\langle \phi_{ij}^{ab}, T_sHT_s^3\phi_0\rangle\\
=& 
\frac{1}{6}\sum_{\substack{k\neq\ell\neq m\\ c\neq d\neq e}} t_i^at_k^ct_{\ell}^dt_m^e\langle \phi_{j}^{b}, H\phi_{k\ell m}^{cde}\rangle
-
\frac{1}{6}\sum_{\substack{k\neq\ell\neq m\\ c\neq d\neq e}} t_i^bt_k^ct_{\ell}^dt_m^e\langle \phi_{j}^{a}, H\phi_{k\ell m}^{cde}\rangle\\
&\quad 
-\frac{1}{6}\sum_{\substack{k\neq\ell\neq m\\ c\neq d\neq e}} t_j^at_k^ct_{\ell}^dt_m^e\langle \phi_{i}^{b}, H\phi_{k\ell m}^{cde}\rangle 
+
\frac{1}{6}\sum_{\substack{k\neq\ell\neq m\\ c\neq d\neq e}} t_j^bt_k^ct_{\ell}^dt_m^e\langle \phi_{i}^{a}, H\phi_{k\ell m}^{cde}\rangle\\
=& 
\sum_{\substack{k\neq\ell< m\\ c\neq d< e}} t_i^at_k^c (t_{\ell}^dt_m^e-t_{\ell}^et_m^d)\langle \phi_{j}^{b}, H\phi_{c.o.(k\ell m)}^{c.o.(cde)}\rangle \\
&\quad
-\sum_{\substack{k\neq\ell< m\\ c\neq d<e}} t_i^bt_k^c (t_{\ell}^dt_m^e-t_{\ell}^et_m^d)\langle \phi_{j}^{a}, H\phi_{c.o.(k\ell m)}^{c.o.(cde)}\rangle\\
&\quad 
-\sum_{\substack{k\neq\ell< m\\ c\neq d< e}} t_j^at_k^c (t_{\ell}^dt_m^e-t_{\ell}^et_m^d)\langle \phi_{i}^{b}, H\phi_{c.o.(k\ell m)}^{c.o.(cde)}\rangle\\
&\quad 
+\sum_{\substack{k\neq\ell< m\\ c\neq d< e}} t_j^bt_k^c (t_{\ell}^dt_m^e-t_{\ell}^et_m^d)\langle \phi_{i}^{a}, H\phi_{c.o.(k\ell m)}^{c.o.(cde)}\rangle
\end{aligned}
\end{equation}
\noindent
We may again express this term using a tensor notation 
\begin{equation}
\begin{aligned}
\frac{1}{6}
\langle \phi_{ij}^{ab}, T_sHT_s^3\phi_0\rangle
&=
t_i^a \begin{bmatrix}
1 \\
{\bf t_{s}}\\
{\bf t_{d}}\\
{\bf t_{x}}
\end{bmatrix}^\top
\begin{bmatrix} 0&0&0&0\\
0 & 0 & 0 &  {\bf H_{(j,b),tr}}\\
0 & 0 & 0& 0\\
0&0&0&0\\
\end{bmatrix}
\begin{bmatrix}
1\\
{\bf t_s}\\ 
{\bf t_d}\\
{\bf t_{x}}
\end{bmatrix}\\
&\quad +t_i^b \begin{bmatrix}
1 \\
{\bf t_{s}}\\
{\bf t_{d}}\\
{\bf t_{x}}
\end{bmatrix}^\top
\begin{bmatrix} 0&0&0&0\\
0 & 0 & 0 &  {\bf -H_{(j,a),tr}}\\
0 & 0 & 0& 0\\
0&0&0&0\\
\end{bmatrix}
\begin{bmatrix}
1\\
{\bf t_s}\\ 
{\bf t_d}\\
{\bf t_{x}}
\end{bmatrix}\\
&\quad+t_j^a \begin{bmatrix}
1 \\
{\bf t_{s}}\\
{\bf t_{d}}\\
{\bf t_{x}}
\end{bmatrix}^\top
\begin{bmatrix} 0&0&0&0\\
0 & 0 & 0 &  {\bf -H_{(i,b),tr}}\\
0 & 0 & 0& 0\\
0&0&0&0\\
\end{bmatrix}
\begin{bmatrix}
1\\
{\bf t_s}\\ 
{\bf t_d}\\
{\bf t_{x}}
\end{bmatrix}\\
&\quad+ t_j^b\begin{bmatrix}
1 \\
{\bf t_{s}}\\
{\bf t_{d}}\\
{\bf t_{x}}
\end{bmatrix}^\top
\begin{bmatrix} 0&0&0&0\\
0 & 0 & 0 &  {\bf H_{(i,a),tr}}\\
0 & 0 & 0& 0\\
0&0&0&0\\
\end{bmatrix}
\begin{bmatrix}
1\\
{\bf t_s}\\ 
{\bf t_d}\\
{\bf t_{x}}
\end{bmatrix} \\
& = \underline{\bf H_{(ij,ab),tr,s}^2}\times_1 \begin{bmatrix}
1\\
{\bf t_s}\\ 
{\bf t_d}\\
{\bf t_{x}}
\end{bmatrix}\times_2 \begin{bmatrix}
1\\
{\bf t_s}\\ 
{\bf t_d}\\
{\bf t_{x}}
\end{bmatrix}\times_3 \begin{bmatrix}
1\\
{\bf t_s}\\ 
{\bf t_d}\\
{\bf t_{x}}
\end{bmatrix}
\end{aligned}
\end{equation}
where $\underline{\bf H_{(ij,ab),tr,s}^2}$ is the 3-tensor
\begin{equation}
\begin{footnotesize}
\begin{aligned}
&\underline{\bf H^2_{(ij,ab),tr,s}}=\\
&\begin{bmatrix} 
...,
\begin{bmatrix} 0&0&0&0\\
0 & 0 & 0 &  {\bf H_{(j,b),tr}}\\
0 & 0 & 0& 0\\
0&0&0&0\\
\end{bmatrix}
, ... ,    
\begin{bmatrix} 0&0&0&0\\
0 & 0 & 0&-{\bf H_{(j,a),tr}} \\
0 & 0 & 0& 0\\
0&0&0&0\\
\end{bmatrix}
,...,    
\begin{bmatrix} 0&0&0&0\\
0 & 0 & 0&-{\bf H_{(i,b),tr}} \\
0 & 0 & 0& 0\\
0&0&0&0\\
\end{bmatrix}
,... ,   
\begin{bmatrix} 0&0&0&0\\
0 & 0 & 0&{\bf H_{(i,a),tr}} \\
0 & 0 & 0& 0\\
0&0&0&0\\
\end{bmatrix}
,...
\end{bmatrix}   
\end{aligned}
\end{footnotesize}
\end{equation}

\noindent
We can combine the last two via
\begin{equation}
\begin{aligned}
&
\frac{1}{2}
\langle\phi_{ij}^{ab}, T_s H T_sT_d \phi_0\rangle
+
\frac{1}{6}
\langle \phi_{ij}^{ab}, T_sHT_s^3\phi_0\rangle\\
&= t_i^a \begin{bmatrix}
1 \\
{\bf t_{s}}\\
{\bf t_{d}}\\
{\bf t_{x}}
\end{bmatrix}^\top
\begin{bmatrix} 0&0&0&0\\
0 & 0 & {\bf H_{(j,b),tr}} &  {\bf H_{(j,b),tr}}\\
0 & 0 & 0& 0\\
0&0&0&0\\
\end{bmatrix}
\begin{bmatrix}
1\\
{\bf t_s}\\ 
{\bf t_d}\\
{\bf t_{x}}
\end{bmatrix}\\
&\quad +t_i^b \begin{bmatrix}
1 \\
{\bf t_{s}}\\
{\bf t_{d}}\\
{\bf t_{x}}
\end{bmatrix}^\top
\begin{bmatrix} 0&0&0&0\\
0 & 0 & {\bf -H_{(j,a),tr}} &  {\bf -H_{(j,a),tr}}\\
0 & 0 & 0& 0\\
0&0&0&0\\
\end{bmatrix}
\begin{bmatrix}
1\\
{\bf t_s}\\ 
{\bf t_d}\\
{\bf t_{x}}
\end{bmatrix}\\
&\quad+t_j^a \begin{bmatrix}
1 \\
{\bf t_{s}}\\
{\bf t_{d}}\\
{\bf t_{x}}
\end{bmatrix}^\top
\begin{bmatrix} 0&0&0&0\\
0 & 0 & {\bf -H_{(i,b),tr}} &  {\bf -H_{(i,b),tr}}\\
0 & 0 & 0& 0\\
0&0&0&0\\
\end{bmatrix}
\begin{bmatrix}
1\\
{\bf t_s}\\ 
{\bf t_d}\\
{\bf t_{x}}
\end{bmatrix}\\
&\quad+ t_j^b\begin{bmatrix}
1 \\
{\bf t_{s}}
{\bf t_{d}}\\
{\bf t_{x}}
\end{bmatrix}^\top
\begin{bmatrix} 0&0&0&0\\
0 & 0 & {\bf H_{(i,a),tr}} &  {\bf H_{(i,a),tr}}\\
0 & 0 & 0& 0\\
0&0&0&0\\
\end{bmatrix}
\begin{bmatrix}
1\\
{\bf t_s}\\ 
{\bf t_d}\\
{\bf t_{x}}
\end{bmatrix}\\
&=
(\underline{\bf H_{(ij,ab),tr,s}^1}+ \underline{\bf H_{(ij,ab),tr,s}^2})\times_1 \begin{bmatrix}
1\\
{\bf t_s}\\ 
{\bf t_d}\\
{\bf t_{x}}
\end{bmatrix}\times_2 \begin{bmatrix}
1\\
{\bf t_s}\\ 
{\bf t_d}\\
{\bf t_{x}}
\end{bmatrix}\times_3 \begin{bmatrix}
1\\
{\bf t_s}\\ 
{\bf t_d}\\
{\bf t_{x}}
\end{bmatrix}
\end{aligned}
\end{equation}

\subsection{Two Electrons in Four Spin Orbitals}
\label{app:2in4}

We shall here derive the individual polynomials for the two electrons in four spin orbitals example.
\paragraph{First polynomial}
The first polynomial is given by 
\begin{equation}
f_1(t)=\langle \Phi_2^3, e^{-T}H e^T \Phi_0  \rangle
\end{equation}

Going through this term by term, only focusing on the degree of the individual monomials we get:

\begin{equation}
\begin{aligned}
\langle \Phi_2^3, H \Phi_0  \rangle = C
\end{aligned}
\end{equation}
which corresponds to $(0,0,0,0,0)$.
The singly nested commutator yields the following
\begin{equation}
\begin{aligned}
\langle \Phi_2^3, [H,T] \Phi_0  \rangle &=
 t_2^3 \langle \Phi_2^3, H \Phi_2^3  \rangle
+t_2^4 \langle \Phi_2^3, H \Phi_2^4  \rangle
+t_1^3 \langle \Phi_2^3, H \Phi_1^3  \rangle\\
&\quad+t_1^4 \langle \Phi_2^3, H \Phi_1^4  \rangle
+t_{12}^{34} \langle \Phi_2^3, H \Phi_{12}^{34}  \rangle
-t_2^3 \langle \Phi_0, H \Phi_0  \rangle
\end{aligned}
\end{equation}
which corresponds to:
\begin{equation}
\begin{aligned}
(1,0,0,0,0)
,(0,1,0,0,0)
,(0,0,1,0,0)
,(0,0,0,1,0)
,(0,0,0,0,1)
\end{aligned}
\end{equation}
The doubly nested commutator yields
\begin{equation}
\begin{aligned}
\langle \Phi_2^3, [[H,T],T] \Phi_0  \rangle &=
  \langle \Phi_2^3, HT^2 \Phi_0 \rangle
-2\langle \Phi_2^3, THT \Phi_0 \rangle
+ \langle \Phi_2^3, T^2H \Phi_0 \rangle\\
&=
(t_1^3t_2^4-t_2^3t_1^4)  \langle \Phi_2^3, H \Phi_{12}^{34}  \rangle 
-2 t_2^3 t_2^3 \langle \Phi_0, H \Phi_2^3  \rangle\\
&\quad
-2 t_2^3 t_2^4 \langle \Phi_0, H \Phi_2^4  \rangle
-2 t_2^3 t_1^3 \langle \Phi_0, H \Phi_1^3  \rangle\\
&\quad
-2 t_2^3 t_1^4 \langle \Phi_0, H \Phi_1^4  \rangle
-2 t_2^3 t_{12}^{34} \langle \Phi_0, H \Phi_{12}^{34}  \rangle
\end{aligned}
\end{equation}
which corresponds to:
\begin{equation}
\begin{aligned}
(1,0,0,1,0)
,(0,1,1,0,0)
,(2,0,0,0,0)
,\\(1,1,0,0,0)
,(1,0,1,0,0)
,(1,0,0,1,0)
,(1,0,0,0,1)
\end{aligned}
\end{equation}
The triply nested commutator yields
\begin{equation}
\begin{aligned}
\langle \Phi_2^3, [[[H,T],T],T] \Phi_0  \rangle 
&=
\langle \Phi_2^3, HT^3 \Phi_0 \rangle
-3\langle \Phi_2^3, THT^2 \Phi_0 \rangle\\
&\quad -3\langle \Phi_2^3, T^2HT \Phi_0 \rangle
+ \langle \Phi_2^3, T^3H \Phi_0 \rangle\\
&=-3t_2^3  \langle \Phi_0, HT^2 \Phi_0  \rangle\\
&=-3t_2^3(t_1^3t_2^4 -t_2^3t_1^4)\langle \Phi_0,  H \Phi_{12}^{34}  \rangle
\end{aligned}
\end{equation}
which corresponds to:
\begin{equation}
\begin{aligned}
(2,0,0,1,0),(1,1,1,0,0)
\end{aligned}
\end{equation}
Overall the Newton polytope that corresponds to the first polynomial is given by
\begin{equation}
\begin{aligned}
{\rm New}_2^3
={\rm conv } \big(&
(0,0,0,0,0),~
(1,0,0,0,0),~
(0,1,0,0,0),~
(0,0,1,0,0),~\\
&(0,0,0,1,0),~
(0,0,0,0,1),~
(1,0,0,1,0),~
(0,1,1,0,0),~\\
&(2,0,0,0,0),~
(1,1,0,0,0),~
(1,0,1,0,0),~
(1,0,0,1,0),~\\
&(1,0,0,0,1),~
(2,0,0,1,0),~
(1,1,1,0,0)
\big)
\end{aligned}
\end{equation}
The first polynomial can be then summarized as
\begin{equation}
\begin{aligned}
f_1(t) &= 
\langle \Phi_2^3, H \Phi_0  \rangle
+t_2^3 \langle \Phi_2^3, H \Phi_2^3  \rangle
+t_2^4 \langle \Phi_2^3, H \Phi_2^4  \rangle\\
&\quad
+t_1^3 \langle \Phi_2^3, H \Phi_1^3  \rangle
+t_1^4 \langle \Phi_2^3, H \Phi_1^4  \rangle
+t_{12}^{34} \langle \Phi_2^3, H \Phi_{12}^{34}  \rangle\\
&\quad
-t_2^3 \langle \Phi_0, H \Phi_0  \rangle
+(t_1^3t_2^4-t_2^3t_1^4)  \langle \Phi_2^3, H \Phi_{12}^{34}  \rangle 
-2 t_2^3 t_2^3 \langle \Phi_0, H \Phi_2^3  \rangle\\
&\quad
-2 t_2^3 t_2^4 \langle \Phi_0, H \Phi_2^4  \rangle
-2 t_2^3 t_1^3 \langle \Phi_0, H \Phi_1^3  \rangle
-2 t_2^3 t_1^4 \langle \Phi_0, H \Phi_1^4  \rangle\\
&\quad
-2 t_2^3 t_{12}^{34} \langle \Phi_0, H \Phi_{12}^{34}  \rangle
-3t_2^3(t_1^3t_2^4 -t_2^3t_1^4)\langle \Phi_0,  H \Phi_{12}^{34}  \rangle
\end{aligned}
\end{equation}

\paragraph{Second polynomial}
The second polynomial is given by 
\begin{equation}
f_2(t)=\langle \Phi_2^4, e^{-T}H e^T \Phi_0  \rangle
\end{equation}

Going through this term by term, only focusing on the degree of the individual monomials we get:

\begin{equation}
\begin{aligned}
\langle \Phi_2^4, H \Phi_0  \rangle = C
\end{aligned}
\end{equation}
which corresponds to $(0,0,0,0,0)$.
The singly nested commutator yields the following
\begin{equation}
\begin{aligned}
\langle \Phi_2^4, [H,T] \Phi_0  \rangle &=
 t_2^3 \langle \Phi_2^4, H \Phi_2^3  \rangle
+t_2^4 \langle \Phi_2^4, H \Phi_2^4  \rangle
+t_1^3 \langle \Phi_2^4, H \Phi_1^3  \rangle\\
&\quad+t_1^4 \langle \Phi_2^4, H \Phi_1^4  \rangle
+t_{12}^{34} \langle \Phi_2^4, H \Phi_{12}^{34}  \rangle
-t_2^4 \langle \Phi_0, H \Phi_0  \rangle
\end{aligned}
\end{equation}
which corresponds to:
\begin{equation}
\begin{aligned}
(1,0,0,0,0)
,(0,1,0,0,0)
,(0,0,1,0,0)
,(0,0,0,1,0)
,(0,0,0,0,1)
\end{aligned}
\end{equation}
The doubly nested commutator yields
\begin{equation}
\begin{aligned}
\langle \Phi_2^4, [[H,T],T] \Phi_0  \rangle &=
  \langle \Phi_2^4, HT^2 \Phi_0 \rangle
-2\langle \Phi_2^4, THT \Phi_0 \rangle
+ \langle \Phi_2^4, T^2H \Phi_0 \rangle\\
&=
(t_1^3t_2^4-t_2^3t_1^4)  \langle \Phi_2^4, H \Phi_{12}^{34}  \rangle 
-2 t_2^4 t_2^3 \langle \Phi_0, H \Phi_2^3  \rangle\\
&\quad
-2 t_2^4 t_2^4 \langle \Phi_0, H \Phi_2^4  \rangle
-2 t_2^4 t_1^3 \langle \Phi_0, H \Phi_1^3  \rangle\\
&\quad
-2 t_2^4 t_1^4 \langle \Phi_0, H \Phi_1^4  \rangle
-2 t_2^4 t_{12}^{34} \langle \Phi_0, H \Phi_{12}^{34}  \rangle
\end{aligned}
\end{equation}
which corresponds to:
\begin{equation}
\begin{aligned}
(1,0,0,1,0)
,(0,1,1,0,0)
,(1,1,0,0,0)
,\\(0,2,0,0,0)
,(0,1,1,0,0)
,(0,1,0,1,0)
,(0,1,0,0,1)
\end{aligned}
\end{equation}
The triply nested commutator yields
\begin{equation}
\begin{aligned}
\langle \Phi_2^4, [[[H,T],T],T] \Phi_0  \rangle 
&=
\langle \Phi_2^4, HT^3 \Phi_0 \rangle
-3\langle \Phi_2^4, THT^2 \Phi_0 \rangle\\
&\quad -3\langle \Phi_2^4, T^2HT \Phi_0 \rangle
+ \langle \Phi_2^4, T^3H \Phi_0 \rangle\\
&=-3t_2^4  \langle \Phi_0, HT^2 \Phi_0  \rangle\\
&=-3t_2^4(t_1^3t_2^4 -t_2^3t_1^4)\langle \Phi_0,  H \Phi_{12}^{34}  \rangle
\end{aligned}
\end{equation}
which corresponds to:
\begin{equation}
\begin{aligned}
(1,1,0,1,0),(0,2,1,0,0)
\end{aligned}
\end{equation}
Overall the Newton polytope that corresponds to the second polynomial is given by
\begin{equation}
\begin{aligned}
{\rm New}_2^4
={\rm conv } \big(&
(0,0,0,0,0),~
(1,0,0,0,0),~
(0,1,0,0,0),~
(0,0,1,0,0),~\\
&(0,0,0,1,0),~
(0,0,0,0,1),~
(1,0,0,1,0),~
(0,1,1,0,0),~\\
&(1,1,0,0,0),~
(0,2,0,0,0),~
(0,1,1,0,0),~
(0,1,0,1,0),~\\
&(1,0,0,0,1),~
(1,1,0,1,0),~
(0,2,1,0,0)
\big)
\end{aligned}
\end{equation}
The second polynomial can be then summarized as
\begin{equation}
\begin{aligned}
f_2(t) &= 
\langle \Phi_2^4, H \Phi_0  \rangle
+t_2^3 \langle \Phi_2^4, H \Phi_2^3  \rangle
+t_2^4 \langle \Phi_2^4, H \Phi_2^4  \rangle\\
&\quad
+t_1^3 \langle \Phi_2^4, H \Phi_1^3  \rangle
+t_1^4 \langle \Phi_2^4, H \Phi_1^4  \rangle
+t_{12}^{34} \langle \Phi_2^4, H \Phi_{12}^{34}  \rangle\\
&\quad
-t_2^4 \langle \Phi_0, H \Phi_0  \rangle
+(t_1^3t_2^4-t_2^3t_1^4)  \langle \Phi_2^4, H \Phi_{12}^{34}  \rangle 
- t_2^3 t_2^3 \langle \Phi_0, H \Phi_2^3  \rangle\\
&\quad
- t_2^4 t_2^4 \langle \Phi_0, H \Phi_2^4  \rangle
- t_2^4 t_1^3 \langle \Phi_0, H \Phi_1^3  \rangle
- t_2^4 t_1^4 \langle \Phi_0, H \Phi_1^4  \rangle\\
&\quad
- t_2^4 t_{12}^{34} \langle \Phi_0, H \Phi_{12}^{34}  \rangle
-t_2^4(t_1^3t_2^4 -t_2^3t_1^4)\langle \Phi_0,  H \Phi_{12}^{34}  \rangle
\end{aligned}
\end{equation}

\newpage

\paragraph{Third polynomial}
The third polynomial is given by 
\begin{equation}
f_3(t)=\langle \Phi_1^3, e^{-T}H e^T \Phi_0  \rangle
\end{equation}

Going through this term by term, only focusing on the degree of the individual monomials we get:

\begin{equation}
\begin{aligned}
\langle \Phi_1^3, H \Phi_0  \rangle = C
\end{aligned}
\end{equation}
which corresponds to $(0,0,0,0,0)$.
The singly nested commutator yields the following
\begin{equation}
\begin{aligned}
\langle \Phi_1^3, [H,T] \Phi_0  \rangle &=
 t_2^3 \langle \Phi_1^3, H \Phi_2^3  \rangle
+t_2^4 \langle \Phi_1^3, H \Phi_2^4  \rangle
+t_1^3 \langle \Phi_1^3, H \Phi_1^3  \rangle\\
&\quad+t_1^4 \langle \Phi_1^3, H \Phi_1^4  \rangle
+t_{12}^{34} \langle \Phi_1^3, H \Phi_{12}^{34}  \rangle
-t_1^3 \langle \Phi_0, H \Phi_0  \rangle
\end{aligned}
\end{equation}
which corresponds to:
\begin{equation}
\begin{aligned}
(1,0,0,0,0)
,(0,1,0,0,0)
,(0,0,1,0,0)
,(0,0,0,1,0)
,(0,0,0,0,1)
\end{aligned}
\end{equation}
The doubly nested commutator yields
\begin{equation}
\begin{aligned}
\langle \Phi_1^3, [[H,T],T] \Phi_0  \rangle &=
  \langle \Phi_1^3, HT^2 \Phi_0 \rangle
-2\langle \Phi_1^3, THT \Phi_0 \rangle
+ \langle \Phi_1^3, T^2H \Phi_0 \rangle\\
&=
(t_1^3t_2^4-t_2^3t_1^4)  \langle \Phi_1^3, H \Phi_{12}^{34}  \rangle 
-2 t_1^3 t_2^3 \langle \Phi_0, H \Phi_2^3  \rangle\\
&\quad
-2 t_1^3 t_2^4 \langle \Phi_0, H \Phi_2^4  \rangle
-2 t_1^3 t_1^3 \langle \Phi_0, H \Phi_1^3  \rangle\\
&\quad
-2 t_1^3 t_1^4 \langle \Phi_0, H \Phi_1^4  \rangle
-2 t_1^3 t_{12}^{34} \langle \Phi_0, H \Phi_{12}^{34}  \rangle
\end{aligned}
\end{equation}
which corresponds to:
\begin{equation}
\begin{aligned}
(1,0,0,1,0)
,(0,1,1,0,0)
,(1,0,1,0,0)
,\\(0,1,1,0,0)
,(0,0,2,0,0)
,(0,0,1,1,0)
,(0,0,1,0,1)
\end{aligned}
\end{equation}
The triply nested commutator yields
\begin{equation}
\begin{aligned}
\langle \Phi_1^3, [[[H,T],T],T] \Phi_0  \rangle 
&=
\langle \Phi_1^3, HT^3 \Phi_0 \rangle
-3\langle \Phi_1^3, THT^2 \Phi_0 \rangle\\
&\quad -3\langle \Phi_1^3, T^2HT \Phi_0 \rangle
+ \langle \Phi_1^3, T^3H \Phi_0 \rangle\\
&=-3t_1^3  \langle \Phi_0, HT^2 \Phi_0  \rangle\\
&=-3t_1^3(t_1^3t_2^4 -t_2^3t_1^4)\langle \Phi_0,  H \Phi_{12}^{34}  \rangle
\end{aligned}
\end{equation}
which corresponds to:
\begin{equation}
\begin{aligned}
(1,0,1,1,0),(0,1,2,0,0)
\end{aligned}
\end{equation}
Overall the Newton polytope that corresponds to the first polynomial is given by
\begin{equation}
\begin{aligned}
{\rm New}_1^3
={\rm conv } \big(&
(0,0,0,0,0),~
(1,0,0,0,0),~
(0,1,0,0,0),~
(0,0,1,0,0),~\\
&(0,0,0,1,0),~
(0,0,0,0,1),~
(1,0,0,1,0),~
(0,1,1,0,0),~\\
&(1,0,1,0,0),~
(0,1,1,0,0),~
(0,0,2,0,0),~
(0,0,1,1,0),~\\
&(0,0,1,0,1),~
(1,0,1,1,0),~
(0,1,2,0,0)
\big)
\end{aligned}
\end{equation}
The third polynomial can be then summarized as
\begin{equation}
\begin{aligned}
f_3(t) &= 
\langle \Phi_1^3, H \Phi_0  \rangle
+t_2^3 \langle \Phi_1^3, H \Phi_2^3  \rangle
+t_2^4 \langle \Phi_1^3, H \Phi_2^4  \rangle\\
&\quad
+t_1^3 \langle \Phi_1^3, H \Phi_1^3  \rangle
+t_1^4 \langle \Phi_1^3, H \Phi_1^4  \rangle
+t_{12}^{34} \langle \Phi_1^3, H \Phi_{12}^{34}  \rangle\\
&\quad
-t_1^3 \langle \Phi_0, H \Phi_0  \rangle
+(t_1^3t_2^4-t_2^3t_1^4)  \langle \Phi_1^3, H \Phi_{12}^{34}  \rangle 
- t_1^3 t_2^3 \langle \Phi_0, H \Phi_2^3  \rangle\\
&\quad
- t_1^3 t_2^4 \langle \Phi_0, H \Phi_2^4  \rangle
- t_1^3 t_1^3 \langle \Phi_0, H \Phi_1^3  \rangle
- t_1^3 t_1^4 \langle \Phi_0, H \Phi_1^4  \rangle\\
&\quad
- t_1^3 t_{12}^{34} \langle \Phi_0, H \Phi_{12}^{34}  \rangle
-t_1^3(t_1^3t_2^4 -t_2^3t_1^4)\langle \Phi_0,  H \Phi_{12}^{34}  \rangle
\end{aligned}
\end{equation}

\newpage

\paragraph{Fourth polynomial}
The fourth polynomial is given by 
\begin{equation}
f_4(t)=\langle \Phi_1^4, e^{-T}H e^T \Phi_0  \rangle
\end{equation}

Going through this term by term, only focusing on the degree of the individual monomials we get:

\begin{equation}
\begin{aligned}
\langle \Phi_1^4, H \Phi_0  \rangle = C
\end{aligned}
\end{equation}
which corresponds to $(0,0,0,0,0)$.
The singly nested commutator yields the following
\begin{equation}
\begin{aligned}
\langle \Phi_1^4, [H,T] \Phi_0  \rangle &=
 t_2^3 \langle \Phi_1^4, H \Phi_2^3  \rangle
+t_2^4 \langle \Phi_1^4, H \Phi_2^4  \rangle
+t_1^3 \langle \Phi_1^4, H \Phi_1^3  \rangle\\
&\quad+t_1^4 \langle \Phi_1^4, H \Phi_1^4  \rangle
+t_{12}^{34} \langle \Phi_1^4, H \Phi_{12}^{34}  \rangle
-t_1^4 \langle \Phi_0, H \Phi_0  \rangle
\end{aligned}
\end{equation}
which corresponds to:
\begin{equation}
\begin{aligned}
(1,0,0,0,0)
,(0,1,0,0,0)
,(0,0,1,0,0)
,(0,0,0,1,0)
,(0,0,0,0,1)
\end{aligned}
\end{equation}
The doubly nested commutator yields
\begin{equation}
\begin{aligned}
\langle \Phi_1^4, [[H,T],T] \Phi_0  \rangle &=
  \langle \Phi_1^4, HT^2 \Phi_0 \rangle
-2\langle \Phi_1^4, THT \Phi_0 \rangle
+ \langle \Phi_1^4, T^2H \Phi_0 \rangle\\
&=
(t_1^3t_2^4-t_2^3t_1^4)  \langle \Phi_1^4, H \Phi_{12}^{34}  \rangle 
-2 t_1^4 t_2^3 \langle \Phi_0, H \Phi_2^3  \rangle\\
&\quad
-2 t_1^4 t_2^4 \langle \Phi_0, H \Phi_2^4  \rangle
-2 t_1^4 t_1^3 \langle \Phi_0, H \Phi_1^3  \rangle\\
&\quad
-2 t_1^4 t_1^4 \langle \Phi_0, H \Phi_1^4  \rangle
-2 t_1^4 t_{12}^{34} \langle \Phi_0, H \Phi_{12}^{34}  \rangle
\end{aligned}
\end{equation}
which corresponds to:
\begin{equation}
\begin{aligned}
(1,0,0,1,0)
,(0,1,1,0,0)
,(1,0,0,1,0)
,\\(0,1,0,1,0)
,(0,0,1,1,0)
,(0,0,0,2,0)
,(0,0,0,1,1)
\end{aligned}
\end{equation}
The triply nested commutator yields
\begin{equation}
\begin{aligned}
\langle \Phi_1^4, [[[H,T],T],T] \Phi_0  \rangle 
&=
\langle \Phi_1^4, HT^3 \Phi_0 \rangle
-3\langle \Phi_1^4, THT^2 \Phi_0 \rangle\\
&\quad -3\langle \Phi_1^4, T^2HT \Phi_0 \rangle
+ \langle \Phi_1^4, T^3H \Phi_0 \rangle\\
&=-3t_1^4  \langle \Phi_0, HT^2 \Phi_0  \rangle\\
&=-3t_1^4(t_1^3t_2^4 -t_2^3t_1^4)\langle \Phi_0,  H \Phi_{12}^{34}  \rangle
\end{aligned}
\end{equation}
which corresponds to:
\begin{equation}
\begin{aligned}
(1,0,0,2,0),(0,1,1,1,0)
\end{aligned}
\end{equation}
Overall the Newton polytope that corresponds to the second polynomial is given by
\begin{equation}
\begin{aligned}
{\rm New}_2^4
={\rm conv } \big(&
(0,0,0,0,0),~
(1,0,0,0,0),~
(0,1,0,0,0),~
(0,0,1,0,0),~\\
&(0,0,0,1,0),~
(0,0,0,0,1),~
(1,0,0,1,0),~
(0,1,1,0,0),~\\
&(1,0,0,1,0),~
(0,2,0,1,0),~
(0,0,1,1,0),~
(0,0,0,2,0),~\\
&(0,0,0,1,1),~
(1,0,0,2,0),~
(0,1,1,1,0)
\big)
\end{aligned}
\end{equation}
The fourth polynomial can be then summarized as
\begin{equation}
\begin{aligned}
f_4(t) &= 
\langle \Phi_1^4, H \Phi_0  \rangle
+t_2^3 \langle \Phi_1^4, H \Phi_2^3  \rangle
+t_2^4 \langle \Phi_1^4, H \Phi_2^4  \rangle\\
&\quad
+t_1^3 \langle \Phi_1^4, H \Phi_1^3  \rangle
+t_1^4 \langle \Phi_1^4, H \Phi_1^4  \rangle
+t_{12}^{34} \langle \Phi_1^4, H \Phi_{12}^{34}  \rangle\\
&\quad
-t_1^4 \langle \Phi_0, H \Phi_0  \rangle
+(t_1^3t_2^4-t_2^3t_1^4)  \langle \Phi_1^4, H \Phi_{12}^{34}  \rangle 
- t_1^3 t_2^3 \langle \Phi_0, H \Phi_2^3  \rangle\\
&\quad
- t_1^4 t_2^4 \langle \Phi_0, H \Phi_2^4  \rangle
- t_1^4 t_1^3 \langle \Phi_0, H \Phi_1^3  \rangle
- t_1^4 t_1^4 \langle \Phi_0, H \Phi_1^4  \rangle\\
&\quad
- t_1^4 t_{12}^{34} \langle \Phi_0, H \Phi_{12}^{34}  \rangle
-t_1^4(t_1^3t_2^4 -t_2^3t_1^4)\langle \Phi_0,  H \Phi_{12}^{34}  \rangle
\end{aligned}
\end{equation}

\newpage

\paragraph{fifth polynomial}

The fifth polynomial is given by 
\begin{equation}
f_5(t)=\langle \Phi_{12}^{34}, e^{-T}H e^T \Phi_0  \rangle
\end{equation}

Going through this term by term, only focusing on the degree of the individual monomials we get:

\begin{equation}
\begin{aligned}
\langle\Phi_{12}^{34}, H \Phi_0  \rangle = C
\end{aligned}
\end{equation}
which corresponds to $(0,0,0,0,0)$.
The singly nest commutator yields
\begin{equation}
\begin{aligned}
\langle \Phi_{12}^{34}, [H,T] \Phi_0  \rangle &=
 t_2^3 \langle \Phi_{12}^{34}, H \Phi_2^3  \rangle
+t_2^4 \langle \Phi_{12}^{34}, H \Phi_2^4  \rangle
+t_1^3 \langle \Phi_{12}^{34}, H \Phi_1^3  \rangle\\
&\quad+t_1^4 \langle \Phi_{12}^{34}, H \Phi_1^4  \rangle
+t_{12}^{34} \langle \Phi_{12}^{34}, H \Phi_{12}^{34}  \rangle\\
&\quad -t_{12}^{34} \langle \Phi_0, H \Phi_0  \rangle
-t_1^3 \langle \Phi_{2}^{4}, H \Phi_0  \rangle
-t_2^4 \langle \Phi_{1}^{3}, H \Phi_0  \rangle\\
&\quad +t_1^4 \langle \Phi_{2}^{3}, H \Phi_0  \rangle
+t_2^3 \langle \Phi_{1}^{4}, H \Phi_0  \rangle
\end{aligned}
\end{equation}
which corresponds to:
\begin{equation}
\begin{aligned}
(1,0,0,0,0),
(0,1,0,0,0),
(0,0,1,0,0),
(0,0,0,1,0),
(0,0,0,0,1)
\end{aligned}
\end{equation}
The doubly nest commutator yields
\begin{equation}
\begin{aligned}
\langle \Phi_{12}^{34}, [[H,T],T] \Phi_0  \rangle 
&=
\langle \Phi_{12}^{34}, HT^2 \Phi_0 \rangle
-2\langle \Phi_{12}^{34}, THT \Phi_0 \rangle
+ \langle \Phi_{12}^{34}, T^2H \Phi_0 \rangle\\
&=
\big(t_1^3t_2^4- t_2^3t_1^4\big)  \langle \Phi_{12}^{34}, H \Phi_{12}^{34}  \rangle \\
&\quad -2\big ( 
t_{12}^{34} \langle \Phi_0, H T \Phi_0  \rangle
+t_1^3 \langle \Phi_{2}^{4}, H T \Phi_0  \rangle
+t_2^4 \langle \Phi_{1}^{3}, H T \Phi_0  \rangle\\
&\qquad-t_1^4 \langle \Phi_{2}^{3}, H T \Phi_0  \rangle
-t_2^3 \langle \Phi_{1}^{4}, H T \Phi_0  \rangle
\big)\\
&\quad+ \big(t_2^3t_1^4- t_1^3t_2^4\big) \langle \Phi_0, H \Phi_0  \rangle
\end{aligned}
\end{equation}
The terms without $T$ yield
\begin{equation}
\begin{aligned}
(1,0,0,1,0),
(0,1,1,0,0)
\end{aligned}
\end{equation}
We now expand each of the individual summands containing $T$:
\begin{equation}
\begin{aligned}
t_{12}^{34}\langle \Phi_0, H T \Phi_0  \rangle
&=
t_{12}^{34}t_2^3 \langle \Phi_0, H \Phi_2^3  \rangle
+t_{12}^{34}t_2^4 \langle \Phi_0, H \Phi_2^4  \rangle\\
&\quad +t_{12}^{34}t_1^3 \langle \Phi_0, H \Phi_1^3  \rangle
+t_{12}^{34}t_1^4 \langle \Phi_0, H \Phi_1^4  \rangle
+t_{12}^{34}t_{12}^{34} \langle \Phi_0, H \Phi_{12}^{34}  \rangle
\end{aligned}
\end{equation}
which yields
\begin{equation}
\begin{aligned}
(1,0,0,0,1),
(0,1,0,0,1),
(0,0,1,0,1),
(0,0,0,1,1),
(0,0,0,0,2)
\end{aligned}
\end{equation}
Next we find
\begin{equation}
\begin{aligned}
t_1^3\langle \Phi_2^4, H T \Phi_0  \rangle
&=
t_1^3t_2^3 \langle \Phi_2^4, H \Phi_2^3  \rangle
+t_1^3t_2^4 \langle \Phi_2^4, H \Phi_2^4  \rangle\\
&\quad +t_1^3t_1^3 \langle \Phi_2^4, H \Phi_1^3  \rangle
+t_1^3t_1^4 \langle \Phi_2^4, H \Phi_1^4  \rangle
+t_1^3t_{12}^{34} \langle \Phi_2^4, H \Phi_{12}^{34}  \rangle
\end{aligned}
\end{equation}
which yields
\begin{equation}
\begin{aligned}
(1,0,1,0,0),
(0,1,1,0,0),
(0,0,2,0,0),
(0,0,1,1,0),
(0,0,1,0,1)
\end{aligned}
\end{equation}
Next we find
\begin{equation}
\begin{aligned}
t_2^4\langle \Phi_1^3, H T \Phi_0  \rangle
&=
t_2^4t_2^3 \langle \Phi_1^3, H \Phi_2^3  \rangle
+t_2^4t_2^4 \langle \Phi_1^3, H \Phi_2^4  \rangle\\
&\quad +t_2^4t_1^3 \langle \Phi_1^3, H \Phi_1^3  \rangle
+t_2^4t_1^4 \langle \Phi_1^3, H \Phi_1^4  \rangle
+t_2^4t_{12}^{34} \langle \Phi_1^3, H \Phi_{12}^{34}  \rangle
\end{aligned}
\end{equation}
which yields
\begin{equation}
\begin{aligned}
(1,1,0,0,0),
(0,2,0,0,0),
(0,1,1,0,0),
(0,1,0,1,0),
(0,1,0,0,1)
\end{aligned}
\end{equation}
Next we find
\begin{equation}
\begin{aligned}
t_1^4\langle \Phi_2^3, H T \Phi_0  \rangle
&=
t_1^4t_2^3 \langle \Phi_2^3, H \Phi_2^3  \rangle
+t_1^4t_2^4 \langle \Phi_2^3, H \Phi_2^4  \rangle\\
&\quad +t_1^4t_1^3 \langle \Phi_2^3, H \Phi_1^3  \rangle
+t_1^4t_1^4 \langle \Phi_2^3, H \Phi_1^4  \rangle
+t_1^4t_{12}^{34} \langle \Phi_2^3, H \Phi_{12}^{34}  \rangle
\end{aligned}
\end{equation}
which yields
\begin{equation}
\begin{aligned}
(1,0,0,1,0),
(0,1,0,1,0),
(0,0,1,1,0),
(0,0,0,2,0),
(0,0,0,1,1)
\end{aligned}
\end{equation}
Next we find
\begin{equation}
\begin{aligned}
t_2^3\langle \Phi_1^4, H T \Phi_0  \rangle
&=
t_2^3t_2^3 \langle \Phi_1^4, H \Phi_2^3  \rangle
+t_2^3t_2^4 \langle \Phi_1^4, H \Phi_2^4  \rangle\\
&\quad +t_2^3t_1^3 \langle \Phi_1^4, H \Phi_1^3  \rangle
+t_2^3t_1^4 \langle \Phi_1^4, H \Phi_1^4  \rangle
+t_2^3t_{12}^{34} \langle \Phi_1^4, H \Phi_{12}^{34}  \rangle
\end{aligned}
\end{equation}
which yields
\begin{equation}
\begin{aligned}
(2,0,0,0,0),
(1,1,0,0,0),
(1,0,1,0,0),
(1,0,0,1,0),
(1,0,0,0,1)
\end{aligned}
\end{equation}

The triply nest commutator yields
\begin{equation}
\begin{aligned}
\langle \Phi_{12}^{34}, [[[H,T],T],T] \Phi_0  \rangle &=
  \langle \Phi_{12}^{34}, HT^3 \Phi_0 \rangle
-3\langle \Phi_{12}^{34}, THT^2 \Phi_0 \rangle\\
&\quad -3\langle \Phi_{12}^{34}, T^2HT \Phi_0 \rangle
+ \langle \Phi_{12}^{34}, T^3H \Phi_0 \rangle\\
&=-3 (t_1^3t_2^4- t_2^3t_1^4 ) \langle \Phi_{12}^{34}, TH \Phi_{12}^{34} \rangle\\
&\quad -3(t_1^3t_2^4- t_2^3t_1^4 ) \langle \Phi_0, HT \Phi_0\rangle
\end{aligned}
\end{equation}
We can now expand the summands containing $T$
\begin{equation}
\begin{aligned}
(t_1^3t_2^4-t_2^3t_1^4 )\langle \Phi_0, H T \Phi_0  \rangle
&=
(t_1^3t_2^4-t_2^3t_1^4 )t_2^3 \langle \Phi_0, H \Phi_2^3  \rangle
+(t_1^3t_2^4-t_2^3t_1^4 )t_2^4 \langle \Phi_0, H \Phi_2^4  \rangle\\
&\quad +(t_1^3t_2^4-t_2^3t_1^4 )t_1^3 \langle \Phi_0, H \Phi_1^3  \rangle
+(t_1^3t_2^4-t_2^3t_1^4 )t_1^4 \langle \Phi_0, H \Phi_1^4  \rangle \\
&\quad +(t_1^3t_2^4-t_2^3t_1^4 )t_{12}^{34} \langle \Phi_0, H \Phi_{12}^{34}  \rangle
\end{aligned}
\end{equation}
which yields
\begin{equation}
\begin{aligned}
(2,0,0,1,0),
(1,1,0,1,0),
(1,0,1,1,0),
(1,0,0,2,0),
(1,0,0,1,1),\\
(1,1,1,0,0),
(0,2,1,0,0),
(0,1,2,0,0),
(0,1,1,1,0)
(0,1,1,0,1)
\end{aligned}
\end{equation}
and
\begin{equation}
\begin{aligned}
(t_1^3t_2^4-t_2^3t_1^4)\langle \Phi_{12}^{34}, T H   \Phi_{12}^{34}  \rangle
&=
(t_1^3t_2^4-t_2^3t_1^4 )t_{12}^{34} \langle \Phi_0, H \Phi_{12}^{34}  \rangle
+(t_1^3t_2^4-t_2^3t_1^4 )t_1^3 \langle \Phi_{2}^{4}, H  \Phi_{12}^{34}  \rangle\\
&\quad+(t_1^3t_2^4-t_2^3t_1^4 )t_2^4 \langle \Phi_{1}^{3}, H  \Phi_{12}^{34}  \rangle
-(t_1^3t_2^4-t_2^3t_1^4 )t_1^4 \langle \Phi_{2}^{3}, H  \Phi_{12}^{34}  \rangle\\
&\quad
-(t_1^3t_2^4-t_2^3t_1^4 )t_2^3 \langle \Phi_{1}^{4}, H  \Phi_{12}^{34} \rangle
\end{aligned}
\end{equation}

which corresponds to:
\begin{equation}
\begin{aligned}
(2,0,0,1,0),
(1,1,0,1,0),
(1,0,1,1,0),
(1,0,0,2,0),
(1,0,0,1,1),\\
(1,1,1,0,0),
(0,2,1,0,0),
(0,1,2,0,0),
(0,1,1,1,0),
(0,1,1,0,1)
\end{aligned}
\end{equation}

The quadruply nested commutator term yields
\begin{equation}
\begin{aligned}
\langle \Phi_{12}^{34}, [[[[H,T],T],T],T] \Phi_0  \rangle 
&=
6 \langle \Phi_{12}^{34}, T^2HT^2 \Phi_0 \rangle\\
&=
6 (t_1^3t_2^4- t_2^3t_1^4 )\left( 
t_2^3 \langle \Phi_1^3, H \Phi_2^3  \rangle
+t_2^4 \langle \Phi_1^3, H \Phi_2^4  \rangle\right .\\
&\qquad\left.
+t_1^3 \langle \Phi_1^3, H \Phi_1^3  \rangle  +t_1^4 \langle \Phi_1^3, H \Phi_1^4  \rangle
+t_{12}^{34} \langle \Phi_1^3, H \Phi_{12}^{34}  \rangle\right)
\end{aligned}
\end{equation}
which corresponds to
\begin{equation}
\begin{aligned}
(2,0,0,1,0),(1,1,0,1,0),(1,0,1,1,0),(1,0,0,2,0),(1,0,0,1,1),\\
(1,1,1,0,0),(0,2,1,0,0),(0,1,2,0,0),(0,1,1,1,0),(0,1,1,0,1)
\end{aligned}
\end{equation}

Overall the Newton polytope that corresponds to the first polynomial is given by
\begin{equation}
\begin{aligned}
{\rm New}_{12}^{34}
={\rm conv } \big(&
(0,0,0,0,0),~(2,0,0,0,0),~(0,2,0,0,0),~(0,0,2,0,0),\\
&(0,0,0,2,0),~(0,0,0,0,2),~
(1,1,0,0,0),~(1,0,1,0,0),\\
&(1,0,0,1,0),~(1,0,0,0,1),~(0,1,1,0,0),~(0,1,0,1,0),\\
&(0,1,0,0,1),~(0,0,1,1,0),~(0,0,1,0,1),~(0,0,0,1,1),\\
&(1,1,1,0,0),~(1,1,0,1,0),~
(0,1,1,1,0),~(0,1,1,0,1),\\
&(1,0,1,1,0),~(1,0,0,1,1),~
(2,0,0,1,0),~(1,0,0,2,0),\\
&(0,1,2,0,0),~(0,2,1,0,0)
\big)
\end{aligned}
\end{equation}

\begin{equation}
\begin{aligned}
f_5(t)
&=
\langle \Phi_{12}^{34}, e^{-T}H e^T \Phi_0  \rangle\\
&=
C
+t_2^3 \langle \Phi_{12}^{34}, H \Phi_2^3  \rangle
+t_2^4 \langle \Phi_{12}^{34}, H \Phi_2^4  \rangle
+t_1^3 \langle \Phi_{12}^{34}, H \Phi_1^3  \rangle
+t_1^4 \langle \Phi_{12}^{34}, H \Phi_1^4  \rangle\\
&\quad
+t_{12}^{34} \langle \Phi_{12}^{34}, H \Phi_{12}^{34} 
-t_{12}^{34} \langle \Phi_0, H \Phi_0  \rangle
-t_1^3 \langle \Phi_{2}^{4}, H \Phi_0  \rangle
-t_2^4 \langle \Phi_{1}^{3}, H \Phi_0  \rangle\\
&\quad 
+t_1^4 \langle \Phi_{2}^{3}, H \Phi_0  \rangle
+t_2^3 \langle \Phi_{1}^{4}, H \Phi_0  \rangle
+\big(t_1^3t_2^4 - t_2^3t_1^4\big)  \langle \Phi_{12}^{34}, H \Phi_{12}^{34}  \rangle \\
&\quad 
-\big ( 
t_{12}^{34} \langle \Phi_0, H T \Phi_0  \rangle
+t_1^3 \langle \Phi_{2}^{4}, H T \Phi_0  \rangle
+t_2^4 \langle \Phi_{1}^{3}, H T \Phi_0  \rangle\\
&\qquad-t_1^4 \langle \Phi_{2}^{3}, H T \Phi_0  \rangle
-t_2^3 \langle \Phi_{1}^{4}, H T \Phi_0  \rangle
\big)+ \big(t_2^3t_1^4- t_1^3t_2^4\big) \langle \Phi_0, H \Phi_0  \rangle\\
&\quad
- (t_1^3t_2^4-t_2^3t_1^4 )t_{12}^{34} \langle \Phi_0, H \Phi_{12}^{34}  \rangle
+(t_1^3t_2^4-t_2^3t_1^4 )t_1^3 \langle \Phi_{2}^{4}, H  \Phi_{12}^{34}  \rangle\\
&\quad+(t_1^3t_2^4-t_2^3t_1^4 )t_2^4 \langle \Phi_{1}^{3}, H  \Phi_{12}^{34}  \rangle
-(t_1^3t_2^4-t_2^3t_1^4 )t_1^4 \langle \Phi_{2}^{3}, H  \Phi_{12}^{34}  \rangle\\
&\quad
-(t_1^3t_2^4-t_2^3t_1^4 )t_2^3 \langle \Phi_{1}^{4}, H  \Phi_{12}^{34} \rangle\\
&\quad
-(t_1^3t_2^4- t_2^3t_1^4 ) \langle \Phi_0, HT \Phi_0\rangle\\
&\quad
+(t_1^3t_2^4- t_2^3t_1^4 )\left( 
t_2^3 \langle \Phi_1^3, H \Phi_2^3  \rangle
+t_2^4 \langle \Phi_1^3, H \Phi_2^4  \rangle\right .\\
&\qquad\left.
+t_1^3 \langle \Phi_1^3, H \Phi_1^3  \rangle  +t_1^4 \langle \Phi_1^3, H \Phi_1^4  \rangle
+t_{12}^{34} \langle \Phi_1^3, H \Phi_{12}^{34}  \rangle\right)
\end{aligned}
\end{equation}

\newpage
In the following we shall list the coefficient matrices for indices $i\geq 2$ corresponding to the quadratic problem describing the CC equations for two electrons in four spin orbitals. 

\begin{align*}
    H(2) = \begin{pmatrix}
    h_{2,0}&h_{2,1}&h_{2,2}&h_{2,3}&h_{2,4}&h_{2,5}&-h_{2,5}\\
     0&0&0&0&0&0&0\\
    -h_{0,0}&-h_{0,1}&-h_{0,2}&-h_{0,3}&-h_{0,4}&-h_{0,5}&h_{0,5}\\
    0&0&0&0&0&0&0\\
    0&0&0&0&0&0&0\\
    0&0&0&0&0&0&0\\
    0&0&0&0&0&0&0
    \end{pmatrix}
\end{align*}\begin{align*}
    H(3) = \begin{pmatrix}
    h_{3,0}&h_{3,1}&h_{3,2}&h_{3,3}&h_{3,4}&h_{3,5}&-h_{3,5}\\
    0&0&0&0&0&0&0\\
    0&0&0&0&0&0&0\\
    -h_{0,0}&-h_{0,1}&-h_{0,2}&-h_{0,3}&-h_{0,4}&-h_{0,5}&h_{0,5}\\
    0&0&0&0&0&0&0\\
    0&0&0&0&0&0&0\\
    0&0&0&0&0&0&0
    \end{pmatrix}
\end{align*}\begin{align*}
    H(4) = \begin{pmatrix}
    h_{4,0}&h_{4,1}&h_{4,2}&h_{4,3}&h_{4,4}&h_{4,5}&-h_{4,5}\\
    0&0&0&0&0&0&0\\
    0&0&0&0&0&0&0\\   
    0&0&0&0&0&0&0\\
    -h_{0,0}&-h_{0,1}&-h_{0,2}&-h_{0,3}&-h_{0,4}&-h_{0,5}&h_{0,5}\\
    0&0&0&0&0&0&0\\
    0&0&0&0&0&0&0
    \end{pmatrix}
\end{align*}\begin{align*}
    H(5) = \begin{pmatrix}
    h_{5,0}&h_{5,1}&h_{5,2}&h_{5,3}&h_{5,4}&h_{5,5}&-h_{5,5}\\
    h_{4,0}&h_{4,1}&h_{4,2}&h_{4,3}&h_{4,4}&h_{4,5}&h_{0,1}\\
    -h_{3,0}&-h_{3,1}&-h_{3,2}&-h_{3,3}&-h_{3,4}&-h_{3,5}&h_{0,2}\\
    -h_{2,0}&h_{2,1}&-h_{2,2}&-h_{2,3}&-h_{2,4}&-h_{2,5}&h_{0,3}\\
    h_{1,0}&h_{1,1}&h_{1,2}&h_{1,3}&h_{1,4}&h_{1,5}&h_{0,4}\\
    -h_{0,0}&-h_{0,1}&-h_{0,2}&-h_{0,3}&-h_{0,4}&-h_{0,5}&h_{0,5}\\
    +h_{0,0}&-h_{4,5}&h_{3,5}&h_{2,5}&-h_{1,5}&h_{0,5}&h_{0,5}
    \end{pmatrix}
\end{align*}

\end{document}